\documentclass[a4paper]{article}
\usepackage[english]{babel}
\usepackage[pdftex]{color,graphicx,hyperref}
\usepackage{amsmath, amssymb}
\usepackage{amsthm}
\usepackage{fancyhdr}
\usepackage[normalem]{ulem} % for strikethrough lines in normal text
\usepackage{cancel} % for math mode

\usepackage[left=3.5cm,right=2.5cm,top=2.5cm,bottom=2.5cm,includeheadfoot]{geometry}

\newcommand{\pib}{\pi^{(\beta)} } 
 
\newcommand{\tk}[2]{t^{(#1)}_{#2}}

\newcommand{\Y}[2]{Y^{(#1)}_{#2}}
\newcommand{\Z}[2]{Z^{(#1)}_{#2}}
\newcommand{\bY}[2]{\bar Y^{(#1)}_{#2}}
\newcommand{\bZ}[2]{\bar Z^{(#1)}_{#2}}

\def \qed{\hfill$\Box$}

\newcommand{\e}{\mathbb{E}}
\newcommand{\E}{\mathbb{E}}
\newcommand{\mal}{\mathbb{D}^{1,2} }

\newcommand{\dom}{\text{dom}(\delta)}

\renewcommand{\L}{\mathbf{L}}
\renewcommand{\P}{\mathbb{P}}
\newcommand{\R}{\mathbb{R}}

\newcommand{\cA}{\mathcal{A}}
\newcommand{\cB}{\mathcal{B}}
\newcommand{\cC}{\mathcal{C}}
\newcommand{\cE}{\mathcal{E}}
\newcommand{\cF}{\mathcal{F}}

\newcommand{\cH}{\mathcal{H}}

\newcommand{\cR}{\mathcal{R}}

\newcommand{\cX}{\mathcal{X}}

\newcommand{\cL}{\mathcal{L}}

\newcommand{\cP}{\mathcal{P}}

\newcommand{\cS}{\mathcal{S}}
\newcommand{\cT}{\mathcal{T}}
\newcommand{\cY}{\mathcal{Y}}
\newcommand{\cZ}{\mathcal{Z}}

\newcommand{\hg}{\mathfrak{g}}
\newcommand{\hH}{\mathfrak{H}}

\newcommand{\1}{\mathbf{1}}

    \def\independenT#1#2{\mathrel{\setbox0\hbox{$#1#2$}%
    \copy0\kern-\wd0\mkern4mu\box0}} 

\newcommand{\HFd}{{$\bf (A_{\partial f})$}}
\newcommand{\HFt}{{$\bf (A_{ f_t})$}}

\newcommand \cuM{C_{u}}
\newcommand \cwM{{\cal C}_{(\ref{eq:iteration:2}a)}}
\newcommand \cuhM{{\cal C}_{(\ref{eq:iteration:2}b)}}
\newcommand \cg[1]{{\cal C}^{(#1)} }

\newcommand{\thetaL}{\theta_L}
\newcommand{\thetaC}{\theta_c}
\newcommand{\thetaP}{\theta_\Phi}

\newcommand{\feps}{f^{(\varepsilon)}}
\newcommand{\Yeps}{Y^{(\varepsilon)}}
\newcommand{\Zeps}{Z^{(\varepsilon)}}
\newcommand{\bZeps}{\tilde Z^{(\varepsilon)}}
\newcommand{\yeps}{y^{(\varepsilon)}}
\newcommand{\zeps}{z^{(\varepsilon)}}
\newcommand{\bzeps}{\tilde z^{(\varepsilon)}}
\newcommand{\nabX}[1]{\nabla X_{#1}}
\newcommand{\nabXM}[3]{\nabla X^{(#1,#2)}_{#3}}

\newcommand{\inXM}[3]{\nabla X^{(#1,#2,-1)}_{#3}}
\newcommand{\normexp}[1]{\|#1\|_2}
\newcommand{\normS}[2]{\|#2\|_{\cS^{#1} }}
\newcommand{\aeps}[1]{a^{(\varepsilon)}_{#1} }
\newcommand{\beps}[1]{b^{(\varepsilon)}_{#1} }
\newcommand{\ceps}[1]{c^{(\varepsilon)}_{#1} }
\newcommand{\Ueps}{U^{(\varepsilon)} }
\newcommand{\Veps}{V^{(\varepsilon)} }
\newcommand{\YM}[2]{Y^{(\varepsilon,#1,#2)}} %clash with paper 3
\newcommand{\yM}[2]{y^{(\varepsilon,#1,#2)}}
\newcommand{\ZM}[2]{Z^{(\varepsilon,#1,#2)}}
\newcommand{\zM}[2]{z^{(\varepsilon,#1,#2)}}
\newcommand{\XM}[2]{X^{(#1,#2)}} %Clash with paper 3
\newcommand{\MW}[2]{H^{(#1,#2)}}

\numberwithin{equation}{section}
\newtheorem{definition}{Definition}[section]
\newtheorem{theorem}[definition]{Theorem}
\newtheorem{proposition}[definition]{Proposition}
\newtheorem{lemma}[definition]{Lemma}
\newtheorem*{lemma*}{Lemma}
\newtheorem{corollary}[definition]{Corollary}
\theoremstyle{definition} \newtheorem{remark}[definition]{Remark}
\theoremstyle{definition} \newtheorem*{remark*}{Remark}

\newcommand{\elip}{\bar\beta}

\newtheorem{plaremark}{Remark (Plamen)}

\newtheorem{emremark}{Remark (Emmanuel)}

\newcommand {\thetac}{{\theta_c}}
\newcommand {\thetal}{{\theta_L}}

\newcommand{\Hs}{$\bf (A_{b,\sigma})$}
\newcommand{\He}{$\bf (A_{\text{u.e.} })$}
\newcommand{\Hg}{$\bf (A_{\Phi })$}
\newcommand{\Hge}{$\bf (A_{b\Phi })$}
\newcommand{\Hgh}{$\bf (A_{h\Phi })$}
\newcommand{\Hgexp}{$\bf (A_{exp\Phi })$}

\newcommand{\Hf}{$\bf (A_{f })$}

\newcommand{\X}[2]{X^{ (#1)}_{ #2}}
\makeindex

\title{Two algorithms for the discrete time approximation of Markovian backward stochastic differential equations under local conditions }
\author{Plamen Turkedjiev\footnote{Currently affiliated to Laboratoire de Finance et des March\'es de l'Energie (FiME).
A significant part of the author's research was done at Humboldt University.
Email: \href{mailto:turkedjiev@cmap.polytechnique.fr}{\tt turkedjiev@cmap.polytechnique.fr}.
The author would like to thank in particular Dirk Becherer and Emmanuel Gobet for their feedback and spotting some early errors in the manuscript.
}\\
Centre de Math\'ematiques Appliqu\'ees\\
Ecole Polytechnique and CNRS\\
Route de Saclay\\
F 91128 Palaiseau cedex, France}
\date{\today}

\begin{document}

\maketitle

\begin{description}
\item{\sc Abstract:} Two discretizations of a  class of  {locally Lipschitz } Markovian backward stochastic differential equations (BSDEs) are studied. 
The first is the classical Euler scheme which approximates a projection of the processes $Z$, and the second a novel scheme based on Malliavin weights which approximates the marginals of the process $Z$ directly.
Extending the representation theorem of Ma and Zhang \cite{ma:zhan:02} leads to advanced a priori estimates and stability results for this class of BSDEs.
These estimates are then used to obtain competitive convergence rates for both schemes with respect to the number of points in the time-grid.
The class of BSDEs considered includes  Lipschitz BSDEs with fractionally smooth terminal condition,  {thus extending the results of \cite{gobe:makh:10}, }  quadratic BSDEs with bounded, H\"older continuous terminal condition  {(for bounded, differentiable volatility), and BSDEs related to proxy methods in numerical analysis}.

%\vspace{0.3cm}
\item{\sc Keywords:}
Backward stochastic differential equation, approximation schemes, Malliavin calculus, representation theorem, a priori estimates.

\item {\sc MSC 2010:} 60H35, 65C30, 60H07, 60H10.

\end{description}

\section{Introduction}
\label{section:intro}

%BSDEs have a myriad of applications,  notably in the theory of mathematical finance, stochastic control, and partial differential 
%equations. Although the theory and applications of BSDEs has grown rapidly over the last twenty years, the development of efficient
%numerical methods is a more recent topic of research. 

{$\blacktriangleright$ \bf Framework.}
Backward stochastic differential equations play an important role in the theory of mathematical finance, stochastic optimal control, and partial differential equations.
In this paper, we study two discrete-time approximations of the for the  {so-called {\em locally  Lipchitz}} Markovian backward stochastic differential equation (BSDE).
 {The purpose is to determine the error induced by these approximations under suitable norms.}
The first is the well-established Euler scheme for BSDEs, and the second is a novel scheme we call the Malliavin weights scheme for BSDEs.
Let $T>0$ be a fixed terminal time and $(\Omega,\cF_{T},\{\cF_t\},\P)$ a filtered probability space, 
where $\{\cF_{t}:0\le t\le T\}$  is the filtration  generated by a $q$-dimensional ($q\ge1$) Brownian motion $W$
and satisfying the usual conditions of right-continuity and completeness.
We  look to approximate the $\R\times(\R^q)^\top$-valued, predictable process $(Y,Z)$ solving the BSDE
\begin{equation}
{Y_t = \Phi(X_{T}) + \int_t^T f(s,X_s,Y_s,Z_s) ds  - \int_t^T Z_s dW_s.}
\label{eq:1:BSDE}
\end{equation}
Here, $(\R^q)^\top$ is the space of $q$-dimensional, real valued row vectors;
$X$ is an $\R^d$-valued ($1 \le d \le q$) diffusion; and  
$\Phi:\R^d\rightarrow\R$ and $f:[0,T)\times\R^d\times\R\times(\R^q)^\top\rightarrow\R$ are deterministic functions that are  termed the {\em terminal condition} and {\em driver}, respectively. 
We focus on the setting in which the terminal condition $\Phi$ is in the space of {\bf fractionally smooth} functions $\L_{2,\alpha}$ for parameter $\alpha \in (0,1]$ - see \Hg \ in Section \ref{sectio:1:ass} for details - 
and the driver is
{\bf locally Lipschitz continuous in $(x,y,z)$} and {\bf locally bounded at $0$} in the sense that
%To be precise, let $f : [0,T] \times \R^d \times \R \times (\R^q)^\top \rightarrow \R$ be continuous and
%assume 
there exist  { exponents $\thetaL,\ \theta_X, \ \thetaC \in (0,1]$, finite constants $L_f,L_X, C_f\ge0$, 
%and a function $L_X : [0,T) \rightarrow [0,\infty)$ dominated by $C_X (T-t)^{(\thetaL - 1)/2}$ 
 such that,
for all $t \in [0,T)$  and $(x,y,z), (x',y',z') \in \R^d \times \R \times (\R^q)^\top$,
\begin{align}
%\left.
%\begin{array}{rl}
|f(t,x,y,z) - f(t,x',y',z')| & \le L_f\frac{  |y-y'| + |z-z'|}{(T-t)^{(1-\thetaL)/2} } + L_X{|x-x'| \over (T-t)^{1 - \theta_X/2}}, \nonumber\\
%\quad
 |f(t,x,0,0)| & \le \frac{C_f}{(T-t)^{1-\thetaC}}.
%\end{array}
%\right\}
\label{eq:loclip:driver}
\end{align}
}
%We say that $(Y,Z)$ solves a BSDE {\bf under local conditions}. 
 {Furthermore, $X$ solves a time-inhomogeneous stochastic differential equation (SDE) with suitable coefficients; see \Hs \ in Section \ref{sectio:1:ass}. }
The existence and uniqueness of this  class of BSDEs -- given in Section \ref{section:mal:BSDE} --   {follows from \cite[Theorem 3.2]{fan:jian:12}}.
Below, we show that this class of BSDEs includes a  section of the important quadratic BSDEs, and  also BSDEs related to so-called proxy schemes used for numerical methods,
so it is of  interest to find good discrete-time approximations for such BSDEs.
We note that fully implementable algorithms -- admitting the full generality of the assumptions considered in this paper -- based on the Euler  and Malliavin weights schemes have been studied in detail in \cite{gobe:turk:13a}\cite{gobe:turk:13b} respectively,
but, to the best of our knowledge, this is the first paper  considering the discretization error under the full generality of the local conditions.

{\noindent $\blacktriangleright$ \bf Summary of results.}
%In this paper, we present two discrete-time approximation schemes for the BSDEs of the form given above. 
In the spirit of \cite{gobe:makh:10}, we make use of non-uniform time-grids $\big \{  \pib_ N := \{ 0 = \tk N0 < \ldots < \tk NN = T\} \ : \ N \ge 1   \big \}$ whose parameter $\beta \in (0,1]$ determines the time-points 
$
\tk {N} i := T - T(1 - i/N) ^{1/\beta}
$.
 {
As in \cite{gobe:makh:10}, the use of these time-grids appears to substantially reduce the error due to disctretization.
}

The first  approximation, studied in Section \ref{section:reg}, is the so-called Euler scheme for BSDEs:
% for each $ N$, set
\begin{align}
\Y NN & := \Phi(X_T), \quad \Z Ni := \frac{1}{\tk N{i+1} - \tk Ni} \E[\Y N{i+1}(W_{\tk N{i+1}}-W_{\tk Ni})^\top|\cF_{\tk Ni}], \nonumber\\
\Y Ni & := \E[\Y N{i+1}+ f(\tk Ni,X_{\tk Ni}, \Y N{i+1}, \Z Ni)(\tk N{i+1} - \tk Ni)|\cF_{t_i}]
\label{eq:1:dp}
\end{align}
for each $i\in\{0,\ldots,N-1\}$. 
The random variable $\Z Ni$ is a discretization of the projection $(\tk N{i+1} - \tk Ni)\tilde Z_{t_i} := \E[ \int_{\tk Ni}^{\tk N{i+1}} Z_s ds | \cF_{\tk Ni}]$.
This approximation has been frequently studied: \cite{zhan:04}\cite{bouc:touz:04}\cite{gobe:laba:07} among others, in the setting where the terminal condition $\Phi$ and the driver are uniformly Lipschitz continuous (i.e. $\thetaL = 1$);
 \cite{gobe:makh:10} in the setting of the fractionally smooth  $\Phi$ but uniformly Lipschitz continuous driver;
 \cite{imke:dosr:10}\cite{rich:11} in the setting of bounded Lipschitz (resp.  H\"older) continuous $\Phi$ and quadratic driver;
 and \cite{rich:12} in the setting of possibly unbounded (locally) Lipschitz continuous $\Phi$ and (super-)quadratic driver.  
 Typically, the discretization error of the Euler scheme is measured by 
\begin{equation}
\cE(N) := \max_{0\le i <N} \E[|Y_{\tk Ni} - \Y Ni|^2] + \sum_{i=0}^{N-1} \int_{\tk Ni}^{\tk N{i+1}} \E[|Z_{t} - \Z Ni|^2] dt.
\label{eq:disc:error:gen}
\end{equation}
%Let us denote the discretization error \eqref{eq:disc:error:gen} associated to the Euler scheme by $\cE_1(N)$.
%By first of all applying the techniques of \cite{gobe:makh:10} directly, 
We show  in Theorem \ref{prop:1:l2 reg:makh} that
 if $\beta < (2\gamma)\wedge\alpha$, where $\gamma := (\frac{\alpha}{2}\wedge \thetaC + \frac{\thetaL}{2})\wedge\thetaC$, then
\[
 \cE(N) \le CN^{-1} \1_{[1,2]}(\alpha + \thetaL) + CN^{-2\gamma}\1_{(0,1)}(\alpha + \thetaL).
\]
%The constant $C$ does not depend on $N$.
The optimal error bound $O(N^{-1})$ is obtained if $\alpha + \thetaL \ge 1$.
This rate is optimal in the sense that it is the same as the rate of convergence obtained in \cite[Theorem 3.2]{gobe:makh:10} in the uniformly Lipschitz driver setting ($\thetaL =1$).
%In the special case where $\thetaL = \alpha$, which is the case of interest for the quadratic BSDE,
%this implies that the optimal rate $N^{-1}$ is obtained for $\alpha \ge 1/2$.
This result can be complimented under the additional assumption that the terminal condition $\Phi$ is $\thetaP$-H\"older continuous:
in Theorem \ref{thm:L2 reg}, we  show that
if $\beta < (2\gamma)\wedge\alpha\wedge\thetaL$, then
\[
\cE(N) \le CN^{-1} \1_{[1,4]}(\thetaP + \beta + 2\gamma) + CN^{-2\gamma} \1_{(0,1)}(\thetaP + \beta + 2\gamma).
\]
%where $\delta_N \ge \1_{[3,\infty)}(N) \ln\ln(N) / \ln(N)$.  
Now $\thetaP + \beta + 2\gamma \ge 1$ is sufficient to obtain the optimal convergence rate 
$O(N^{-1})$.
Although the complex relationship between $\thetaP$, $\alpha$ and $\gamma$ make it difficult to compare the two results in full generality, 
 {  the latter result relaxes the constraint $\alpha +\thetaL \ge 1$ in order to obtain the optimal  error bound $O(N^{-1})$ if $\thetaC \ge 1/2$ -- see \eqref{eq:loclip:driver} to recall the definition of $\thetaC$. }

The second approximation, studied in Section \ref{section:num}, is the so-called Malliavin weights scheme. 
Rather than approximating the projections of the process $Z$, this algorithm is used to approximate the version of $Z$,  {determined by the Malliavin integration-by-parts formula of Theorem \ref{thm:repr}}, at the points of the time grid directly:
for each $N\ge 1$, set
\begin{align}
 \bY NN & := \Phi(X_T), \quad \bY N{i} :=  \e[\Phi(X_{T}) + \sum_{j=i}^{N-1}f(\tk Nj, X_{\tk Nj},\bY N{j+1}, \bZ N{j})(\tk N{j+1} - \tk Nj)|\cF_{\tk Ni}], \nonumber \\
 \bZ N{i} & := \e[\Phi(X_{T})H^{i}_{N} + \sum_{j=i+1}^{N-1} f(\tk Nj , X_{\tk Nj},\bY N{j+1}, \bZ N{j})H^{i}_{j}(\tk N{j+1} - \tk Nj)|\cF_{\tk Ni}]
\label{eq:1:disc:malscheme}
\end{align}
for $i\in\{0,\ldots,N-1\}$, where 
$(H^i_j)_{i,j}$ is a suitable random variable.
Due to the connection between BSDEs and quasilinear partial differential equations (PDEs) -- see \cite{rich:12}\cite{cris:dela:12} and references therein -- it may be of interest to approximate the marginals of the process $Z$ rather than the projections.
Other schemes that make use of Malliavin calculus  are available \cite{bria:laba:13}\cite{hu:nual:song:11}, but this is, to the best of our knowledge, the first scheme  {which makes use of the Malliavin integration-by-parts formula (Theorem \ref{thm:repr}}).
Convergence results are given 
 {-- for weaker norms than those used in $\cE(N)$ for the Euler scheme -- }
 in Theorem \ref{thm:1:disc er}.
  {Although one is able to prove results under stronger norms than for the Euler scheme, there are several disadvantages (regardless of the norm used to measure the error) of the Malliavin weights scheme over the Euler scheme.}
Our results are proven under stronger conditions than for the Euler scheme
because the use of stronger a priori estimates -- Proposition \ref{prop:veps:bd} --  is essential in the proof:
one requires that either the terminal condition has exponential moments or that it is H\"older continuous.
%, and unfortunately
 We have not yet been able to  weaken the conditions on these a priori estimates.
 One also requires a greater constraint $\beta \le \gamma \wedge \thetaL \wedge \alpha$  { (where $\gamma := (\frac{\alpha}{2}\wedge \thetaC + \frac{\thetaL}{2})\wedge\thetaC$)} on the time-grid than for the Euler scheme.
 The rate of convergence again depends on the parameters $(\alpha,\thetaL,\thetaC,\beta)$.
 In the more general setting  {of exponential moments on the terminal condition}, $\beta + 2\gamma \ge 1$ is required for the  {optimal error bounds $O(N^{-1})$}, whereas in the setting of $\thetaP$- H\"older continuous terminal condition, 
 $\beta + \thetaP + 2\gamma \ge 1$ is sufficient.
 { One may ask, given the additional constraints, why it is of interest to study the Malliavin weights scheme over the Euler scheme?
The reason has to do with the approximation of the conditional expectation.
}
It is shown in \cite{gobe:turk:13b} that, using Monte Carlo  least-squares regression to approximate the conditional expectation, one can theoretically gain an order one improvement  {with respect to the number of time-steps $N$} on the algorithm complexity using the Malliavin weights scheme compared to the multi-step forward implementation of the Euler scheme \cite{gobe:turk:13a}.
 {Such a complexity reduction is substantial, given that $N$ may be very large.}

 {In order to obtain the results on discretization, we extend some basic tools from the literature of BSDEs. 
These results are interesting in their own right.
}
 {
Firstly, we extend stability estimates for Lipschitz BSDEs to the class of BSDEs satisfying local Lipschitz continuity and boundedness conditions \eqref{eq:loclip:driver}. 
This enables us to make estimates on the basis of constructing approximating sequences, a key technique used throughout the paper.
A natural consequence of stability estimates are a priori estimates, which we also frequently require.
These results are contained in Section \ref{section:regularity}.
}
 {
Secondly, we  obtain  dynamical representations of the process $Z_t$ in the form of the product  $U_t \sigma(t,X_t)$, where $(U,V)$ is the solution of a linear BSDE.
Such representations are very valuable for making estimates on the increments $\E[|Z_t - Z_s|^2]$, because one can make use of a priori estimates on the linear BSDE and the process $X$.
In fact, it is not possible to obtain the results for $Z$ directly, but for a suitable sequence $\{\Zeps_t \ : \ \varepsilon >0\}$ of approximating BSDEs.
A priori estimates for the approximation are computed and play an important role in the overall convergence rate of the numerical schemes.
To obtain this result, we extend the method and results of \cite[Section 2]{gobe:makh:10}, who consider the setting \eqref{eq:loclip:driver} with $\thetaL =\thetaC = 1$ only, to our more general setting.
The key results are contained in Lemma \ref{lem:mal:ueps veps}.
%In order to obtain a priori estimates on the linear BSDEs, 
%we must take more care compared to \cite{gobe:makh:10} due to the local Lipschitz continuity and boundedness of the driver (i.e. the case $\thetaC <1$ and $\thetaL<1$).
%, and for this we extend the a priori results of  \cite[Proposition 3.2]{bria:dely:hu:pard:stoi:03}.
}
 {Thirdly, we extend the classical representation theorem of Ma and Zhang \cite[Theorem 4.2]{ma:zhan:02} for the $Z$ process to our class of BSDEs.
This theorem is proved in  Section \ref{section:repr}
and is a key result in this paper.
One the one hand, it is the basis for the Malliavin weights scheme.
}
On the other hand, we use the representation theorem to obtain stability estimates directly on the marginals of the process $Z$ -- see Proposition \ref{prop:1:apriori} -- which are key to the analysis.
These stability estimates lead in turn to a priori estimates of the form 
\begin{align*}
|Z_t|  \le C \sqrt{(T-t)^{-1}\E_t[|\Phi(X_T) - \e_t[\Phi(X_T)]|^2]} + C(T-t)^{\thetaC - 1/2 } + C \E_t[\Phi(X_T)^2]^{1/2}(T-t)^{\thetaL/2} 
%\\
%\E_s[|\Delta Z_t|^2]^{1/2}  \le C\frac{\E_s\big[(\Delta\Phi - \E_t\Delta\Phi)^2]^{1/2} }{\sqrt{T-t}} 
%+ C \int_t^T \frac{\E_s[\Delta f_r^2]^{1/2} }{\sqrt{r-t}} dr 
% + C \E_s[\Delta\Phi^2]^{1/2}(T-t)^{\thetaL/2} 
\end{align*}
for all $t\in[0,T)$ almost surely.
%enabling estimates on the norm and moments of $Z_t$ depending on the regularity of the terminal condition; 
%see Proposition \ref{cor:mom bd}.
Such estimates are,  {to the best of our knowledge,} novel and allow us to study the impact of the regularity of the terminal condition on a priori estimates -- see Proposition \ref{cor:mom bd}. 
% {Fourthly,
%The fractionally smooth terminal condition causes a host of technical problems due to time dependancy, especially when combined with a complex driver, and this paper extends the tool-box for dealing with these problems;
%see also \cite{geis:geis:04}\cite{geis:hujo:07} \cite{gobe:makh:10}\cite{geis:gobe:11}\cite{geis:geis:gobe:12} for fractionally smooth terminal conditions. 
%}
 {
Finally, in Proposition \ref{prop:veps:bd}, we obtain a priori estimates for the process $\Veps_t$ --
the solution $(\Ueps,\Veps)$ to  the linear BSDE such that the approximating BSDE solution satisfies $\Zeps_t = \Ueps_t \sigma(t,X_t)$ -- under additional regularity conditions on the terminal condition.
These estimates are essential to analyse the error due to the Malliavin weight scheme.
Rather than considering a second Malliavin derivative of the process $Y_t$, as for example do \cite{cris:dela:12}, we make use of a functional representation that comes from the Markov property of $X$ and determine regularity properties of the said functional representation.
A  consequence of this is the Lipschitz continuity of the functional representation of the process $Z_t$ under suitable conditions -- see Corollary \ref{cor:proxy:lip}.
To our knowledge, this result is novel.
Since regularity properties are very useful for the calibration of numerical schemes -- see for example \cite[Section 4.4]{gobe:turk:13a} -- this result may have some impact on reducing the cost of fully implementable algorithms.
}

{\noindent $\blacktriangleright$ \bf Contributions to quadratic BSDEs and proxy methods.}
We consider the setting where $\Phi$ is a bounded, $\thetaP$-H\"older continuous function.
To make the contributions of the numerical results in this paper clearer, we consider two important examples.
Note that these examples have also been given some attention in  \cite[Section 2]{gobe:turk:13a}.
 { We emphasize that the forward process $X$ is  a diffusion with bounded, twice continuously differentiable coefficients, whose partial derivatives are bounded and H\"older continuous; this assumption stands throughout this paper -- see \Hs.}

%In the first application, we treat a large class of quadratic BSDEs.
Quadratic BSDEs have powerful applications in financial mathematics, for example to solve utility optimization problems in incomplete markets \cite{elka:roug:01}\cite{hu:imke:mull:05}.
Let $q=d$ and the measurable function $F:[0,T)\times\R^d\times\R\times\R^d\rightarrow \R$ satisfy
\begin{align*}
|F(t,x,y,z)|&\leq c\ (1+|y|+|z|^2), \\
|F(t,x,y,z)-F(t,x,y',z')|&\leq c\ (1+|z|+|z'|)(|y-y'|+|z-z'|).
\end{align*}
It is known \cite{dela:guat:06} that the solution $(Y,Z)$ of the BSDE with terminal condition $\Phi$ and driver $F(t,x,y,z)$ exists and is unique and that  there is a constant $\theta \in (0,1]$ and finite $C_u>0$ such that $|Z_t| \le C_u (T-t)^{(\theta-1)/2}$
for all $t\in[0,T)$ almost surely.
This implies that $(Y,Z)$ also solves the BSDE under local conditions with terminal condition $ \Phi$ and driver $ f(t,x,y,z) := F(t,x,y, \cT_{C_u(T-t)^{(\theta-1)/2}}(z))$,
where $\cT_L(z) := (-L \vee z_1\wedge L, \ldots , -L \vee z_q \wedge L)$.
Indeed, $C_{ f}=c$, $\thetac=1$, $L_{ f}=c(T^{(1-\theta)/2}+2\sqrt d C_u)$, and $\thetal= \theta$.
The terminal condition is fractionally smooth with parameter $\alpha$ at least as large as $\thetaP$ - see Remark \ref{rem:holder}.
It is shown in Corollary \ref{cor:mom bd} that $|Z_t| \le C(T-t)^{(\thetaP - 1)/2}$, so $\theta_L$ is at least as large as $\thetaP$.
Therefore, the error $\cE(N)$ of the Euler scheme is bounded above by $C_\beta N^{-1} \1_{[1,4]}(3\thetaP + \beta ) + CN^{-2\thetaP} \1_{(0,1)}(3\thetaP + \beta )$ for any $\beta < \thetaP$.
In \cite{rich:11}, the Euler scheme for bounded, H\"older continuous is also considered, 
but with a different non-uniform time-grid and a transformation of the terminal condition;
 {there is a further modelling difference in that the author requires no uniform elliptic condition, but sacrifices state-dependence in the volatility matrix}.
The author obtains a rate of convergence $C_\eta N^{\eta - \thetaP}$ for any $\eta >0$, so we have obtained an improvement in this work;
This improvement is likely due to the use of the time-grids $\pib_N$ in our scheme --
indeed, \cite{gobe:makh:10} show a rate of convergence $O(N^{-\alpha})$ in the uniformly Lipschitz continuous driver setting if only a uniform time-grid is used.
 {It is important to remark that this work is a complement to the recent papers \cite{rich:12}\cite{chas:rich:14}, in which the authors consider  weaker assumptions on the drift and the volatility of the SDE -- only Lipschitz continuity and linear growth are required -- however stronger assumptions are required on the terminal function $\Phi$, which must be locally Lipschitz continuous.}

Next we consider a particular instance of the proxy method. 
Let $F(t,x,y,z)$ satisfy \eqref{eq:loclip:driver} with exponents $\theta_{L,F} \le 1$, $\theta_{X,F} =1$ and $\theta_{c,F} = 1$,  and constants $L_{F}$, $L_{F,X}$ and $C_{F}$.
Let $(\cY,\cZ)$ satisfy the BSDE with terminal condition $\Phi$ and driver $F(t,x,y,z)$.
%For functions $b$ and $\sigma$ -- the drift and volatility of $X$ -- satisfying the conditions of \Hs \ in Section \ref{sectio:1:ass}, 
 {
Let the function $\bar F(t,x,y,z)$ satisfies \eqref{eq:loclip:driver} with exponents $\theta_{L,\bar F} = \theta_{X,\bar F} = \theta_{c,\bar F} = 1$,  and constants $L_{\bar F}$, $L_{\bar F,X}$ and $C_{\bar F}$, and $\bar \Phi(x)$ is $\theta_\Phi$-H\"older continuous
}
and suppose that the  parabolic PDE  
 {
\begin{equation*}
%\label{eq:lin:PDE}
%\left.
\begin{array}{l}
0 = \partial_t v +  {\bar \cL_{t,x} v+ \bar F(t,x, v(t,x), \nabla_x v(t,x) \sigma(t,x) )}, 
%\\ \\
\qquad
v(T,x) = \bar \Phi(x)
\end{array}
%\right\}
\end{equation*}
}
has a unique strong solution $v$, and, for every $t\in[0,T)$, the $k$-th order ($k\le3$) partial derivatives in $x$ of $v$  are bounded by $C_u (T-t)^{(\thetaP -k)/2}$.
 {We assume also that the parabolic operator $\bar \cL_{t,x}$ satisfies the property that, for any $i\in\{1,\ldots,d\}$, $\|\partial_{x_i} \{ \bar \cL_{t,x} -\cL_{t,x} \} v(t,\cdot) \|_{\infty} \le C_u (T-t)^{(\thetaP -2)/2}$, where 
$\cL_{t,x}$ is the parabolic operator given by
\[
\cL_{t,x} u(t,x) :=  \left\{\frac12 \sum_{i,j=1}^d (\sigma(t,x) \sigma(t,x)^\top)_{i,j} \frac{\partial^2}{\partial_{x_i}\partial_{x_j}}
+ \sum_{i=1}^d b_i (t,x)\frac{\partial}{\partial_{x_i}} \right\} u(t,x) ;
\]
this is stronger than  the previous assumption on the third order partial derivatives of $v(t,\cdot)$, which asks for the upper bound $C_u(T-t)^{(\thetaP-3)/2}$.
}
%It is well known that 
Then  $(Y_t,Z_t) := (\cY_t - v(t,X_t), \cZ_t - \nabla_xv(t,X_t)\sigma(t,X_t) )$ solves a BSDE 
with  terminal condition $\Phi(x) - \bar\Phi(x)$ and driver 
\[f(t,x,y,z) :=  F(t,x,v(t,x)+y,\nabla_x v(t,x)\sigma(t,x) + z)  {-\bar F(t,x,v(t,x),\nabla_x v(t,x)\sigma(t,x)) + (\cL_{t,x} - \bar \cL_{t,x}) v(t,x) }.\] 
 {The driver $f(t,x,y,z)$ satisfies \eqref{eq:loclip:driver} with exponents $\thetaL = \theta_{L,F}$,  {$\theta_X = \thetaL + \thetaP -1$,}  $\thetaC = \thetaP$, and constants 
$L_f = L_F$, $L_X := L_F C_u + \sqrt d C_uT^{(1-\thetaL)/2}(1+L_{\bar F}) $, and $C_f = \sqrt d (L_F + L_{\bar F}) C_u + C_F + C_{\bar F}$.}
The idea is that it may be numerically advantageous to simulate the BSDE  $(Y,Z)$ as opposed to the original BSDE $(\cY,\cZ)$.
A simple example of a  proxy is given by  {$\bar \Phi(x) \equiv \Phi(x)$, $\bar F \equiv 0$,} and
$
 {\bar \cL_{t,x} u(t,x) = \cL_{t,x} u(t,x);}
% :=  \left\{\frac12 \sum_{i,j=1}^d (\sigma(t,x) \sigma(t,x)^\top)_{i,j} \frac{\partial^2}{\partial_{x_i}\partial_{x_j}}
%+ \sum_{i=1}^d b_i (t,x)\frac{\partial}{\partial_{x_i}} \right\} u(t,x) ;}
$
see Lemma \ref{lem:lin:pde:bds} for the gradient bounds.
 {We show in Corollary \ref{cor:proxy:lip} that the process $(Y,Z)$ brought about by this proxy  may lead to some regularity improvements for the process $Z$ compared with the original process $\cZ$.
This may lead to an improvement of the numerical complexity for fully implementable algorithms that approximate the conditional expectation, where regularity is extremely important;
moreover, \cite{gobe:turk:13a}\cite{gobe:turk:13b} both demonstrate that there will an improvement in the constants for the error estimates when using Monte Carlo least-squares regression on this proxy compared to the same algorithm on the original BSDE $(\cY,\cZ)$.
}

\vspace{0.3cm}

%Throughout this paper, we actually compute the discretization error of $(Y,Z)$ in order to obtain that of $(\cY,\cZ)$. 
%This strategy is used for both Euler and Malliavin weights scheme --  see \eqref{} and \eqref{} -- although we do not explicitly state it.
%, although it is technically not in the scope of \eqref{eq:loclip:driver}.
%For $\thetaL \equiv1$, it is known that this error has a better rate of convergence \cite{gobe:makh:10}, but it is not necessarily so for $\thetaL < 1$.
%On the other hand, it is shown in the proof of Proposition \ref{prop:veps:bd} -- equation \eqref{eq:z0} -- 
%that $|Z_t| \le C (T-t)^{(\thetaP + (2\thetaC)\wedge\thetaL-1)/2}$, whereas $|\cZ_t| \le C (T-t)^{((2\thetaC)\wedge\thetaP -1)/2}$. 
%In fact, smaller almost sure bounds yield improved rates of convergence for numerical methods \cite{gobe:turk:13b}\cite{gobe:turk:13a}, so it is numerically beneficial to simulate $(\cY,\cZ)$ rather than $(Y,Z)$.

\noindent {$\blacktriangleright$ \bf  Remarks on extensions.}
In this paper, we work with one of the simplest time-inhomogeneous  SDE models with stochastic volatility,
which, in particular, allows us to make use of  results from the theory of parabolic PDEs \cite{frie:64}
-- see Lemma \ref{lem:lin:pde:bds}.
The representation theorem for $Z$ in Theorem \ref{thm:repr} also makes use of the uniform ellipticity condition. 
Our application to quadratic BSDEs  requires these conditions, and additionally that $\Phi$ is H\"older continuous and bounded, because we make use of the results of \cite{dela:guat:06} to introduce local Lipschitz continuity.
There are already several directions that may help us to avoid the uniformly elliptic condition.
The  results of \cite{kusu:03}\cite{cris:dela:12}\cite{nee:11}, offer suitable PDE results  under UFG conditions.
Also, a representation theorem beyond the uniformly elliptic setting has been found by \cite{zhan:05} and \cite{gobe:muno:05} (although only for the zero driver case in the second reference).
%It is the goal of this paper to extend this literature and provide solutions that may be stepping stones to further extension.
Another interesting aspect of our general results is that we require neither BMO results nor (local)-Lipschitz continuity of $\Phi$.
Combined with the connection to quadratic BSDEs already discussed here, this suggests 
%\xout{there may be potential to apply the results of this} 
 {the results of this paper may be an important stepping-stone to}  obtain
novel representation theorems, a priori estimates, existence and uniqueness results for (super-)quadratic BSDEs with possibly unbounded and discontinuous terminal conditions.
 {It would also be interesting to combine the results of this paper with those of \cite{rich:12} to handle the setting of unbounded, state-dependent $\sigma$ with non-Lipschitz continuous terminal condition.}
%but it is not clear that the Malliavin weights obtained are directly suitable for the analysis undertaken in this paper.
Unfortunately, all of these extensions are beyond the scope of this paper.

\subsection{Notation and conventions}
\label{section:notation}
$\blacktriangleright${\bf Time-grids.}
Since each result is given for a fixed number of time-points $N$, we denote the points $\{\tk Ni\}$ of the time-grid simply by $\{t_i\}$.
Let $\Delta_i := t_{i+1} - t_i$ and $\Delta W_i := W_{t_{i+1}} - W_{t_i}$.
We also suppress the superscript $(N)$ in the Euler and Malliavin weights scheme.

\noindent $\blacktriangleright$
{\bf Expectations and norms.}
For $p \ge1$, we denote by $\| \cdot\|_p$ the norms $(\E[| \cdot|^p])^{1/p}$;
in particular, we  make  use of the  norm $\sqrt{\e[| \cdot|^2]}$ denoted by $\normexp{\cdot}$.

% $n$, $k$, and $l$ be non-zero integers whose value may change depending on the context.
%
%\vspace{0.3cm}
\noindent $\blacktriangleright$
{\bf Conditional expectations.}
The conditional expectation $\E[\cdot|\cF_{t}]$ is denoted by $\E_{t}[\cdot]$, and $\E_{t_i}[\cdot]$ is denoted $\E_i[\cdot]$.
%The norm $\sqrt{\e[| \cdot|^2]}$ by $\normexp{\cdot}$.
We  make use of a conditional version of Fubini's theorem, stated in Lemma \ref{lem:cond:fubini}.
 {
We slightly abuse  notation by writing $\int_0^T \E_t[f_s] ds := \int_0^T F_t(\cdot,s) ds$, (likewise $\int_0^T g(\E_t[f_s]) ds :=  \int_0^T g(F_t(\cdot,s) )ds $ for any measurable function $g$) where $F_t$ is the process defined in  Lemma \ref{lem:cond:fubini}, because we believe this notation to be somewhat clearer -- in particular, this formal definition indicates more clearly that the inner integral comes from a conditional expectation than strictly mathematically correct version using the process $F_t(\cdot,s)$.
}

%\vspace{0.3cm}
\noindent
$\blacktriangleright$
{\bf Lebesgue measure}
For any Euclidean space $E$, $\cB(E)$ denotes the Borel measurable sets in $E$,
and the Lebesgue measure on the measurable space $\big(E, \cB(E) \big)$ is denoted by $m$.
%For a $\cB(E)$-measurable function $f$, the integral $m(f)$ of $f$ with respect to the Lebesgue measure is denoted by $\int_E f(x) dx$. 
%$\L_2(E; (\R^k)^\top )$ denotes the space of $\cB(E)$-measurable functions $f : E \rightarrow (\R^k)^\top $ such that $\int_E|f(t)|^2 dx$
%is finite.
%
%
%\vspace{0.3cm}
%$\L_{2}(\cF_T)$ is the space of random variables that are square integrable
%with respect to $\E$. 
%For any sub-$\sigma$-algebra $\cG\subset\cF_T$, 
%$\L_2(\cG) \subset \L_2(\cF_T)$  is the space of square integrable $\cG$-measurable random variables. 
%%For a given Hilbert space $H$, define by $\L_2(\cF_T ; H)$ the space of random variables $X$ taking values in $H$ such that $\E[ \| X\|_H^2]$ is finite.
%
%\vspace{0.3cm}
%% $dt$, $dr$, $ds$, or $dv$, depending on the context.
%The measure $m\times\P$  on the measurable space $\big( [0,T) \times \Omega, \cB([0,T)) \otimes \cF_T \big)$ denotes the product measure of $m$  and $\P$.
%$\L_{2}([0,T]\times\Omega)$ is the space of $\cB\big([0,T) \big)\otimes\cF_T$-measurable processes that are square integrable with respect to 
%$m\times \P$.
%$\L_2([0,T)\times\Omega; \R^k)$ denotes the processes in $\L_2([0,T)\times\Omega)$ taking values in $\R^k$.

\noindent $\blacktriangleright$
{\bf Processes and spaces.}
For two processes $X$ and $Y$ in $\L_0([0,T]\times \Omega; \R^k)$, $Y$ is said to be a version of $X$ if $X = Y$ $m\times\P$-a.e.
$\cP\subset\cB([0,T])\otimes\cF_T$ is the predictable $\sigma$-algebra, generated by the continuous, adapted processes,
and $\cH^2$ is the subspace of $\L_{2}([0,T]\times\Omega)$ containing only predictable processes.
For $p\ge 2$,  $\cS^p$ is the subspace of $\cH^2$ of continuous processes $Y$ such that 
$\normS{p}{Y} := (\E[\sup_{0\le s \le T}|Y_{s}|^p])^{\frac1p}$ is finite 
for all $Y\in\cS^p$;  %\gmen{\red{This is awkward! why not $\E[\sup_t |Y_t|^p]$}?}
$\normS{p}{\cdot}$ is a norm for this space.

 %on $\cS^p$ and $\cS^p$ is complete under this norm.
% The solutions of BSDEs $(Y,Z)$ lie in the space $\cS^2 \times\cH^2$.

%\vspace{0.3cm}
% The function space $\L_{2,\alpha}$ is the space of $\alpha$-fractionally smooth functions.
%\noindent $\blacktriangleright$
%We say that the function
% $x\mapsto\Phi(x)$ is in the space of fractionally smooth functions, $\L_{2,\alpha}$ if

%\vspace{0.3cm}
\noindent $\blacktriangleright$
{\bf Linear algebra}
We identify the space of $k\times n$ dimensional, real valued matrices with $\R^{k\times n}$.
$x^\top$ denotes the transpose of the vector $x$.
 {$I_n$ denotes the identity matrix in $\R^{n\times n}$.}
For any $A \in \R^{k\times n}$, let $A_j$ denote the $j$-th column vector of $A$.
%For any vector or matrix $A$, $A^\top$ denotes its transpose. For any two vector $V$ and $U$, the dot product $V^\top U$ is denoted
%by $V\cdot U$.
%In stochastic integrals where the integrator $V$ is $\R^k$ valued and integrand $U$ is $(\R^k)^\top$-valued, 
% we write $\int U_t dV_t$ to mean $\sum_{i=1}^k \int U_{i,t} dV_{i,t}$.
% In particular, the stochastic integral $\int Z_t dW_t$ means $\sum_{l=1}^q \int Z_{l,t} dW_{l,t}$.
%When the  integrator $V$ is $\R^k$ valued and integrand $U$ is $(\R^{n\times k})^\top$-valued, 
% we write $\int U_t dV_t$ to mean the $\R^n$-valued random variable whose $i$-th component is $\int U_{i,t} dV_{t}$,
% where $U_i$ is the $i$-th row of $U$.
%
%\vspace{0.3cm}
For any vector $x\in\R^n$, $|x|$ is the vector 2-norm, defined by $(\sum_{i=1}^n |x_i|^2)^{1/2}$, and
for any matrix $A$, $|A|$ is the matrix 2-norm, defined by $\max_{|x|=1} |Ax|$, where $|Ax|$ is the vector 2-norm of the vector Ax.

%\vspace{0.3cm}
\noindent $\blacktriangleright$
{\bf Functions and regularity.}
Let $\gamma \in (0,1]$ and  $A(\cdot)$ be a function in the domain $[0,T)\times \R^l$ taking values in $ \R^{k\times n}$ (resp. $\R^k$). 
We say that $A(t,\cdot)$ is $\gamma$-H\"older  continuous uniformly in $t$ with H\"older constant $L_A $ if, 
for all $(x,y) \in (\R^l)^2$ and $t\in[0,T)$,  $|A(t,x) - A(t,y)| \le L_A |x-y|^{\gamma}$; in the case that $\gamma = 1$, we say that 
$A(t,\cdot)$ is Lipschitz continuous uniformly in $t$ with Lipschitz constant $L_A$. 
Likewise, we say that $A(\cdot,x)$ is $\gamma$-H\"older continuous uniformly in $x$ with H\"older constant $L_A $ if,
for every $(t_1,t_2) \in [0,T)^2$ and $x\in \R^l$, $|A(t_1,x) - A(t_2,x)| \le L_A |t_1-t_2|^{\gamma}$.
For a given multi-index $\alpha = (i_1, \dots, i_{|\alpha|})$ with no zero entries, we define by $\partial^\alpha_x A(t,\cdot)$ the multiple derivative $\partial_{x_{i_1}}\dots \partial_{x_{i_{|\alpha|}}} A(t,\cdot)$.
%We say that $A(t,\cdot)$ is continuously differentiable if the partial derivative $\partial_{x_j} A_{u,v}(t,x)$ (resp. $\partial_{x_j} A_u(t,x)$ )
%of the $(u,v)$-th component $A_{u,v}(t,x)$  (reps. $u$-th component $A_u(t,x)$) 
%of $A(t,x)$ exists and is continuous for every $(u,v) \in \{1,\dots,k\}\times\{1,\dots,n\}$, $j \in \{1,\dots,l\}$ and $(t,x) \in [0,T)\times \R^l$.
If $A(t,\cdot)$ takes values in $\R^k$ and is differentiable, we define by $\nabla_x A(t,\cdot)$ the $\R^{k \times l}$ valued function
whose $(u,v)$-th component is $\partial_{x_v} A_u(t,\cdot)$.
If $A(t,\cdot)$ takes values in $(\R^k)^\top$ and is differentiable, we define by $\nabla_x A(t,\cdot)$ the $\R^{l \times k}$-valued function whose $(u,v)$-th component is $\partial_{x_u} A_v(t,\cdot)$.
Define by $\|A \|_\infty$ the infinity norm 
\[ \max_{u,v} \  \sup_{(t,x) \in [0,T) \times\R^l }|A_{u,v}(t,x)| \quad  
\text{(resp. } \ \max_{u} \  \sup_{(t,x) \in [0,T) \times\R^l } |A_{u}(t,x)| \text{).} \]

\vspace{0.3cm}
% Finally, we will often make the transition between a process and it's
% {\em predictable projection} \cite[Paper IV, Theorem 5.6]{revu:yor:99}.
% For example, for a predictable process $P$, we will make use of the equality
% $\E_{t}[\int_{0}^TP_{s} ds] = \int_{0}^TV_{s}ds$, where 
% $V$ is the unique predictable modification of $(s,\omega)\mapsto \E_{t}[P_{s}](\omega)$.
% However, to take advantage of this equality, we will be making computations on the marginals $V_{s} = \E_{t}[P_{s}]\; a.s.$, whence it will be notationally more convenient to write
% $\E_{t}[\int_{0}^T P_{s} ds] = \int_{t}^T \E_{t}[P_{s}]ds$ and allow the context to determine whether $\E_{t}[P_{s}]$ or it's predictable projection is in use.

%\vspace{0.3cm}
%For a given time-grid $\pi = \{0 = t_0<\ldots<t_N=T\}$, define the $i$-th time-increment by $\Delta_i := t_{i+1} - t_i$.
%% and the number of points by $\#\pi := N$.
%For a parameter $\beta\in(0,1]$, the time grids $\pib_N$ are those with $N$ time-points whose $i$-th time point is $t_i = T - T(1-i/N)^{1/\beta}$ for 
%$i\in\{0,\ldots,N\}$.

%\vspace{0.3cm}
%At several points, it will be necessary to apply a mollification procedure to continuous functions. 
\noindent $\blacktriangleright$
{\bf Mollifiers.}
The following definitions will come in handy.
\begin{definition}
\label{def:1:moll}
Let $n$ be a non-zero integer. A mollifier is a smooth function  $\phi: \R^n \rightarrow [0,\infty)$ with compact support on 
$\{x : \in \R^n \ : \ |x| \le 1\}$ such that $\int_{\R^n}\phi(x) dx = 1$ 
and $\lim_{
%M
R \rightarrow\infty} R
%M
^{n}\phi(R
%M
x) = \delta(x)$ for all $x\in \R^n$, where  $\delta(x) $ is the Dirac delta function. 
For $R > 0$, define the function $\phi_{R
%,M
} : \R^n \rightarrow [0,\infty)$ be the function $ x \mapsto R
%M
^{n}\phi(R
%M
x)$.
\end{definition}
An example of a mollifier is  {$\phi(x) = e^{-1/(1-|x|)} \1_{|x| < 1}/ \int_{|x|<1} e^{-1/(1-|y|)}dy$.}
% where $B(1) := \{y \in R^n \ : \ |y| \le 1 \}$.} 
The following lemma, which is standard, shows how a mollifier can be used to generate a smooth function from a continuous one.
%, and shows that all continuous functions can be approximated by smooth functions.
\begin{lemma}
\label{lem:1:moll}
Let $F:\R^n \rightarrow\R$ be continuous, 
and define the function
$F_R
%M
(x) := \int_{\R^n} F(x - y) \phi_R
%M
(y) dy$. Then the function $F_R
%M
(x)$ is smooth and $\lim_{R
%M
\rightarrow\infty} F_R
%M
(x) = F(x)$ for
all $x\in \R^n$.
\end{lemma}
%The proof is standard, but we include it for completeness.
%
%\noindent {\bf Proof.} Let $E_{R
%%,M
%}$ be the compact support of $\phi_{R
%%,M
%}$. Fix $x\in\R^n$.
%\begin{align*}
%| F(x) - F_R
%%M 
%(x)| = | \int_{\R^n} (F(x) - F(x - y) ) \phi_{R
%%,M
%}(y) dy | \le  \int_{E_{R
%%,M
%}} | F(x) - F(x - y) | \phi_{R
%%,M
%}(y) dy.
%\end{align*}
%%where $E_{R,M}(x)$ is the set $\{y \in \R^n \ : \ y = z + x, \ \ z \in E_{R,M} \}$.
%Since $\sup_{z,y \in E_{R
%%,M
%} } |y-z|$ converges to zero as $R$
%%$M$ 
%increases to infinity, and since $E_{R_1} \subset E_{R_2}$ whenever $R_2 > R_1$,
%the result follows from the continuity of $F$.

\subsection{Assumptions}
\label{sectio:1:ass}
The following assumptions will hold throughout this paper.
\begin{enumerate}
\item[\Hs] %The coefficients of the SDE \eqref{sde} satisfy the following properties.
$X$ is a solution to the stochastic differential equation (SDE)
\begin{equation}
\label{eq:1:SDE}
X_0 = x_0, \quad X_t = x_0 + \int_0^t b(s,X_s) ds + \int_0^t \sigma(s,X_s) dW_s \quad t > 0, 
\end{equation}
where $x_0 \in \R^d$ is fixed and $b$ and $\sigma$ satisfy
\begin{enumerate}
 \item $(t,x) \in [0,T]\times\R^d \mapsto b(t,x)$ is $\R^d$-valued, measurable and uniformly bounded.
Moreover, $b(t,\cdot)$ is twice continuously differentiable with uniformly bounded derivatives and H\"older continuous second derivative,
and $b(\cdot,x)$ is $1/2$-H\"older continuous uniformly in $x$.
\item $(t,x) \in [0,T]\times\R^d \mapsto \sigma(t,x)$ is $\R^{d\times q}$-valued, measurable and uniformly bounded.
Moreover, $\sigma(t,\cdot)$ is twice continuously differentiable with uniformly bounded derivatives and H\"older continuous second derivative,
and $\sigma(\cdot,x)$ is $1/2$-H\"older continuous uniformly in $x$.
\item $\sigma(\cdot)$ satisfies a uniformly elliptic condition: there exists some finite $\elip>0$ such that, for any $\zeta\in \R^d$, 
$\zeta^\top \sigma(t,x) \sigma(t,x)^\top \zeta \ge \elip|\zeta|^2$ for all $(t,x) \in [0,T]\times\R^d$.
\end{enumerate}
\item[\Hg] 
The terminal condition
$\Phi:\R^d \rightarrow \R$ is a measurable function and there exists a constant $\alpha\in(0,1]$ such that $K^\alpha(\Phi) < \infty$, where
\begin{equation}
\left.
\begin{array}{rl}
K^\alpha(\Phi)^2 & := \E[|\Phi(X_T)|^2] + \sup_{0\le t < T} \frac{V_{t,T}(\Phi)^2}{(T-t)^\alpha}\\ \\
\text{for } V_{t,T}(\Phi)^2 & : = \E[|\Phi(X_T) - \E_t[\Phi(X_T)]|^2].
\end{array}
\right\}
\label{eq:cdn:var}
\end{equation}
We say that $\Phi$ is fractionally smooth, and that it belongs to the space $\L_{2,\alpha}$.
%We remark that any $\alpha$-H\"older continuous function $g$ belongs to $\L_{2,\alpha}$, but that there are $\alpha$-H\"older continuous functions $g$
%for which there is a $\beta \in (\alpha,1]$ such that $g\in\L_{2,\beta}$; see \cite[Section 2]{geis:geis:gobe:12} for examples. 
We refer to  \cite{gobe:makh:10} for further discussion of and references for the space $\L_{2,\alpha}$.
% There is an $\alpha\in(0,1]$ such that the terminal function $\Phi$ is in the space $\L_{2,\alpha}$.
% satisfies the following properties.
% \begin{itemize} 
%  \item[ i)] There exists $\alpha\in(0,1]$ such that $ \Phi$ is in the space $\L_{2,\alpha}$.
% \item[ii)] $\Phi$ is uniformly bounded by $C_\Phi$. 
% \end{itemize}
\item[\Hf] The driver $f : [0,T) \times \R^d \times \R \times (\R^q)^\top \rightarrow \R$ satisfies \eqref{eq:loclip:driver}.
\end{enumerate}

\vspace{0.3cm}
 {
The following condition will be required for both the Euler scheme and the Malliavin weights scheme convergence results; this is a standard assumption for BSDE approximation schemes in order to obtain a convergence bounded from above by $O(N^{-1})$.
\begin{enumerate}
\item[\HFt] The driver $f(t,x,y,z)$ is  $\frac12$-H\"older continuous in its $t$ uniformly in $(x,y,z)$ with H\"older constant $L_f$.
\end{enumerate}
}

\vspace{0.3cm}
 {
Our convergence results for the Malliavin weights scheme require stronger conditions than those of the Euler scheme;
one of the following assumptions will be necessary to obtain the main result, Theorem \ref{thm:1:disc er}, of Section \ref{section:num}.
\begin{enumerate}
\item[\Hgexp] The terminal condition has exponential bounds in the sense that there is a finite $C_\xi >0$ such that $\E[e^{|\Phi(X_T)|}] \le C_\xi$.
\item[\Hgh] The function $\Phi$ is H\"older  continuous: there exists a finite constants $K_\Phi$ and $\thetaP \in (0,1]$ such that $|\Phi(x_1) - \Phi(x_2)| \le K_\Phi |x_1 - x_2|^{\thetaP}$
for any $x_1,x_2 \in \R^d$.
\end{enumerate}
}

\vspace{0.3cm}
The following assumptions will be needed for partial results only. They will  hold only when specifically stated.
\begin{enumerate}
\item[\HFd] The driver $(t,x,y,z) \mapsto f(t,x,y,z)$ is continuously differentiable with respect
$(x,y,z)$ for all $t\in[0,T)$. The partial derivatives  {in $(y,z)$} are bounded by $L_f(T-t)^{(\thetaL -1)/2}$  {and the partial derivatives in $x$ are bounded above by $L_X(T-t)^{1-\theta_X/2}$}. 
\item[\Hge] The function $\Phi$ is uniformly bounded: $\|\Phi\|_\infty < \infty$.
\end{enumerate}
\begin{remark}
\label{rem:holder}
Due to \Hs, \Hgh \ implies \Hgexp \ and \Hg. Note that it is possible that $\thetaP < \alpha$: see \cite[page 2086, e.g. (i)]{geis:geis:gobe:12}.
\end{remark}

\vspace{0.3cm}
In the proofs below, it will be necessary to compute a right-inverse to the matrix $\sigma(\cdot)$, i.e., 
for every $(t,x) \in [0,T) \times \R^d$, it will be necessary to find a
$(q , d)$-dimensional matrix $\sigma^{-1}(t,x)$ such that $\sigma(t,x) \sigma^{-1}(t,x) = I_d$.
 In the case where the dimensions $d$ and $q$ are equal, this is uniquely defined by usual matrix inverse of $\sigma(t,x)$,
whose existence is guaranteed by the uniform ellipticity condition \He.
If the dimensions $d$ and $ q$ are not equal, $\sigma^{-1}(t,x)$ is defined by the pseudoinverse 
$\sigma(t,x)^\top \big( \sigma(t,x) \sigma(t,x)^\top \big)^{-1}$;
this is well defined because the uniform ellipticity condition \He \ guarantees the existence of the inverse of $\sigma \sigma ^\top$.

% Recalling that $\sigma$ is uniformly elliptic, there is a $d\times d$-dimensional, symmetric, invertible  matrix $\Sigma$
% \cite[Lemma 5.2.1]{stro:vara:79}
% such that $\Sigma^2 = \sigma\sigma^\top$. 
% % One can then replace $\sigma^{-1}$ with $\Sigma^{-1}$ where it is required without changing the results.
% Since the diffusion $X$ depends only on $\sigma\sigma^\top$ \cite[Theorem 5.3.2]{stro:vara:79}, we replace 
% $\sigma$ by $\Sigma$ in these cases; by slight abuse of notation, and to avoid confusion in the notation,
% we will always refer to $\sigma$.

\section{Key preliminary results}

\subsection{Malliavin calculus}
\label{section:mal:cal}
We recall briefly some properties and definitions of Malliavin calculus.  For details, we refer the reader to \cite{nual:95}.

\vspace{0.3cm}
For any $m \ge 1$, define $C^\infty_p (\R^m)$ to be the space of functions taking values in $\R$ which are infinitely differentiable such that all partial derivatives have at most polynomial growth,
and denote by $W(h) := \int_0^T h_t dW_t$ the It\^o integral of the $(\R^q)^\top$-valued, deterministic function $h \in \L_2([0,T) ; (\R^q)^\top )$.
Let $\cR \subset \L_2(\cF_T)$ be the subspace containing all random variables $F$ of the form $f(W(h_1),\ldots,W(h_m))$ for $h_i \in \L_2([0,T); \R^q) $ and any finite $m$.
%\[F = f(W(h_1),\ldots,W(h_m)), \quad f \in C^\infty_p(\R^m),\; h_i \in \L_2([0,T); \R^q) \quad \forall m\in \N.\]
Define the derivative operator $D : \cR \mapsto \L_2([0,T]\times \Omega)$ by 
$D_t F := \sum_{i=1}^m \partial_i f(W(h_1),\ldots,W(h_m)) h_i(t).$
The derivative operator is extended to $\mal\subset\L_2(\cF_T)$, the closure of $\cR$ in $\L_2(\cF_T)$ under the norm
$\|F\|_{1,2}^2 := \normexp{F}^2 + \E[\int_0^T|D_tF|^2 dt]$.
Define by $\mal(\R^k)$ (resp. $\mal( (\R^k)^\top)$)  by the space of random variables $F = (F_1,\dots, F_k)^\top$ 
(resp. $F = (F_1,\dots, F_k)$)
such that 
$F_i \in \mal$ for each $i \in\{0,\dots,k\}$.
The Mallivin derivative $D F$ is denoted by the $\R^{k\times q}$- (resp. $\R^{q\times k}$-) valued process whose $i$-th row 
(resp. column) is $D F_i$ (resp. $(D F_i)^\top$).
%For non-zero integer $k$, define $C^\infty_p (\R^m;\R^k)$ to be the space of functions $f = (f_1,\dots,f_k)^\top$ taking values in $\R^k$ 
%such that, for each $i\in\{1,\dots,k\}$, $f_i$ is infinitely differentiable such that all partial derivatives have at most polynomial growth.
%In the case that $f(\cdot)$ takes values in $\R^k$, for some non-zero integer $k$, $D_tF$ is defined by the $k\times m$ random matrix
%whose $(u,v)$-th component is
%$\sum_{i=1}^m \partial_i F_u(W(h_1),\ldots,W(h_m)) h_{i,v}(t)$,
%where $F_u(\cdot)$ is the $u$-th component of $F(\cdot)$ and $h_{i,v}(\cdot)$ is the $v$-th component of $h_i(\cdot)$.
%If $F(\cdot)$ takes values in $(\R^k)^\top$, then $D_tF$ is defined by $\big( D_t F^\top)^\top$.

The following lemma, termed the {\em chain rule} of Malliavin calculus, is proved in \cite[Proposition 1.2.3]{nual:95}.
\begin{lemma}[Chain rule]
\label{lem:mal:chain rule}
Let $(F_1,\dots,F_m)\in(\mal)^m$. For any continuously differentiable function $f : \R^m \rightarrow \R$ with bounded partial derivatives,
and $F = f (F_1,\ldots,F_m) \in \mal$, the random variable $f(F) \in \mal$ and
$Df(F) = \sum_{i=1}^m \partial_i f(F) DF_i = \nabla_x f(F) DF$.
\end{lemma}
\begin{remark*}
In the case that $F$ takes values in $(\R^m)^\top$, the result of Lemma \ref{lem:mal:chain rule} hold with $Df(F) =  \nabla_x f(F) (DF)^\top$.
In the case that $f$ takes values in $(\R^k)^\top$, applying Lemma \ref{lem:mal:chain rule} component-wise yields that $f(F)$ is in 
$\mal((\R^k)^\top)$ and $Df(F) = (DF)^\top \nabla_x f(F)$.
\end{remark*}

For the space \[\dom := \{ u \in \L_{2}([0,T]\times\Omega ; (\R^q)^\top ) : \exists c\in \R \; s.t \;
\forall F\in\mal \; | \E[\int_{0}^T ( u_{s} \cdot D_{s}F ) ds] | \le c \normexp{F}^2\;\}\]
define the Skorohod integral operator $\delta : \dom \rightarrow \L_{2}(\Omega)$ as the dual operator to the Malliavin derivative in the sense that
$\E[\int_{0}^T ( u_s \cdot D_{s}F ) ds] = \E[F\delta(u)].$
%The following Lemma, relating the Skorohod integral to It\^o's integral, is given in \cite[Proposition 1.3.11]{nual:95}
Below are the key properties of the Skorohod integral used in this paper.
%\begin{lemma}[Skorohod integral versus It\^o integral]
%\label{lem:skor:adapted}
% If $u\in\dom$ is adapted, $\delta(u) = \int_{0}^T u_{s}dW_{s}$, the It\^o integral of $u$.
%\end{lemma}

%The following Lemma, termed the {\em integration by parts formula} of Malliavin calculus, is proved  in \cite[Paper 1.3.1 (4)]{nual:95}.
\begin{lemma}[Integration-by-parts]
 \label{lem:mal:ibp}
 Suppose that $u\in\dom$ and $F\in\mal$ are such that 
$\E[F^2\int_{0}^T |u_{s}|^2 ds] < \infty$. Then, the {\em integration by parts formula}  holds:
$\int_{0}^T (u_{s} \cdot  D_{s} F) ds = F\delta(u) - \delta(Fu).$
\end{lemma}

%This following result is used and proved in the proof of \cite[Theorem 4.2]{ma:zhan:02}. 
%%We include the proof here in greater detail  for the convenience of the reader. 
%\begin{lemma}
%\label{lem:mal:cdn exp}
%Let $u\in\dom$. Then, for all $t\in[0,T)$,  $\E_{t}[\delta(\1_{(t,T]}(\cdot)u)] = 0$.
%\end{lemma}
%
%{\bf Proof.}
%Let $F$ be $\cF_{t}$-measurable, square integrable random variable.
%We will show that \\ $\E[F\delta(\1_{(t,T]}(\cdot)u_\cdot)] = 0$, whence the proof is completed by taking $F = \1_{A}$ for any $A\in\cF_{t}$.
%Using \cite[Theorem 1.1.1]{nual:95}, it suffices to prove that $\E[H_m\big(W(h) \big)\delta(\1_{(t,T]}(\cdot)u_\cdot)] = 0$ for the Hermite
%polynomial $H_m(\cdot)$ of degree $m$ for any $m \ge 0$
%and $h(\cdot) $ of the form $ \tilde h(\cdot) \1_{[0,t]}(\cdot)$ for any bounded $\tilde h\in \L_2([0,T]; (\R^q)^\top)$.
%Suppose first that $m\ge 1$.
%The chain rule, Lemma \ref{lem:mal:chain rule} -  yields $D_s H_m\big( W(h) \big) = H_{m-1} \big( W(h) \big) h(s)$ for every $s\in [0,T]$.
%Using the definition of the Skorohod integral,
%\[
% \E[H_m\big(W(h) \big)\delta(\1_{(t,T]}(\cdot)u_\cdot)] = \E[ \int_t^T u_s D_s H_m \big( W(h) \big) ds] 
%= \E[ \int_t^T u_s  H_{m-1} \big( W(h) \big) h(s) ds] = 0
%\]
%as required. The case $m=0$ is simpler, because $H_0\big(W(h)\big)$ is $1$ and so its Malliavin derivative is 0.
%\qed

\begin{remark}
\label{rem:skor:multidim}
Suppose that the process $u$ takes values in $\R^{q\times k}$ is such that $u_i^\top$ is in $\dom$ for each 
$i\in\{0,\dots,k\}$, where $u_i$ is the $i$-th column of $u$. The Skorohod integral of $u$, denoted by $\delta(u)$, is defined by
\begin{equation}
\label{skor:multidim}
\delta(u) := \big( \delta(u_1^\top)^\top , \dots , \delta(u_k^\top )^\top \big).
\end{equation}
%In the case where $u$ is adapted, Lemma \ref{lem:skor:adapted} implies that $\delta(u)$ is equal to $(\int_0^T u_t^\top dW_t )^\top$.
The integration by parts formula, Lemma \ref{lem:mal:ibp}, is applied column-wise in the case of matrix valued $u$.
%It follows that, if $F \in \mal$ and $u$ satisfy $\E[|F|^2 \int_0^T |u_t|^2 dt]$ is finite, then
%\[
%\int_0^T ( D_s F u_s ) ds = F\delta(u) + \delta(Fu)
%\]
where $D_sF u_s$ is understood as a matrix-matrix multiplication, and the Skorohod integrals are defined in the multidimensional sense
of equation \eqref{skor:multidim}.
\end{remark}
%Finally, we define the space
%\begin{align*}
%\mala  := \bigg\{ & u \in \L_2([0,T]\times\Omega) \; : u_t\in\mal \ m-a.e., 
%%\\& 
%\exists \ \text{measurable version of } (s,t,\omega)\mapsto D_s u_t (\omega) \ s.t. \\
%&\quad \E[\int_0^T \int_0^T |D_su_t|^2dsdt] < \infty \bigg\}.
%\end{align*}

\subsection{SDEs and Malliavin calculus}
\label{section:mal:SDE}
Fix $t \in [0,T)$ and $x \in \R^d$.
We recall some standard properties on the Malliavin calculus applied to SDEs $\XM tx$ of the form
\begin{equation}
\label{eq:flow}
\XM tx_s  =  x +  \int_t^s b(r,\XM tx_r)\1_{(t,T]}(r) dr + \int_t^s \sigma(r,\XM tx_r)\1_{(t,T]}(r) dW_r .
\end{equation}
Observe that the SDE $X$ defined in \eqref{eq:1:SDE} is equal to $\XM 0{x_0}$.
%
%\vspace{0.3cm}
%\begin{lemma}
%\label{lem:mal:SDE:1}
%The random variables $\XM tx_r$ - the marginals of the process $X$ - are in $\mal(\R^d)$ for each $r \in [0,T]$.
%For all $s$, $(D_s \XM tx_r)_{r\ge s}$ solves the linear SDE
%\[V_r = \1_{[s,T]}(r)\Big\{\sigma(s,\XM tx_s) + \int_s^r \1_{[t,T]}(u)\nabla_xb(u,X_u)V_u du 
%+ \sum_{j=1}^q \int_s^r \1_{[t,T]}(u) \nabla_x \sigma_j(u,\XM tx_u)V_u dW_{j,u} \Big\} .\]
%where $\sigma_j(\cdot)$ is the $j$-th column of $\sigma(\cdot)$.
%Moreover, for every $p\ge 2$, there is a constant $C_p$ depending only on $\|\nabla_x b\|_\infty$, $\|\sigma\|_\infty$, $\|\nabla_x\sigma_j\|_\infty$, $T$ and $p$ such that  $\E[\sup_{s \le r \le T} |D_sX_r|^p] \le C_p$.
%\end{lemma}
%{\bf Proof.}
%The relation to the SDE is proved in \cite[Theorem 2.2.1]{nual:95}.
%The bound is proved using the proof method of \cite[Theorem IX.2.4]{revu:yor:01} (essentially using Gronwall's inequality).
%\qed
%
%
%\vspace{0.3cm}
First, we recall the flow  $\nabXM tx{ }$ and its inverse $\inXM tx{ }$, which are respectively 
defined as the solutions to the SDEs
\begin{align*}
\nabXM tx{r}& = I_d + \int_{t}^r 
%b'_u
\nabla_x b(u,\XM tx_{u}) 
\nabXM tx u du + \sum_{j=1}^q \int_{0}^r 
%\sigma'_{j,u}
\nabla_x\sigma_j(u,\XM tx_{u})
% \nabla_x\sigma_j( u ,\XM tx_{u}) 
 \nabXM txu dW_{j,u}, \\
\inXM tx{r} & = I_d + \int_{t}^r \inXM tx u \big( \sum_{j=1}^q (
\nabla_x\sigma_j(u,\XM tx_{u})
%\nabla_x\sigma_j(r,\XM tx_{r}
)^2 - 
%b'_u
\nabla_xb(u,X_{u})
\big)  du \\
& \qquad  - \sum_{j=1}^q \int_{t}^r  \inXM tx u
\nabla_x\sigma_j(u,\XM tx_{u})
% \nabla_x\sigma_j(r,X_{r})
  dW_{j,u},
\end{align*}
%where $b'_u$ denotes $ \nabla_xb( u ,\XM tx_{u}) $ and $\sigma'_{j,u}$ denotes $\nabla_x\sigma_j(r,\XM tx_{r})$ 
 {where $\sigma_j$ is the $j$-th column of $\sigma$. } 
These processes are linear SDEs, and we list some standard properties used throughout this paper  in the following Lemma.%; the proof is standard from the properties of linear SDEs.
%where $\sigma_j$ is the $j$-th column of $\sigma$.
%The notation and results of \cite[p325 - 327]{prot:04} are required to make the connection between the SDEs: for a $\R^{d\times d}$-valued continuous semimartingale $(\cZ_t)_t$,
%recall that the {\em left} (resp. {\em right}) stochastic exponential $\cE(\cZ)_t$ is the solution to the linear SDE
%\[
%\cE(\cZ)_t = I_d + \int_0^t \cE(\cZ)_r d \cZ_r \qquad \Big(\text{resp. } \ \cE^R(\cZ)_t = I_d + (\int_0^t \cE(\cZ)_r^\top d \cZ_r^\top)^\top \Big).
%\]
%Now, taking the $\R^{d\times d}$-valued, continuous semimartingales $(\cW_t)_t$ given by \[ \cW_t = \int_0^t \nabla_{x} b(r,X_{r}) dr 
%+ \sum_{j = 1}^q \int_{0}^t \nabla_{x}\sigma_j(r,X_{r}) dW_{j,r}.\]
%it follows that
%\begin{equation}
%% \left.
%% \begin{array}{l}
% \nabX t  = \cE^R(\cW)_t, \quad 
% \inX{t}  = \cE(- \cW + <\cW>)_t, 
%% I_d + \int_{0}^t (\nabla_x\sigma(r,X_{r})^2 - \nabla_xb(r,X_{r})) \inX r dr \\
%% \qquad \qquad \qquad - \sum_{j=1}^q \int_{0}^t \nabla_x\sigma_j(r,X_{r}) \inX r dW_{j,r}
%% \end{array}
%% \right\}
%\label{eq:grad:proc}
%\end{equation}
%where $<\cW>$ is the quadratic variation process of $\cW$.
\begin{lemma}
\label{lem:mal:SDE:2}
For every $p> 1$, $\nabXM tx{ }$ and $\inXM tx{ }$ are in $\cS^p$, and 
there is a constant $C_p$ depending only on $ \| \sigma\|_\infty$, $\|\nabla_x b\|_\infty$, $\|\nabla_x\sigma_j\|_\infty$, $T$ and $p$ such that  
\[\normS{p}{\nabXM tx{ } } + \normS{p}{\inXM tx{ } } \le C_p.\]
Moreover, 
\[\normexp{\nabXM tx{r} - \nabXM tx{s}}^2 + \normexp{\inXM tx{r} - \inXM tx{s}}^2 \le C_2 |r-s|\] for all  $(t,s) \in [0,T]^2$,
  \[ \nabXM txr \inXM txr = I_d \] 
 for all $ r \in[t,T] $ almost surely, and, for any $r<u<s$, 
\[\E_r[|\nabXM txs \inXM txr - \nabXM txu \inXM txr  |^2] \le C_2(s-u) \quad \P-\text{a.s.} \]
\end{lemma}
%Lemma \ref{lem:mal:SDE:2} follows  from standard estimates on linear SDEs with bounded coefficients. The bound on $\E_r[|D_s X_t|^p]$ follows from the fact that $\nabX t \inX s$ solves the SDE
%\[
%U_t = I_d + \int_{s}^t \nabla_xb(r,X_{r}) U_r dr + \sum_{j=1}^q \int_{s}^t  \nabla_x\sigma_j(r,X_{r}) U_r dW_{j,r} \quad \text{for all } t \in [s,T].
%\]
%{\bf Proof.}
%That $\nabX{ }$ and $\inX{ }$ are in $\cS^p$, and the bounds, can be proved using  
%the proof method of \cite[Theorem IX.2.4]{revu:yor:01} (essentially using Gronwall's inequality).
%That $\nabX{t}\inX{t}$ is the identity for all $t$ almost surely follows immediately from \cite[Paper 5 Theorem 48]{prot:04}.
%To prove the second property, observe from \eqref{eq:grad:proc} that $\nabX{t} = \nabX{s} + (\int_s^t ( \nabX{r} )^\top \cW_r^\top ) ^\top $.
%Using the invertibility property, it follows that $(\nabX{t}\inX{s}  \sigma(s,X_s ) )_{t\ge s}$ solves the SDE
%\begin{align*}
%\cX_t & = \sigma(s,X_s) + (\int_s^t \big( \cX_r \big)^\top \cW_r^\top ) ^\top.
%\end{align*}
%If one expands out the definition of $\cW_r$, one sees that the process $\big( \nabX{t}\inX{s}\sigma(s,X_s) \1_{[s,T]}(t) \big)_{t\ge s}$
%solves the same linear SDE as $(D_sX_t)_{t \ge s}$ given in Lemma \ref{lem:mal:SDE:1}, whence the result follows from
%the uniqueness of solutions of linear SDEs.
%% Since $D_sX_t$ satisfies a linear SDE, Lemma \ref{lem:mal:SDE:1}, this
%%  is due the representation of $D_sX_t$ as a stochastic exponential.
%\qed

\vspace{0.3cm}
%A direct consequence of Lemma \ref{lem:mal:SDE:2} is the following result.% bound on $\sup_s\normS{2}{D_sX}$.
The Malliavin derivative of the marginals of $\XM tx$ is strongly related to the flow and its inverse, as shown in the following Lemma.
The proof of the estimates follows directly from Lemma \ref{lem:mal:SDE:2}.
\begin{lemma}
\label{lem:mal:SDE:3}
%For any $p\ge 2$,
For all $r \in [0,T]$, $\XM tx_r$ is in $\mal (\R^d)$
 {
and there is a version $D_s \XM tx_r$ satisfying the SDE
\[
D_s \XM tx_r \1_{[s,T]}(r) = \Big\{\sigma(s,\XM tx_s) + \int_s^r \nabla_xb(\tau,\XM tx_\tau)D_s\XM tx_\tau d\tau 
+ \sum_{j=1}^q \int_s^r \nabla_x \sigma_j(\tau,\XM tx_\tau)D_s\XM tx_r dW_{j,r} \Big\} .
\]
}
Moreover,
for all $0\le s,r\le T$,
\[D_{s}\XM tx_{r} = \nabXM tx r \inXM tx s \sigma(s,\XM tx_{s})\1_{[s,T]}(r) \1_{[t,T]}(s) \quad a.s.\]
whence there exists a constant $C_p$ depending only on $ \| \sigma\|_\infty$, $\|\nabla_x b\|_\infty$, $\|\nabla_x\sigma_j\|_\infty$,  $T$ and $p$ such that
%for all $0 \le r \le s < t \le T $, 
$\E_s[|D_s \XM tx_r|^p] \le C_p$, 
and
 $\sup_s\E[\sup_{s\le r\le T} |D_s\XM tx_r|^2]^{1/2} \le C_2$;
 moreover, for any $x_1,x_2 \in \R^d$, $\E[|D_s \XM t{x_1}_r  - D_s \XM t{x_2}_r|^p ] \le C_p |x_1 - x_2|^p $
 and, for any $r<u<s$,  $\E_r[|D_r \XM tx_s  - D_r \XM tx_u  |^2] \le C_2(s-u)$.
% \pmen{need to include continuity in space}
% Moreover, if $s,t\le u$, $\normexp{D_sX_u - D_tX_u}^2 \le C|t-s|$ for some constant depending on $\|\sigma\|_\infty$, $\|\nabla_xb\|_\infty$, 
% $\|\nabla_x\sigma\|_\infty$, and the H\"older coefficient of $t\mapsto\sigma(t,x)$ only.
\end{lemma}
% {\bf Proof.}
% This is due to Lemma \ref{lem:mal:SDE:2}, the $1/2$-H\"older regularity of $t\mapsto\sigma(t,x)$, the uniform boundedness of $\sigma$,
% and the H\"older regularity of the paths of $(\inX{s})_s$. \qed

\subsection{Existence, uniqueness, approximation and decomposition  of the BSDE}
\label{section:mal:BSDE}
Since the class of BSDEs under local conditions has, to the best of our knowledge, not been studied in full generality,
we now include a proof of the existence and uniqueness of solutions.
 {We remark that the existence and uniqueness follows also from \cite[Theorem 3.2]{fan:jian:12}. 
The proof below is simpler, since a simpler class of BSDEs is considered, and different, so we include for the interest of the reader.}

\begin{theorem}
\label{cor:1:exist and uni}
There exists a unique pair of process $(Y,Z)$ in $\cS^2 \times \cH^2$ solving the BSDE \eqref{eq:1:BSDE} with terminal condition $\Phi(X_T)\in \L_2(\cF_T)$
and driver $f$ satisfiying the locally Lipschitz continuous and boundedness of \eqref{eq:loclip:driver}.
\end{theorem}
{\bf Proof.}
Let $(\phi,\psi)$ be in $\cH^2 \times \cH^2$, and define the random function
%\[
$f(r,y,z) = f(r) := f(r,X_r,\phi_r,\psi_r).$
%\]
 {We show that there exists a unique solution $(Y^{(\phi,\psi)},Z^{(\phi,\psi)})$  to the BSDE
\[
Y^{(\phi,\psi)}_t = \Phi(X_T) + \int_t^T f(r) dr - \sum_{j=1}^q \int_t^T Z^{(\phi,\psi)}_{j,r} dW_{j,r}.
\]
in $ \cH^2 \times \cH^2$ (in fact, $Y^{(\phi,\psi)}$ is in $\cS^2$). 
This will imply the function $\Xi : \cH^2 \times \cH^2 \rightarrow \cH^2 \times \cH^2$ mapping $(\phi,\psi)$ to $(Y^{(\phi,\psi)},Z^{(\phi,\psi)})$
 is well defined. 
 For this, we use  \cite[Theorem 4.2]{bria:dely:hu:pard:stoi:03}.
 The function $f$ is predictably measurable; we must show that $f$ satisfies assumptions (H1)-(H5) of \cite[Section 4]{bria:dely:hu:pard:stoi:03}.
Since $f$ takes no argument in $(y,z)$, it is only necessary to check (H1), which follows readily 
 the local Lipschitz continuity and local boundedness of the driver \eqref{eq:loclip:driver}.
 Therefore, from \cite[Theorem 4.2]{bria:dely:hu:pard:stoi:03}, $(Y^{(\phi,\psi)},Z^{(\phi,\psi)})$ exists and is unique. }
As in the proof of \cite[Theorem 2.1]{elka:peng:quen:97}, we prove that $\Xi$ is a contraction.
For $k \in \{1,2\}$, let $(\phi_k , \psi_k) \in \cH^2 \times \cH^2$ and define the BSDE $(Y_k,Z_k) := \Xi(\phi_k,\psi_k)$.
Define the differences $\delta Y = Y_1 - Y_2$, $\delta Z = Z_1 - Z_2$, $\delta \phi = \phi_1 - \phi_2$ and $\delta \psi = \psi_1 - \psi_2$.
It then follows from H\"older's inequality that 
\begin{align*}
\normexp{\delta Y_t}^2 + \int_t^T \normexp{\delta Z_r }^2 dr & \le \normexp{ \int_t^T | f(r,X_r,\phi_{1,r},\psi_{1,r}) - 
f(r,X_r,\phi_{2,r},\psi_{2,r})| dr }^2 \\
& \le L_f^2 (T-t)^{\thetaL}  \int_t^T \{ \normexp{\delta \phi_r}^2 + \normexp{\delta \psi_r}^2 \} dr 
\end{align*}
for all $t\in [0,T)$. 
Setting $t_0 = (T - 1/(4 L_{f}^2)^{1/\thetaL} \wedge 1)\vee 0$ ensures, on the one hand, that $L_{f}^2(T-t_0)^{\thetaL} \le 1/4$, and,
on the other hand, that $T-t_0 \le 1$. Integrating the above inequality on the interval $t\in[t_0,T)$ then yields
%\begin{align*}
%&
$4 \int_{t_0}^T \{\normexp{\delta Y_r }^2 + \normexp{\delta Z_r }^2 \} dr 
\le \int_{t_0}^T \{ \normexp{\delta \phi_r}^2 + \normexp{\delta \psi_r}^2 \} dr $
and
%\\
%&
 $4\normexp{\delta Y_t}^2 \le  \int_{t_0}^T \{ \normexp{\delta \phi_r}^2 + \normexp{\delta \psi_r}^2 \} dr \quad \text{for all } t\in[t_0,T)$.
%\end{align*}
On the interval $[0,t_0)$, the function $f(t,x,\cdot)$ is Lipschitz continuous with a uniform Lipschitz constant for all $(t,x)$, 
so we proceed as in the proof of Theorem \cite[Theorem 2.1]{elka:peng:quen:97} to show that, for sufficiently large $\eta >0$,
\begin{align*}
\int_{0}^{t_0} e^{\eta r} \{\normexp{\delta Y_r }^2 + \normexp{\delta Z_r }^2 \} dr 
& \le e^{\eta t_0} \normexp{\delta Y_{t_0} }^2 +  \frac12 \int_{0}^{t_0} e^{\eta r} \{ \normexp{\delta \phi_r}^2 + \normexp{\delta \psi_r}^2 \} dr
%  \\
% & \le \frac12 \int_{0}^{T} e^{\beta_r} \{ \normexp{\delta \phi_r}^2 + \normexp{\delta \psi_r}^2 \} dr
\end{align*}
Combining this with the above estimates on $\int_{t_0}^T \{\normexp{\delta Y_r }^2 + \normexp{\delta Z_r }^2 \} dr$ 
and $\normexp{\delta Y_{t_0} }^2$ then yields
\begin{align*}
\int_{0}^{T} e^{\eta_r} \{\normexp{\delta Y_r }^2 + \normexp{\delta Z_r }^2 \} dr 
& \le   \frac12 \int_{0}^{T} e^{\eta_r} \{ \normexp{\delta \phi_r}^2 + \normexp{\delta \psi_r}^2 \} dr 
\end{align*}
where $\eta_r = \eta (r\wedge t_0 )$. This is sufficient to prove that $\Xi$ is a contraction.
\qed

\vspace{0.3cm}
 {
We now introduce an approximation procedure that will be  used repeatedly in this paper; we introduce intermediate BSDEs by ``cutting" the tail of the driver close to the time horizon $T$, prove our results for these BSDEs, then extend the result to the BSDE we're interested by limiting procedures. This technique was used extensively in \cite{gobe:makh:10}, and
we shall frequently take advantage of it throughout this work.
}
\begin{definition}
\label{def:1:intermediate BSDEs}
Let $(t,\varepsilon) \in [0,T)^2$ and define 
$\feps(t,x,y,z):=f(t,x,y,z)\1_{[0,T-\varepsilon)}(t)$. 
Let $(\Yeps,\Zeps)$ be the solution of the BSDE
\begin{equation}
  \Yeps_t = \Phi(X_T) + \int_t^T \feps(s,X_s,\Yeps_s,\Zeps_s)ds - \int_t^T \Zeps_sdW_s.
\label{eq:1:yeps:def}
\end{equation}
Additionally, let $(y,z)$ be the solution of the  BSDE with zero driver $y_t = \Phi(X_T) - \int_t^T z_s dW_s$ and 
$(\yeps,\zeps)$ the solution of the  BSDE with zero terminal condition 
\begin{equation}
 \yeps_t = \int_t^T \feps(s,X_s,y_s + \yeps_s, z_s +\zeps_s)ds 
- \int_t^T \zeps_sdW_s.
\label{eq:1:def:zeps small}
\end{equation}

\end{definition}

Since $\feps(t,x,y,z)$ is Lipschitz continuous uniformly in $t$ with Lipschitz constant $L_f\varepsilon^{(\thetaL-1)/2}$,
the solutions of the BSDEs in Definition \ref{def:1:intermediate BSDEs}
exists in $\cS^2\times\cH^2$ and are unique for all $\varepsilon \in [0,T)$ \cite[Theorem 2.1]{elka:peng:quen:97}. 
We shall also make use of the decomposition $(\Yeps,\Zeps) = (y + \yeps,z+\zeps)$, which is standard in BSDE literature \cite{gobe:makh:10}.

%\vspace{0.3cm}
%We now briefly outline some Malliavin calculus and path properties of the BSDEs $(y,z)$ and $(\yeps,\zeps)$. 
%The results are proven in \cite{gobe:makh:10} and shall be used extensively throughout this work.

\vspace{0.3cm}
We first treat the linear BSDE $(y,z)$. %This BSDE is strongly related to the linear PDE \eqref{eq:lin:PDE}.
The following Lemma relates the linear BSDE $(y,z)$ to the PDE in \eqref{eq:lin:PDE}
and gives some boundedness properties for the function $u$ and its derivatives; these bounds will be used throughout this paper.
%given in the proof of \cite[Lemma 1.1]{gobe:makh:10} and is based on a representation formula for the derivatives given in \cite[Theorem 2.11]{gobe:muno:05}.
%; we obtain stronger bounds on the partial derivatives than 
%\cite[Equation (0.11)]{gobe:makh:10} because we make stronger assumptions (uniformly bounded rather than only exponentially bounded) on the terminal 
%condition:
\begin{lemma}
\label{lem:lin:pde:bds}
Let \Hge \ be in force and consider the PDE
% and let $u : [0,T] \times \R^d \to \R$ be the solution
\begin{equation}
\label{eq:lin:PDE}
\left.
\begin{array}{l}
0 = \partial_t u + \frac12 \sum_{i,j=1}^d (\sigma \sigma^\top)_{i,j} \frac{\partial^2}{\partial_{x_i}\partial_{x_j}}u 
+ \sum_{i=1}^d b_i \frac{\partial}{\partial_{x_i}}u, \\
u(T,x) = \Phi(x).
\end{array}
\right\}
\end{equation}
Then, 
 {for all $(t,x) \in [0,T]\times\R^d$,}
\[u(t,x) = \E[\Phi(X_T)|X_t = x]\]
 is a classical solution of the PDE \eqref{eq:lin:PDE} (the so-called Feynman-Kac representation).
The derivatives
$\partial^\alpha_x u$ ($|\alpha| \le 3$), 
%$\nabla_x u$, $\nabla_x^2 u$, $\nabla_x^3 u$, 
$\partial_t u$, $\partial_t \nabla_x u$ exist and are continuous. There is a constant $C$ 
depending only on the bound on $b$ and it's derivatives, the bound on $\sigma$ and it's derivatives,  and $\elip$ such that 
 %$|\nabla_x u(t,x)| \le C \|\Phi\|_\infty (T-t)^{-1/2}$ and $|\nabla_x^2 u(t,x)| \le C \|\Phi\|_\infty (T-t)^{-1}$ for all 
\[ \|\partial^\alpha_x u(t,\cdot)\|_\infty \le C \| \Phi\|_\infty (T-t)^{-|\alpha|/2}\]
for all $(t,x) \in [0,T)\times\R^d$.
Moreover, $\big( u(t,X_t), ( \nabla_x u(t,X_t)\sigma(t,X_t))^\top \big)$ is the solution to 
the linear BSDE $(y,z)$.
For any $x_1,x_2 \in \R^d$, $t \in [0,T)$ recall from \eqref{eq:flow} the SDEs $\XM t{x_1}$ and $\XM t{x_2}$, and for  $\alpha, \beta \in [0,1]$ define $\bar X := \alpha \XM t{x_1} + \beta \XM t{x_2}$;
then 
\[|\nabla_x u(r,\bar X_r)|^2 \le  {C \E_r[ | \Phi(\bar X_T) - \E_r [\Phi(\bar X_T)]|^2] \over (T-r)}
%^{-1}
\ \text{ and } \ | \nabla_x^2 u(r,\bar X_r)|^2 \le {C \E_r[ | \Phi(\bar X_T) - \E_r [\Phi(\bar X_T)]|^2] \over (T-r)^{2}}
\]
for all $r\in[0,T)$.
\end{lemma}
{\bf Proof.} The Feynman-Kac representation of the solution is well known, see \cite{gobe:muno:05} among  others.
 { To obtain the gradient bounds, recall that  $X$ is a Markov process and denote its transition density by $p(t,x;s, \xi)$. }
 {For some $C_1$ and $\beta$ finite, the following gradient bounds hold on $p(t,x;s,\xi)$:
\begin{align*}
|\partial^\alpha_x p(t,x; s,\xi)| \le {C_1 e^{\beta|x-\xi|^2/(s-t)} \over (s-t)^{(d+|\alpha|)/2}} \quad \text{for } |\alpha| \le 3, \\
|\partial_t \partial^\alpha_x p(t,x; s,\xi)| \le {C_1 e^{\beta|x-\xi|^2/(s-t)} \over (s-t)^{(d+2+|\alpha|)/2}} \quad \text{for } |\alpha| \le 1.
\end{align*}
We obtained these bounds from \cite[Appendix A]{gobe:laba:10}, who provide references for proofs.
}
%Moreover, $X$ is a Markov process with transition density $p(t,x; \xi,s)$ and, under \Hs, the derivatives of the transition density $\partial^\alpha_x p(t,x; \xi,s)$ exist -- for all $\alpha$ with $|\alpha| \le 3$ --
%and are bounded \cite[Chapter 9, Theorem 7]{frie:64}.
 { The bounds on the derivatives of $u(t,\cdot)$ then follow from Lebesgue's differentiation theorem  (differentiation with respect to $t$ and $x$) applied to 
 \[
\partial^{\alpha_0}_t \partial^\alpha_x u(t,x) =  \partial^{\alpha_0}_t \partial^\alpha_x \int_{\R^d} \Phi(\xi) p(t,x ; T,\xi) d\xi =   \int_{\R^d} \Phi(\xi) \partial^{\alpha_0}_t \partial^\alpha_x p(t,x ; T,\xi) d\xi 
 \]
for multiindices $\alpha_0$ and $\alpha$;
we apply the gradient bounds on the transition density above and the boundedness of $\Phi$ to obtain the result on $|\partial^{\alpha_0}_t \partial^\alpha_x u(t,x) |$.
}
 
 To show the bound on $|\nabla_x u(r,\bar X_r)|$, let us recall first that the result in the case $\alpha =1$ and $\beta =0$ is given in \cite[Lemma 1.1]{gobe:makh:10}.
 The authors use the tools of \cite[Lemma 2.9]{gobe:muno:05} to show that, for every $r\in[0,T)$ and $x\in \R^d$, there is a $\cF_T$-measurable random variable $H_{r,x}$ such that 
 \[\nabla_x u(r,\XM tx_r) = \E_r[(\Phi(X_T) - \E_r[\Phi(X_T)]) H_{r,x}].\]
 { This result follows largely from the integration-by-parts formula of Malliavin calculus -- Lemma \ref{lem:mal:ibp} -- and martingale arguments; see the proof of \cite[Lemma 2.9]{gobe:muno:05} for details.}
  $H_{r,x}$ satisfies $\E_r[H_{r,x} = 0]$ and $\E_r[|H_{r,x}|^2] \le C(T-r)^{-1}$. 
   { 
  The result for $(\alpha,\beta) =(1,0)$ then follows by the Cauchy-Schwarz inequality.
  (Note that we in fact don't need \Hge \ to obtain this result.)
  }
One can follow the proof method of \cite[Lemma 2.9]{gobe:muno:05}, using additionally the linearity of the Malliavin derivative, to show that 
\[\nabla_x u(r,\bar X_r) = \E_r[ \Phi(\bar X_T) \bar H_r],\] 
where $\bar H_{r} := \alpha H_{r,x_1} + \beta H_{r,x_2}$, 
whence the result follows. The proof for the bound on $|\nabla_x^2 u(r,\bar X_r)|$ is similar.
 \qed

\vspace{0.3cm}
We move onto the non-linear BSDE $(\yeps,\zeps)$. 
The following representations and a priori estimates will be critical throughout this paper.
\begin{lemma}
\label{lem:mal:ueps veps}
 {Let \HFd \ and \Hge \ hold.
Recall the function $u : [0,T] \times \R^d \to \R$ solving the PDE \eqref{eq:lin:PDE} and  that it is differentiable  (Lemma \ref{lem:lin:pde:bds}), }
 define $\Theta_r = (r,X_r,\Yeps_r,\Zeps_r)$, and set
\begin{align}
&\aeps{r} := \nabla_x \feps(\Theta_r) +  \nabla_y \feps(\Theta_r) \nabla_x u(r,X_r) + \nabla_z \feps(\Theta_r) U(r,X_r)^\top, \nonumber \\
&\beps{r} := \nabla_y \feps(r,X_r,\Yeps_r,\Zeps_r), \quad \ceps r := \nabla_z \feps(r,X_r,\Yeps_r,\Zeps_r)
\label{eq:1:def:aeps}
\end{align}
where  {the gradients $\nabla_\xi f(\Theta_r)$ is given by $\nabla_\xi f(r,x,y,z)|_{(r,x,y,z) = \Theta_r}$ for  $\nabla_\xi f(r,x,y,z)$  defined as in Section \ref{section:notation}}, and $U(r,x)$ is defined by
% in \eqref{def:deriv du s} in Lemma \ref{lem:lin bsde:mal}.
\begin{equation}
U(t,x) := \nabla_x^2 u(t,x) \sigma(t,x) + \sum_{j=1}^d (\nabla_x u)_j (t,x) \nabla_x\sigma^\top_j(t,x), 
\label{def:deriv du s}
\end{equation}
Then there a finite constant $C$ depending only on $T$, $d$, $K^\alpha(\Phi)$, the bounds on $b$ and $\sigma$ and their derivatives, $L_f$, and $\thetaL$ such that
\begin{equation}
 \normexp{\aeps{r}} \le C\1_{[0,T-\varepsilon)}(r) (T-r)^{(\alpha + \thetaL - 3 )/2}
\label{eq:1:bds:aeps}
\end{equation}
There exists a unique solution $(\Ueps,\Veps) \in \cS^2 \times \cH^2$ of the BSDE
\begin{align}
\Ueps_t = \int_t^T \aeps{r} & + \Ueps_r \big\{ \beps{r} I_d + \nabla_x b(r,X_r) + \sum_{j=1}^q \ceps{j,r}\nabla_x\sigma_j(r,X_r) \big\} dr \nonumber  \\
&  + \int_t^T \sum_{j=1}^q (\Veps_{j,r})^\top \big\{ \ceps{j,r} I_d + \nabla_x\sigma_j(r,X_r) \big\} dr - \sum_{j=1}^q\int_t^T (\Veps_{j,r} )^\top dW_{j,r}
\label{eq:1:def:ueps veps}
\end{align}
where $\sigma_j(\cdot)$ is the $j$-th column of $\sigma(\cdot)$, $\ceps{r,j}$ is the $j$-th component of $\ceps{r}$,
and $\Veps_{j,r}$ is the $j$-th column of $\Veps_r$.
There is a (possibly different) constant $C$  such that, for any $0\le t  < T$ and $\varepsilon>0$,
% and there is a constant $C$ depending only on $L_f$, $T$, and $\thetaL$ such that
\begin{align}
\E[\sup_{t\le r <T} |\Ueps_r|^2] + \int_t^T \normexp{\Veps_r}^2 dr 
\le  C \normexp{\int_t^{T-\varepsilon} |\aeps{r}| dr }^2 \le \frac{C}{\varepsilon^{1-(\thetaL + \alpha)\wedge1}} .
\label{eq:mal:zeps:lin}
\end{align}
Let us consider $(\nabla \yeps, \nabla \zeps)$ solving the BSDE
 {
\begin{align}
\nabla \yeps_t & = \int_t^T \nabla_x \feps(\Theta_r)\nabX{r} + \nabla_y\feps(\Theta_r) \{ \nabla_x u(r,X_r) \nabX{r}
+ \nabla \yeps_r\}  dr  \nonumber \\ 
& \quad + \int_t^T \nabla_z\feps(\Theta_r)U(r,X_r)^\top \nabX{r} + \sum_{j=1}^q \nabla_z\feps_{j}(\Theta_r) (\nabla \zeps_{j,r})^\top  dr 
- \sum_{j=1}^q \int_t^T  (\nabla \zeps_{j,r})^\top dW_r.
\label{eq:grad:y}
\end{align}
}
The processes $\zeps$ and $\nabla \zeps$ satisfy the representations
\begin{align}
\zeps_t &= \Ueps_t\sigma(t,X_t)\quad m \times\P-a.e.
 \label{eq:1:zeps:repr:ueps}\\                                                
(\Veps_{j,t})^\top &=  (\nabla\zeps_{j,t})^\top \sigma^{-1}(t,X_t) - \Ueps_t \nabla_x \sigma_j(t,X_t) \; m \times\P-a.e.
\label{eq:1:veps:repr:mall zeps}
\end{align}
where $\nabla \zeps_{j,t}$ is the $j$-th column of $\nabla \zeps_{t}$.
\end{lemma}
{\bf Proof.} 
In what follows, $C$ may change from line to line.
From \cite[Lemma 1.1]{gobe:makh:10}, $\normexp{\nabla_xu(t,X_t)}\le C(T-t)^{(\alpha-1)/2}$ and 
$\normexp{\nabla_x^2u(t,X_t)} \le C(T-t)^{(\alpha-2)/2}$.
Therefore, $\normexp{\aeps r} \le C(T-r)^{(\thetaL + \alpha -3)/2}$ for all $r \in [0,T-\varepsilon]$, which is the bound \eqref{eq:1:bds:aeps}, whence
\[ (\int_0^{T-\varepsilon} \normexp{\aeps{r}}dr)^2 < \frac{C}{\varepsilon^{1-(\thetaL + \alpha)\wedge1}} <\infty.\]
This is the second inequality in \eqref{eq:mal:zeps:lin}.
Additionally, for all $t\in[0,T)$, $|\beps{t}| + \max_j |\ceps{j,t}| \le C(T-t)^{(\thetaL -1)/2}$ almost surely.
The first inequality in \eqref{eq:mal:zeps:lin} follows.
% and the existence and uniqueness of the solution now follow from Lemma \ref{lem:makhlouf:A1}.
Let $(\phi,\psi)$ be a $(\R^d)^\top \times \R^{d \times q} - $ valued process in $\cH^2 $, and define the random function
\begin{align*}
g(r,y,z) = g(r) := \aeps r & +  \phi_r((\beps{r} I_d  + \nabla_x b(r,X_r) + \sum_{j=1}^q \ceps{j,r}\nabla_x\sigma_j(r,X_r)) \\
& + \sum_{j=1}^q (\psi_{j,r})^\top (\ceps{j,r} I_d + \nabla_x\sigma_j(r,X_r)) .
\end{align*}
The function $g$ is progressively measurable and satisfies assumptions (H1)-(H5) of \cite[Section 4]{bria:dely:hu:pard:stoi:03}.
Since $f$ takes no argument in $(y,z)$, it is only necessary to validate (H1): using the triangle inequality, Jensen's inequality, the Cauchy-Schwarz inequality,
and assumptions \HFd \ and \Hs, it follows that 
%there exists a finite constant $C$ depending only on $d$, $T$, $L_f$, $\|\nabla_x b\|_\infty$ and $\max_j \|\nabla_x \sigma_j\|_\infty$ such that
\begin{align*}
\E[ ( \int_0^T |g(r)|dr )^2 ]^{1/2}  \le \E[(\int_0^T  |\aeps r| & dr )^2]^{1/2} %+ \E[\sup_{0\le r <T }  ] ^{1/2} \int_0^T \mu_r dr 
% \\&
+  C (\int_0^T\E[|\phi_r|^2 ] dr)^{1/2}( \int_0^T  \frac{dr}{(T-r)^{1 - \theta} } )^{1/2} 
\\ & \qquad 
+ C\sum_{j=1}^q (\int_0^T\E[|\psi_{j,r}|^2] dr)^{1/2}( \int_0^T  \frac{dr}{(T-r)^{1 - \theta} } )^{1/2} 
 < \infty.
\end{align*}
Thanks to \cite[Theorem 4.2]{bria:dely:hu:pard:stoi:03}, there exists a unique solution $(u,v)$  to the BSDE
\[
u_t =  \int_t^T g(r) dr - \sum_{j=1}^q \int_t^T v_{j,r} dW^j_r \qquad t \in [0,T).
\]
in $ \cS^2 \times \cH^2$. The remainder of the proof of existence and uniqueness  follows exactly as the proof of 
Theorem \ref{cor:1:exist and uni}. 
To prove the first inequality in  \eqref{eq:mal:zeps:lin}, observe that the driver $g(r)$ satisfies \eqref{eq:bria et al} from Proposition 
\ref{prop:briand:et:al} with $f_r = | \aeps r|$ and $\lambda_r = \mu_r = C(T-r)^{(\theta -1)/2}$.

The proofs of \eqref{eq:1:zeps:repr:ueps} and \eqref{eq:1:veps:repr:mall zeps} are given in \cite[Theorem 2.1]{gobe:makh:10}.
The inclussion of the local Lipschitz continuity assumptions \eqref{eq:loclip:driver} 
make no difference, because the driver $f(t,x,y,z) \1_{[0,T-\varepsilon)}(t)$
is  Lipschitz continuous uniformly in  $t$ in $(x,y,z)$ with Lipschitz coefficient $L_f \varepsilon^{(\thetaL-1)/2}$.
\qed

\subsection{A priori estimates}
\label{section:regularity}
%Recall from \eqref{eq:malweights} the Malliavin weights
For $0\le s < r \le T$, we define the Malliavin weights by
\begin{equation}
H^s_r := \frac{1}{r-s} \big( \int_s^r (\sigma^{-1}( t,X_t) 
D_sX_t% \nabla X_t \nabla X_s^{-1} \sigma(s,X_s) 
)^\top dW_t \big)^\top
\label{eq:malweights}
\end{equation}
where $D_sX_t$ is the Malliavin derivative of $X_t$ at $s$ defined in Section \ref{section:mal:SDE}. 
It was shown in Lemma \ref{lem:1:inv sigma:lip} that $|\sigma^{-1}(t,x)|$ is uniformly bounded in $(t,x)$.
% Moreover, $(D_sX_t)_{s\le t \le T}$ has a modification which
% solves a linear SDE with bounded coefficients. This implies that $\E[\sup_{s\le t\le T}|D_s X_t|^2] <\infty$.
% These properties allow us to determine the following bounds on the second moments of the Malliavin weights:
The following constant appears throughout this paper
\begin{equation}
C_M := \|\sigma^{-1}\|_\infty^2 \sup
%_{0\le s <T} 
_{s\in[0,T)}
 \sup
% _{s < t \le T} 
_{t\in(s,T]}
 \e_s[ | D_s X_t|^2]. 
\label{eq:1:constant:M}
\end{equation}
It is known from Lemma \ref{lem:mal:SDE:2} that $\sup_{ s \in [0,T)}  \sup_{t \in (s,T] } \e_s[ | D_s X_t|^2]$ is bounded.
The following result is used in the proof of \cite[Lemma 1.1]{gobe:makh:10}; 
we include it here for completeness.
\begin{lemma}
\label{lem:mal:weight}
For any $0\le s\le r\le T$,
% for all $u\le s$,
\[ \E_s[|H^s_r|^2] \le \frac{C_M}{r-s} \qquad \text{almost surely.}\]
Moreover, for every $p\ge 2$, there is a finite $C_p \ge 0$ depending only on $p$, $\|\sigma\|_\infty$, $\|\nabla_xb\|_\infty$, $\max_j \|\nabla_x \sigma_j\|_\infty$, and $T$
such that $\| H^s_r \|_p \le C_p(r-s)^{-p/2 }$.
\end{lemma}
{\bf Proof.} 
%Start with the case $s=0$. Applying It\^o's isometry yields
%\begin{align*}
%\normexp{H^s_r}^2 & = \frac{\int_s^r\normexp{\sigma^{-1}( t,X_t) D_0X_t}^2 dt}{r^2} 
%\le \frac{\|\sigma^{-1}\|_\infty^2 \E[ \sup_{0 \le r\le T}|D_0X_r|^2]}{r}
%\end{align*}
%and the proof is completed using Lemma \ref{lem:mal:SDE:3}. For $s >0$, 
Observe, using Lemma \ref{lem:mal:SDE:3} and the fact that $(s-r)^2|H^s_r|^2 -  \int_s^r |\sigma^{-1}(t,X_t) D_s X_t|^2 dt$ is a (local) martingale, that 
\[
\E_s[|H^s_r|^2] = (r-s)^{-2} \E_s[  \int_s^r |\sigma^{-1}(t,X_t) D_s X_t|^2 dt ] \le \frac{ \|\sigma^{-1}\|_\infty}{(r-s)^2} \E_s[ \int_s^r |D_s X_t|^2 dt].
\]
%Moreover, since $\nabX{t}\inX{s}   =I_d + (\int_s^t \big( \nabX{u} \inX{s} \big)^\top \cW_u^\top ) ^\top$, as shown in the proof of Lemma \ref{lem:mal:SDE:2},
%and since the partial derivatives of $b$ and $\sigma$ are uniformly bounded, it can easily be shown that $\E_s[| \nabX{t}\inX{s}|^2]$ is less than or equal to $\normS{4}{ \nabX { } }^2 \normS{4}{ \inX { } }^2$, whence the result follows. 
One then applies the conditional Fubini's lemma, Lemma \ref{lem:cond:fubini}, and the uniform bound on $\e_s[ |D_s X_t|^2 ]$ from Lemma \ref{lem:mal:SDE:3} to complete the proof.
The bound on $\| H^s_r \|_p $ is proved using the Burkholder-Davis-Gundy inequality on the continuous local martingale $(t-s)H^s_t$.
\qed

\vspace{0.3cm}
 {
The Malliavin weight is a critical element of this work. We use it to obtain a priori estimates in this section, to obtain the representation theorem in Section \ref{section:repr}, and for the Malliavin weights scheme of Section \ref{section:num}. The following elementary corollary indicates an important technique in which we make use of the Cauchy-Schwarz inequality in conditional form in order to obtain upper bounds:
}
 {
\begin{corollary}
\label{cor:cs:ineq}
Let $G \in \L_2(\cF_T)$ and $g \in \L_2([0,T]\times\Omega)$. Then 
\[
| \E_t[G H^t_T] | \le {\sqrt{C_M} (\E_t[|G|^2])^{1/2} \over \sqrt{T-t}} \quad \text{and} \quad |\E_t[\int_t^T \{ g_s H^t_s \}  ds] | \le \sqrt{C_M} \int_t^T  { \hg_s \over \sqrt{s-t}}  ds
\]
\end{corollary}
where $\hg \in  \L_2([0,T]\times\Omega)$ is a version of $\big( (\E_t[|g_s|^2])^{1/2} \big)_{s \in [0,T)}$.
}
 {
\begin{remark*}
We leave the implementation of the conditional Fubini theorem, Lemma \ref{lem:cond:fubini}, in its full form in the above lemma, without using the notation given in Section \ref{section:notation}.
We do this to be absolutely clear about how the conditional Fubini theorem is used in this paper, before returning to the -- in our opinion -- much more clear, if slightly abusive, notation $\int_t^T (\e_t[g_s])^{1/2} ds$.
\end{remark*}
}

 {
{\bf Proof.}
The first inequality follows from application of the conditional Cauchy-Schwarz inequality $|\E_t[G H^t_T]| \le (\E_t[|G|^2])^{1/2} (\E_t[|H^t_T|])^{1/2} $, then using Lemma \ref{lem:mal:weight} to upper bound the conditional expectation $(\E_t[|H^t_T|])^{1/2}$.
The second inequality is a little more intricate to obtain due to the Lebesgue integral.
First, apply the conditional Fubini theorem, Lemma \ref{lem:cond:fubini}, to obtain
\[
\E_t[ \int_t^T | g_s H^t_s | ds  ] = \int_t^T \hH_s ds  
\]
where $\hH \in  \L_2([0,T]\times\Omega)$ is a version of $(\e_t[|g_s H^t_s|])_{s \in[0,T)}$.
Now, applying the conditional Cauchy-Schwarz inequality and Lemma \ref{lem:mal:weight} to $\e_t[|g_s H^t_s|]$, it follows that 
\[
\hH_s \le C_M{ \hg_s \over \sqrt{s-t} } \quad \text{for almost all } s\in[0,T) \quad \P-\text{a.s.} 
\]
as required. \qed
}

We now state and prove a priori results on the solutions of BSDEs with drivers satisfying \eqref{eq:loclip:driver}.
 {
These estimates are in the spirit of \cite[Proposition 2.1]{elka:peng:quen:97} with two extensions: 
firstly, we allow the drivers of the BSDEs to satisfy locally Lipschitz continuity  like condition \Hf;
secondly, we prove point-wise (in time) a priori estimates on the $Z$ processes assuming the existence of a representation formula.
The latter estimates will be extremely useful, as we shall prove the this representation formula for our BSDEs in Section \ref{section:repr} and use the below proposition extensively in subsequent sections.
}

\begin{proposition}
\label{prop:1:apriori}
Let $x\mapsto\Phi_{1},\Phi_{2} \in \L_2(\cF_T)$% be measurable functions such that $\big(\Phi_1(X_T),\Phi_2(X_T)\big) \in \big(\L_2(\cF_T)\big)^2$.
%in $\L_{2,\alpha}$. 
%Additionally, let
 and $(\omega , t,y,z) \mapsto f_{1}(\omega, t,y,z),f_{2}(\omega, t,y,z)$ be 
$\cP\otimes\cB(\R) \otimes\cB( (\R^q)^\top)$-measurable functions
%  satisfying \eqref{eq:loclip:driver} 
for which there are constants $(\theta_{1,L},\theta_{2,L}) \in (0,1]^2$ and   $(L_{f_{1}},L_{f_{2}}) \in (0,\infty)^2$ such that 
\[
|f_i(\omega,t,y,z) - f_i(\omega,t,y',z')| \le \frac{L_{f_i}\{ |y - y'| + |z - z'| \} }{(T-t)^{(1-\theta_{i,L})/2} }  \quad m \times\P-\text{almost everywhere,}
\]
%   with parameters $L_{f_{i}}$, $C_{f_{i}}$, $\theta_{i,L}$
and $f_i(\omega,t,0,0) \in \cH^2$ for $i \in\{ 1,2\}$.
Let $(Y_{i},Z_{i})$ be a solution to the FBSDE with terminal condition
$\Phi_{i}$ and driver $f_{i}(t,y,z)$ ($i=1,2$ respectively). 

Define 
\begin{align*}
&\Delta Y_{t} := Y_{1,t} -  Y_{2,t}, \qquad \Delta Z_{t} := Z_{1,t} -  Z_{2,t}, \\
& \Delta f_{t} :=  f_{1}(t,Y_{1,t},Z_{1,t}) -  f_{2} (t,Y_{1,t},Z_{1,t}),
%\\ & 
\quad \Delta \Phi := \Phi_{1} - \Phi_{2}.
\end{align*}
Then there is a finite constant  $C \ge 0$ depending only on $T$, $L_{f_{2}}$ and $\theta_{2,L}$
 such that, for all $s<t < T$,
\begin{align}
& \E_s[\Delta Y_t^2] + \E_s[ \int_t^T |\Delta Z_s|^2 ds]
 \le C \E_s[\Delta\Phi^2] + C \big(\int_{t}^T \E_s[\Delta f_r^2]^{1/2} dr\big)^2
\label{eq:1:apriori:1} 
%\label{eq:1:apriori:3} 
\end{align}
Moreover, suppose that
$Z_{i,t} := \E_t[\Phi_i(X_T) H^t_T + \int_t^T f_i(r,X_r,Y_{i,r},Z_{i,r}) H^t_r dr]$ for all $t\in[0,T)$ almost surely  ($i=1,2$).
Then there is a (possibly different) finite constant $C\ge 0$ depending only on $T$, $C_M$, $L_{f_{2}}$,  and $\theta_{2,L}$ such that,
\begin{align}
(\E_s[|\Delta Z_t|^2])^{1/2} \le C\frac{(\E_s\big[(\Delta\Phi - \E_t\Delta\Phi)^2])^{1/2} }{\sqrt{T-t}} 
+ C \int_t^T \frac{(\E_s[\Delta f_r^2])^{1/2} }{\sqrt{r-t}} dr 
 + C (\E_s[\Delta\Phi^2])^{1/2}(T-t)^{\thetaL/2} 
\label{eq:1:apriori:2} 
%\label{eq:1:apriori:4} 
\end{align}
for all $t\in[0,T)$ almost surely.
\end{proposition}

{\bf Proof.} 
In what follows, $C$
 may change from line to line. 
 We start by proving the result for $s = 0$; the general case is proved analogously,
the only difference is that one must use the conditional version of the Minkowski, Cauchy-Schwarz  {(Corollary \ref{cor:cs:ineq})}, and H\"older inequalities  {in the place of the usual version of these with the regular expectation}.
%Lemma \ref{lem:1:apriori:cond} is particularly useful when the terminal condition is bounded and H\"older continuous, because
%this allows us to make almost sure estimates on $|Z_t|$.
Using the definition of the BSDE \eqref{eq:1:BSDE},
\begin{align*}
\Delta Y_{t} + \int_{t}^T \Delta Z_{s} dW_{s} & = \Delta\Phi + \int_{t}^T \Delta f_{s } ds 
%\\& 
+ \int_{t}^T f_{2} (s,Y_{1,s},Z_{1,s}) -f_{2} (s,Y_{2,s},Z_{2,s}) ds.
\end{align*}
% Set $\beta_s := C(T\vee1) \int_0^s \frac{dr}{(T-r)^{1-\thetaL} }  = C(T^{\thetaL} - (T-s)^{\thetaL})  \ge 0$.
Using \eqref{eq:loclip:driver} and H\"older's inequality,
\begin{align}
\normexp{\Delta Y_{t}}^2 + \int_{t}^T& \normexp{\Delta Z_{s}}^2 ds
% \nonumber \\ & 
 \le 3\normexp{\Delta\Phi}^2 + 3 \normexp{ \int_{t}^T \Delta f_{s } ds }^2 
%\nonumber\\&
+ 3  \normexp{ \int_{t}^T f_{2} (s,Y_{1,s},Z_{1,s}) -f_{2} (s,Y_{2,s},Z_{2,s}) ds }^2 
\nonumber\\& 
\le 3 \normexp{\Delta\Phi}^2 + 3\normexp{ \int_{t}^T| \Delta f_{s } | ds }^2
+ 3L_{f_2}^2 \normexp{ \int_{t}^T \frac{|\Delta Y_{s}| + |\Delta Z_{s}|}{(T-s)^{(1-\theta_{2,L})/2 } } ds }^2 
\nonumber\\
& \le 3 \normexp{\Delta\Phi}^2 + 3 \normexp{ \int_{t}^T |\Delta f_{s } |ds }^2 
%\nonumber \\ & 
+ 3L_{f_2}^2 \int_{t}^T \frac{ 1 }{(T-s)^{1-\theta_{2,L}} } ds 
\int_{t}^T \{\normexp{\Delta Y_{s}}^2 + \normexp{\Delta Z_{s}}^2 \}ds \nonumber \\
& \le 3 \normexp{\Delta\Phi}^2 + 3 \normexp{ \int_{t}^T |\Delta f_{s } |ds }^2
+ 3L_{f_2}^2 (T-t)^{\theta_{2,L}} \int_{t}^T \{\normexp{\Delta Y_{s}}^2 +\normexp{\Delta Z_{s}}^2 \}ds 
\label{eq:1:apriori:y}
\end{align}
Setting $t_0 = (T - 1/(6L_{f_2}^2)^{1/\theta_{2,L}} )\vee 0$ ensures that $3L_{f_2}(T-t_0)^{\theta_{2,L}} \le 1/2$, and,
on the other hand, that $T-t_0 \le 1$.
Integrating \eqref{eq:1:apriori:y} over $(t_0,T)$, we obtain
\begin{align}
\label{eq:1:apriori:integrals}
 \int_{t_0}^T \normexp{\Delta Y_{t}}^2 + \normexp{\Delta Z_{s}}^2 ds 
\le 6 \normexp{\Delta\Phi}^2 + 6 L_{f_2}^2\normexp{ \int_{t_0}^T |\Delta f_{s } |ds }^2
\end{align}
Substituting \eqref{eq:1:apriori:integrals} into \eqref{eq:1:apriori:y} then yields
\[
 \sup_{t_0 \le t < T} \normexp{\Delta Y_{t}}^2 \le 6 \normexp{\Delta\Phi}^2 + 6 \normexp{ \int_{t_0}^T |\Delta f_{s } |ds }^2
\]
and this gives the result in the interval $[t_0,T]$.  

In the interval $[0,t_0)$, the function $(y,z) \mapsto f_2(\omega,t,y,z)$ is $m\times \P$ Lipschitz continuous with Lipschitz constant 
$ { \tilde L} := L_f (T-t_0)^{(\theta_{2,L}-1)/2}$. It then follows from \cite[Proposition 2.1]{elka:peng:quen:97} that
\[
 \sup_{0 \le t < t_0}  \normexp{\Delta Y_{t}}^2 + \int_{t_0}^T \normexp{\Delta Z_{s}}^2 ds 
\le C \normexp{\Delta Y_{t_0}}^2 + C \normexp{ \int_0^{t_0} |\Delta f_{s } |ds }^2
\]
and the proof of \eqref{eq:1:apriori:1} is complete by substituting the bounds on $\normexp{\Delta Y_{t_0}}^2$ from above.

% The proof of \eqref{eq:1:apriori:1} is completed by pasting. 
% Since the interval length is always $T-t_0$, there will only be a finite number
% $T(C/2)^{1/(1+\theta_{2,L})}$ of pasting steps, so there are no subtle arguments to be made here. As usual, the dependence on time and $L_{f_2}$
% of the final estimate is exponential.

\vspace{0.3cm}
Next, we prove  \eqref{eq:1:apriori:2}. 
 {Recall that $\E_t[H^t_s] = 0$ for all $(t,s)$,
which implies that $\E_t[\Phi_i H^t_T] = \E_t[(\Phi_i - \E_t[\Phi_i] ) H^t_T] $.}
Using the representation $Z_{i,t} :=
\E_t[\Phi_i H^t_T + \int_t^T f_i(r,Y_{i,r},Z_{i,r}) H^t_r dr]$, 
it follows from Minkowski's inequality, the Cauchy-Schwarz inequality  {(i.e. Corollary \ref{cor:cs:ineq})}, and Lemma \ref{lem:mal:weight} that 
\begin{align}
\normexp{\Delta Z_t} & \le \frac{C V_{t,T}(\Delta\Phi)}{\sqrt{T-t} } + C \int_t^T \frac{\normexp{\Delta f_r}}{\sqrt{r-t}}dr 
+C \int_t^T \frac{\normexp{\Delta Y_r} + \normexp{\Delta Z_r}}{(T-r)^{(1-\theta_{2,L})/2}\sqrt{r-t}} dr.
\label{eq:1:apriori:z:1}
\end{align}
where we define $V_{t,T}(\Delta\Phi)$ by $\E[|\Delta\Phi - \E_t[\Delta\Phi]|^2]^{1/2}$. Defining $\Theta_r := \normexp{\Delta Y_r} + \normexp{\Delta Z_r}$ and recalling \eqref{eq:1:apriori:y}, it follows that
\begin{align}
\Theta_t & \le C\normexp{\Delta\Phi} + \frac{C V_{t,T}(\Delta\Phi)}{\sqrt{T-t} } + C \int_t^T \frac{\normexp{\Delta f_r}}{\sqrt{r-t}}dr 
+C \int_t^T \frac{\Theta_r}{(T-r)^{(1-\theta_{2,L})/2}\sqrt{r-t}} dr.
\label{eq:1:apriori:theta:1}
\end{align}
Applying Lemma \ref{lem:iteration:2} with $u_t := \Theta_t$ and 
\[
w_t :=  
 C\normexp{\Delta\Phi} + \frac{C V_{t,T}(\Delta\Phi)}{\sqrt{T-t} } + C \int_t^T \frac{\normexp{\Delta f_r}}{\sqrt{r-t}}dr,
\]
it follows that
\begin{align*}
\Theta_r & \le  C w_t + C \int_t^T \frac{w_r  }{(T-r)^{(1-\theta_{2,L})/2}\sqrt{r-t}} dr + 
C \int_t^T \frac{\Theta_r}{(T-r)^{(1-\theta_{2,L})/2}} dr
\end{align*}
whence it follows from Lemma \ref{lem:iteration:3} that 
\[
\int_t^T \frac{\Theta_r}{(T-r)^{(1-\theta_{2,L})/2}\sqrt{r-t}} dr \le C \int_t^T \frac{w_r  }{(T-r)^{(1-\theta_{2,L})/2}\sqrt{r-t}} dr
\]
Substituting this into \eqref{eq:1:apriori:z:1} and applying Lemma \ref{lem:integration:1} leads to
\begin{align*}
\normexp{\Delta Z_t} & \le \frac{C V_{t,T}(\Delta\Phi)}{\sqrt{T-t} } + C \int_t^T \frac{\normexp{\Delta f_r}}{\sqrt{r-t}}dr 
+ C\int_t^T \frac{w_r  }{(T-r)^{(1-\theta_{2,L})/2}\sqrt{r-t}} dr \\
& = \frac{C V_{t,T}(\Delta\Phi)}{\sqrt{T-t} } + C \int_t^T \frac{\normexp{\Delta f_r}}{\sqrt{r-t}}dr 
+ C\int_t^T \frac{V_{r,T}(\Delta\Phi)  }{(T-r)^{(2-\theta_{2,L})/2}\sqrt{r-t}} dr \\
& \quad + C\int_t^T \frac{\int_r^T\normexp{\Delta f_s}(s-r)^{-1/2} ds }{(T-r)^{(1-\theta_{2,L})/2}\sqrt{r-t}} dr 
+ C\normexp{\Delta\Phi}(T-t)^{\theta_{2,L}/2} \\
& = \frac{C V_{t,T}(\Delta\Phi)}{\sqrt{T-t} } + C \int_t^T \frac{\normexp{\Delta f_r}}{\sqrt{r-t}}dr 
+ C\int_t^T \frac{V_{r,T}(\Delta\Phi)  }{(T-r)^{(2-\theta_{2,L})/2}\sqrt{r-t}} dr \\
& \quad + C\int_t^T \normexp{\Delta f_s} \{\int_r^s(s-r)^{-1 +\theta_{2,L}}(r-t)^{-1/2} dr \} ds 
+ C\normexp{\Delta\Phi}(T-t)^{\theta_{2,L}/2} .
\end{align*}
The proof is completed by observing that $V_{r,T}(\Delta\Phi)$ is non-increasing in $r$.
\qed

\vspace{0.3cm}
%One can adapt in a straightforward manner the proof of \cite[Theorem 2.1]{elka:peng:quen:97} with the methods used in the proof of Proposition \ref{prop:1:apriori} to obtain the existence and
%uniqueness of the locally Lipschitz BSDE.

%\vspace{0.3cm}
The estimates \eqref{eq:1:apriori:2} allow us to determine a priori estimates on the conditional second moments of the solution of the BSDE $(Y,Z)$.
\begin{corollary}%[Moment bounds]
\label{cor:mom bd}
Assume that 
%$f$ satisfies \eqref{eq:loclip:driver},  and 
$Z_t = \E_t[\Phi(X_T) H^t_T + \int_t^T f(r,X_r,Y_r,Z_r) H^t_r dr]$ for all $t\in[0,T)$
almost surely. 
Then there is a constant $C$ depending only on $L_f$, $\thetaL$, $C_f$, $\thetaC$, $K^\alpha(\Phi)$ and $T$ 
 such that, for all $t\in[0,T)$ and $s \in [0,t]$, we have
\begin{align*}
& \sup_{s\le t \le T} (\E_s[|Y_t|^2])^{1/2}  \le C( 1 +  (\E_s[|\Phi(X_T) - \E_s[\Phi(X_T)]|^2])^{1/2}) , 
%\label{eq:2:mom bounds}
\\ 
& (\E_s[|Z_t|^2])^{1/2}  \le \frac{C(\E_s[|\Phi(X_T) - \E_s[\Phi(X_T)]|^2])^{1/2}}{ \sqrt{T-t} } + \frac{C}{(T-t)^{(1-2\thetaC)/2} } 
 + C(\E_s[|\Phi(X_T) |^2])^{1/2} (T-t)^{\thetaL/2}. % \quad dt-a.e.
%\label{eq:2:mom bounds}
\end{align*}
In particular, $\normexp{Y_t} \le C$ and $\normexp{Z_t} \le C(T-s)^{((2\thetaC)\wedge\alpha -1)/2}$ for all $t \in [0,T)$, and
%Moreover, for all $t\in[0,T)$,
\begin{equation}
\normexp{f(s,X_s,Y_s,Z_s)} \le \frac{C}{(T-s)^{1-((2\thetaC)\wedge\alpha + \thetaL)/2}} +\frac{C}{(T-s)^{1-\thetaC}}.% \quad ds-a.e.
\label{eq:bd:f} 
\end{equation}
If \Hgh \ is in force, we have additionally that $|Z_t| \le CK_\Phi(T-s)^{((2\thetaC)\wedge\thetaP -1)/2}$ for all $t\in[0,T)$ almost surely.
\end{corollary}
%\begin{remark}
%\label{rem:holder:thetaL}
%If \Hgh \ is in force, one sees from the above absolute bounds that $(Y,Z)$ solves the BSDE with terminal condition $\Phi(X_T)$
%and driver $f(t,X_t, y, \cT_{\thetaP} (z) )$, where $\cT_{\thetaP}(z):=  CK_\Phi(T-s)^{((2\thetaC)\wedge\thetaP -1)/2}$
%\end{remark}
{\bf Proof.}
In what follows, $C$  may change from line to line.
As in Proposition \ref{prop:1:apriori}, we only prove the result for $s =0$;
 {the general case is proved using the conditional version of the Minkowski, Cauchy-Schwarz  {(Corollary \ref{cor:cs:ineq})}, and H\"older inequalities  {in the place of the usual version of these with the regular expectation}.}
 {Recalling $V_{t,T}(\Phi)$ from \Hg,}
apply \eqref{eq:1:apriori:2} from Proposition \ref{prop:1:apriori} with $(Y_1,Z_1) := (0,0)$ and $(Y_2,Z_2) := (Y,Z)$ to obtain
(for all $t\in[0,T)$)
\begin{align*}
\normexp{ Z_t} &\le C\frac{V_{t,T}(\Phi) }{\sqrt{T-t}} 
+ C \int_t^T \frac{\normexp{ f(r,X_r,0,0)} }{\sqrt{r-t}} dr 
 + C \normexp{\Phi}(T-t)^{\thetaL/2} \\
&\le \frac{C}{(T-t)^{(1-\alpha)/2} } + C \int_t^T \frac{dr }{(T-r)^{1-\thetaC}\sqrt{r-t}} + C (T-t)^{\thetaL/2}.
\end{align*}
%and the \eqref{eq:1:mom bounds} follows. 
Combining the local Lipschitz continuity and boundedness of $f$ in \eqref{eq:loclip:driver} leads to the required bound on the conditional second moments of $Z_t$.
The estimate on the conditional moments of $Y_t$ is obtained similarly starting from \eqref{eq:1:apriori:1}.
The remaining bounds are obtained by taking into account  \eqref{eq:loclip:driver} and the regularity of the terminal condition (\Hg \ or \Hgh).
%and the bounds in \eqref{eq:1:mom bounds}
%leads to \eqref{eq:bd:f}.
\qed

\vspace{0.3cm}
Recall $(\Yeps,\Zeps)$ from Definition \ref{def:1:intermediate BSDEs} in Section \ref{section:mal:BSDE}, 
the BSDE with terminal condition $\Phi$ and driver  $\feps(t,x,y,z) := 
f(t,x,y,z)\1_{[0,T-\varepsilon)}(t)$.
The following corollary of Proposition \ref{prop:1:apriori} will be used extensively throughout this paper;
it provides a stability results between the BSDEs $(Y,Z)$ and $(\Yeps,\Zeps)$ that are controlled by $\varepsilon$.
\begin{corollary}
\label{cor:cut:error}
Let $\gamma := (\thetaC\wedge\frac\alpha2 + \frac\thetaL2)\wedge\thetaC $ and assume that 
$Z = \e_t[\Phi(X_T) H^{t}_T + \int_{t}^T  f(s,X_s,Y_s,Z_s) H^{t}_s ds] $ and 
$\Zeps_t = \e_t[\Phi(X_T) H^{t}_T + \int_{t}^T \feps(s,X_s,\Yeps_s,\Zeps_s) H^{t}_s ds] $
for all $t\in[0,T)$ almost surely.
Then there is a constant $C$  such that
\begin{align}
&\sup_{0\le t\le T}\normexp{Y_t - \Yeps_t}^2 + \int_0^T\normexp{Z_t - \Zeps_t}^2 dt \le C\varepsilon^{2\gamma}, \label{eq:cut:error:1}\\
& \normexp{Z_t - \Zeps_t} \le C\int_{t\vee(T-\varepsilon)}^T \frac{ds }{(T-s)^{1-\gamma}\sqrt{s-t}}  \label{eq:cut:error:2}
\end{align}
for all $t\in[0,T)$.
In particular, $(\Yeps,\Zeps) \rightarrow (Y,Z)$ as $\varepsilon \rightarrow 0$ in $\cS^2\times\cH^2$.
\end{corollary}
{\bf Proof.}
In what follows, $C$  may change from line to line.

It follows from \eqref{eq:1:apriori:1} in Proposition \ref{prop:1:apriori} that
\begin{align}
\sup_{0\le t\le T}\normexp{Y_t - \Yeps_t}^2 + \int_0^T \normexp{Z_s - \Zeps_s}^2 ds & 
\le C \Big(\int_{T-\varepsilon}^T \normexp{f(s,X_s,Y_s,Z_s) } ds\Big)^2.
\label{eq:1:Z- Zeps}
\end{align}
Substituting \eqref{eq:bd:f} into \eqref{eq:1:Z- Zeps} 
combined with $\Big(\int_{T-\varepsilon}^T \frac{ds}{(T-s)^{(1-\gamma)}}\Big)^2 \le C\varepsilon^{2\gamma}$ completes the proof
of \eqref{eq:cut:error:1}.
Next, it follows from \eqref{eq:1:apriori:2}  that
\begin{align*}
\normexp{Z_t - \Zeps_t} & \le C\int_{t\vee(T-\varepsilon)}^T \frac{\normexp{f(s,X_s,Y_s,Z_s)} }{\sqrt{s-t}} ds \quad 
\text{for all } t\in[0,T).
% & \le C\int_{t\vee(T-\varepsilon)}^T \frac{ds }{(T-s)^{1-\gamma}\sqrt{s-t}} 
% \label{eq:cut:error:2}
\end{align*}
Substituting \eqref{eq:bd:f} above proves \eqref{eq:cut:error:2}.
\qed

\vspace{0.3cm}
To end this section, we present a mollification procedure that will be used frequently to allow us to extend results under
the assumptions \HFd \ and \Hge \ to the same results without these assumptions. The following corollary is a trivial consequence of Proposition \ref{prop:1:apriori} and the properties of mollifiers.

\begin{corollary}
\label{cor:1:moll}
%Let $(Y,Z)$ be the solution of the BSDE with terminal condition $\Phi(X_T)$ and driver $f(t,X_t,y,z)$. 
Let $M>0$ be finite, and $M \mapsto R(M) \ge 1$ be 
%finite for all $M\in  \R_+$ and 
increasing w.r.t. $M$. Define  $\Phi_M(x) := -M\vee \Phi(x) \wedge M$  and,
recalling 
%the function $\phi_{R,M}$ from
the mollifier $\phi$ of Definition \ref{def:1:moll},  
\[f_M(t,x,y,z) :=\int_{\R^d\times\R\times (\R^q)^\top} f(t,x - x',y - y', z - z')
\phi_{R(M)} ( x',y',z' ) d(x',y',z'),\]
Let $(Y_M,Z_M)$ be the solution of the BSDE with terminal condition $\Phi_M$ and driver $f_M(t,x,y,z)$.
Then $\Phi_M$ satisfies \Hge, $f_M$ satisfies \HFd, and $(Y_M,Z_M) \rightarrow (Y,Z)$ as $M\rightarrow \infty$ in $ \cS^2\times\cH^2$.
% 
% Suppose that $Z$ satisfies the representation \eqref{eq:1:Z} and $Z_M$ satisfies 
% \[
%  Z_{M,t} = \E_t[\Phi_M(X_T)H^t_T + \int_t^T f_M(r,X_r,Y_{M,r}, Z_{M,r}) H^t_r dr] \quad dt \times \P-a.e.
% \]
% Then
\end{corollary}

\subsection{Representation theorem}
\label{section:repr}

 {In this section, we prove that BSDEs satisfying the local Lipschitz continuity and local  boundedness conditions \Hf \ also satisfy the a representation theorem in the spirit of \cite[Theorem 3.1]{ma:zhan:02}.
Following on from Section \ref{section:regularity}, we see that this representation is very valuable, as it gives us additional access to a priori results.
We use these a priori results in the sections that follow, so it is essential that we also establish the representation result.
Unlike in the proof of \cite[Theorem 3.1]{ma:zhan:02}, we do not prove the representation result on $Z$ directly.
The strategy is rather to take the approximative BSDE $(\Yeps,\Zeps)$, for which we  already know that $\Zeps$ satisfies the representation from \cite[Theorem 3.1]{ma:zhan:02}, then to prove it converges in $\cH^2$ to the process that we claim is a version of $Z$ as $\varepsilon$ converges to $0$ by classical $(\varepsilon,\delta)-$arguments, and to finally conclude using the fact that $\Zeps$ also converges to $Z$ in $\cH^2$ and because $Z$ is unique.
}

%To conclude this section, we state and prove the representation theorem for locally Lipschity continuous and locally bounded BSDEs.
\begin{theorem}
\label{thm:repr}
 {Recall $\L_{2,\alpha}$ from \Hg,}  suppose that $\Phi \in \L_{2,\alpha}$ and $(t,x,y,z) \mapsto f(t,x,y,z)$ satisfies \Hf. % \eqref{eq:loclip:driver}.
 Then, there is a predictable version $\cZ$ of $Z$ which satisfies  
\begin{equation}
\cZ_t = \e_t[\Phi(X_T) H^{t}_T + \int_{t}^T f(s,X_s,Y_s,Z_s) H^{t}_s ds] \quad 
\text{for all } t\in[0,T) \quad \P-a.s.
\label{eq:1:Z}
\end{equation}
where $H^t_s$ are the Malliavin weights given in  \eqref{eq:malweights}.
\end{theorem}
%From now on, we work with the the version $\cZ$ of $Z$.

% \begin{theorem}
% \label{thm:representation:unif}
% Assume that $\Phi \in \L_{2,\alpha}$ and $f$ satisfies the uniform Lipschitz continuous conditions
% \eqref{eq:loclip:driver} with $\thetaL =1$. 
% Moreover, assume that $(t,x)\mapsto b(t,x), \; \sigma(t,x)$ are bounded and twice continuously differentiable with bounded,
% $\gamma$-H\"older continuous partial derivatives in $x$, and are $\frac12$-H\"older continuous in time.
% Then $(Z_t)_{0\le t<T}$ has the representation
% \begin{equation*}
% Z_t = \e_t[\Phi(X_T) H^{t}_T + \int_{t}^T f(s,X_s,Y_s,Z_s) H^{t}_s ds] \quad dt\times\P-a.e.
% \end{equation*}
% for $H^s_r$ given by \eqref{eq:malweights}. Moreover, $Z$ has a continuous version.
% \end{theorem}
{\bf Proof.} 
%To prove the representation formula for $Z$, we will work in dimension 1 to keep notation simple.
In the following,   $C$ is a constant whose value may change from line to line.

\vspace{0.3cm}
To start with, 
let assume \HFd \ and \Hge \ be in force.
 {We prove the representation theorem first under these conditions, and then extend to the general result by means of mollification.}
Recall the BSDEs  $(\Yeps,\Zeps)$, $(y,z)$ and $(\yeps,\zeps)$ from  Section \ref{section:mal:BSDE}, and  the decomposition 
$(\Yeps,\Zeps) = (y + \yeps,z + \zeps)$.
 {
We first prove the that there is a predictable version of $\Zeps$ equalling
\begin{equation}
\E_t[\Phi(X_T)H^t_T + \int_t^T \feps(r,X_r, \Yeps_r,\Zeps_r) H^t_r dr] \quad \text{for all } t\in[0,T) \quad \P-a.s. 
\label{eq:repr:eps}
\end{equation}
In fact, this is an application of  \cite[Theorem 4.2]{ma:zhan:02}; this is not immediately clear, so we make the calculations explicit for the benefit of the reader.
}
%We obtain the representation theorem for the BSDEs $(\yeps,\zeps)$ and use the decomposition to obtain representation theorem for $(\Yeps,\Zeps)$;
%the representation theorem for $(y,z)$ is well known, see the proof of \cite[Lemma 1.1]{gobe:makh:10}.
%Firstly, 
%recall from 
 {
Definition \ref{def:1:intermediate BSDEs} and 
Lemma
\ref{lem:lin:pde:bds} give us that $(\yeps,\zeps)$ solves the BSDE with terminal condition $0$ and driver 
\[F(t,x,y,z) := \feps(t,x,u(t,x) + y, \nabla_xu(t,x) \sigma(t,x)  + z )\]
on the time interval $[0,T-\varepsilon]$.
%on this interval,
Due to the bounds on $u$ and its derivatives given in Lemma \ref{lem:lin:pde:bds}, the Lipschitz constant of $(x,y,z) \mapsto F(t,x,y,z)$ is bounded from above (for all $t\in [0,T-\varepsilon]$) by
\[
{L_f \over (T-t)^{(1-\thetaL)/2}} \{1 + \|\nabla_x u(t,\cdot)\|_\infty + \|\nabla_x^2 u(t,\cdot)\|_\infty \} \le {C \over (T-t)^{(3-\thetaL)/2} } \le C\varepsilon^{-(3-\thetaL)/2} =: L_F.
\]
Using this Lipschitz constant, we also show that $F(t,x,0,0)$ is bounded (for all $(t,x) \in [0,T-\varepsilon] \times\R^d$ by
\[
{C_f \over (T-t)^{1-\thetaC} } + L_F \{1 + \|u(t,\cdot)\|_\infty + \|\nabla_x u(t,\cdot)\|_\infty\} \le C_f \varepsilon^{-(1-\thetaC)} + L_F\varepsilon^{-(1 - \thetaL)/2} =: C_F.
\]
Therefore, the driver $F$ is uniformly Lipschitz continuous in $(x,y,z)$ and uniformly bounded at $(y,z) = (0,0)$,
i.e. it satisfies \Hf \ with $\theta_{L,F} \equiv 1$, $\theta_{C,F} \equiv 1$, and constants $L_F$ and $C_F$ (given above). 
$F$ is also continuous in $t$.}
Therefore, \cite[Theorem 4.2]{ma:zhan:02} applies to the BSDE in the interval $[0,T-\varepsilon]$, i.e. there is a version of $\zeps$ equalling
\[
 \e_t[\int_t^{T-\varepsilon} F(r,X_r,  \yeps_r, \zeps_r) H^t_r dr] \ \text{ for all } t\in[0,T-\varepsilon] \ \text{ almost surely.}
 \]
%for all $t \in [0,T-\varepsilon]$ almost surely. 
On the other hand, $\zeps_t $ and $F(t,x,y,z)$ are $0$ for all $t\in (T-\varepsilon, T]$ almost surely, so the representation holds trivially in the interval $(T-\varepsilon,T]$, whence it follows that there is a version of $\zeps$ equalling
\[
 \e_t[\int_t^{T-\varepsilon} F(r,X_r,  \yeps_r, \zeps_r) H^t_r dr] \ \text{ for all } t\in[0,T] \ \text{ almost surely.}
\]
 {Now, it is well known 
-- see for example \cite[Page 1116]{gobe:makh:10}, where our $H^t_T$ is given by $H^{(1)}_{t,T} \sigma(t,X_t)$ in their notation -- 
that there is predictable version of $(z_t)_{t\in[0,T)}$ equalling 
\[
\E_t[\Phi(X_T) H^t_T] \ \text{ for all } t\in[0,T) \ \text{ almost surely,}
\]
and this implies the version of $\Zeps$ given by \eqref{eq:repr:eps} thanks to the the decomposition  $(\Yeps,\Zeps) = (y + \yeps,z + \zeps)$.
}
%Therefore, including the representation of $(y,z)$, it follows that % since $z_t = \E_t[\Phi(X_T) H^t_T]\quad m \times\P-a.e.$, we have that

%To complete the proof for $(\Yeps,\Zeps)$, define by $\cZ^{(\varepsilon)}$ the predictable projection
%\cite[Theorem 2.28]{jaco:shir:03}
% of the process $ (\cX^{(\varepsilon)}_t := \Phi(X_T)H^t_T + \int_t^T \feps(r,X_r, \Yeps_r,\Zeps_r) H^t_r dr)_{t\in[0,T)}$, 
% and observe from \eqref{eq:repr:eps} that $\Zeps_t = \cZ^{(\varepsilon)}_t$ $m \times \P$-almost everywhere.

Define by $\cZ$ the predictable projection \cite[Theorem 2.28]{jaco:shir:03}
of the process $(\cX_t := \Phi(X_T)H^t_T + \int_t^T f(r,X_r, Y_r,Z_r) H^t_r dr)_{t\in[0,T)}$.
%We take the version of $\Zeps$ given by the the predictable representation of the process  $ (\cX^{(\varepsilon)}_t := \Phi(X_T)H^t_T + \int_t^T \feps(r,X_r, \Yeps_r,\Zeps_r) H^t_r dr)_{t\in[0,T)}$.
In what follows, we show that $\normexp{\Zeps_t - \cZ_t} \rightarrow 0$ as $\varepsilon\rightarrow0$ for almost all $t\in[0,T)$.
This implies, by the dominated convergence theorem, that $\Zeps \rightarrow \cZ$ in $\cH^2$.
Since $\Zeps \rightarrow Z$ in $\cH^2$ was determined in 
Corollary \ref{cor:cut:error}, this implies that $Z_t = \cZ_t \ m \times\P-a.e.$, which completes the proof
under the assumptions \HFd \ and \Hge.

\vspace{0.3cm}
We first need some intermediate upper bounds.
Analogously to Corollary \ref{cor:mom bd}, 
we have that
%the following bounds on the second moments of $\Yeps$ and $\Zeps$ hold:
%\begin{align*}
%\sup_{0\le r \le T}\normexp{\Yeps_r}^2   \le C, \quad \text{and} \quad
%\normexp{\Zeps_r}^2  \le \frac{C}{(T-r)^{1- \alpha\wedge(2\thetaC)}}\quad \text{for all } r\in[0,T).
%\end{align*}
%We will use the notation $\gamma:=(\frac\alpha2\wedge\thetaC +\frac\thetaL2)\wedge \thetaC$ hereafter.
%It follows analogously to \eqref{eq:bd:f} that 
\begin{equation}
 \normexp{f(r,X_r,\Yeps_r,\Zeps_r)} \le \frac{C}{(T-r)^{1-\gamma}} \quad \text{for all } r\in[0,T).
\label{eq:f:bound:repr} 
\end{equation}
Fix $t\in[0,T)$ and $\eta >0$.
Using the representation formula \eqref{eq:repr:eps}, %\[\Zeps_t = \E_t[\Phi(X_T) H^t_T + \int_t^T \feps(r,X_r,\Yeps_r,\Zeps_r)H^t_rdr],\]
it follows from Minkowski's inequality, the conditional Cauchy-Schwarz inequality  {(Corollary \ref{cor:cs:ineq})}, and Lemma \ref{lem:mal:weight} that
\begin{align}
 \normexp{\Zeps_t - \cZ_t} & = \normexp{ \E_t [ \int_t^T \big(\feps(r,X_r,\Yeps_r,\Zeps_r) - f(r,X_r, Y_r,Z_r) \big) H^t_r dr ] } \nonumber\\
& \le  \normexp{ \E_t[\int_t^T \big(\feps(r,X_r,\Yeps_r,\Zeps_r) - f(r,X_r, \Yeps_r,\Zeps_r) \big) H^t_r dr ] } \nonumber\\
& \quad +  C_M^{1/2} \int_t^T \frac{\normexp{ f(r,X_r,\Yeps_r,\Zeps_r) - f(r,X_r, Y_r,Z_r)} }{\sqrt{r-t} } dr .
\label{eq:1:repre:est:1}
\end{align}
Taking $\varepsilon < (T-t)/2$ and using \eqref{eq:f:bound:repr}, it follows that
\begin{align*}
\normexp{ \E_t[& \int_t^T \big(\feps(r,X_r,\Yeps_r,\Zeps_r) - f(r,X_r \Yeps_r,\Zeps_r) \big) H^t_r dr ] } 
%\\&
\le C_M^{1/2}  \int_{T-\varepsilon}^T \frac{\normexp{f(r,X_r,\Yeps_r,\Zeps_r)}}{\sqrt{r-t} }dr  \\
&
\le \frac{C_M^{1/2}}{\sqrt{T-t -\varepsilon}}   \int_{T-\varepsilon}^T \normexp{f(r,X_r,\Yeps_r,\Zeps_r)} dr  
% \\&
 \le \frac{C}{\sqrt{T-t} } \int_{T-\varepsilon}^T  \frac{dr}{(T-r)^{1-\gamma} } 
 = \frac{C\varepsilon^{\gamma}}{\sqrt{T-t}}.
\end{align*}
Taking $\varepsilon < \eta^{1/\gamma}(T-t)^{1/(2\gamma)}/C$, where $C$ is the last constant in the inequality above, is sufficient to bound the above term by $\eta$.
On the other hand, letting $\delta < (T-t)/2$,
\begin{align}
&\int_t^T \frac{\normexp{ f(r,X_r,\Yeps_r,\Zeps_r) - f(r,X_r, Y_r,Z_r)} }{\sqrt{r-t} } dr\nonumber \\
& \le  C_M^{1/2}  \frac{\int_{t+\delta}^T \normexp{ f(r,X_r,\Yeps_r,\Zeps_r) - f(r,X_r, Y_r,Z_r)} dr }{\sqrt\delta}  \nonumber \\
& \quad + C_M^{1/2} \int_{t}^{t+\delta} \frac{\normexp{f(r,X_r,\Yeps_r,\Zeps_r) - f(r,X_r, Y_r,Z_r)} }{\sqrt{r-t} } dr
\label{eq:rep:conv}
\end{align}
%To bound the first integral term in \eqref{eq:rep:conv}, we observe that 
%\[|f(r,X_r,\Yeps_r,\Zeps_r) - f(r,X_r, Y_r,Z_r)| \le L_f\{|Y_r - \Yeps_r| + |Z_r - \Zeps_r|\}(T-r)^{(\thetaL-1)/2}.\]
To bound the first integral term on the right hand side above, we  apply H\"older's inequality and the Lipschitz continuity of $f(t,\cdot)$ to obtain
\begin{align*}
C_M^{1/2} \int_{t+\delta}^T & \normexp{ f(r,X_r,\Yeps_r,\Zeps_r) - f(r,X_r ,Y_r,Z_r)} dr  \\
& \le C_M^{1/2}L_f \Big( \int_0^T \frac{dr}{(T-r)^{1-\thetaL} }  \Big)^{1/2} 
\Big( \sup_{0\le s\le T}\normexp{Y_s-\Yeps_s}^2 + \int_0^T \normexp{Z_r - \Zeps_r}^2 dr \Big) ^{1/2}
\end{align*}
Using that $(\Yeps,\Zeps) \rightarrow (Y,Z)$ in $\cS \times \cH^2$ as $\varepsilon \rightarrow 0$ (Corollary \ref{cor:cut:error}), 
set $\varepsilon$ sufficiently small so that the above is bounded above by $\sqrt{\delta}\eta$.
To bound the second integral term on the right hand side of \eqref{eq:rep:conv},  
we use \eqref{eq:bd:f} and \eqref{eq:f:bound:repr} combined with the triangle inequality to show that  
\begin{align*}
C_M^{1/2}\int_{t}^{t+\delta} & \frac{\normexp{f(r,X_r,\Yeps_r,\Zeps_r) - f(r,X_r, Y_r,Z_r)}}{\sqrt{r-t} } dr \\
& \le \frac{C}{(T-t-\delta)^{1-\gamma} } \int_t^{t+\delta} \frac{dr}{\sqrt{r-t} } 
 \le \frac {C\sqrt{\delta} }{(T-t)^{1-\gamma}} 
\end{align*}
and set $\delta$ sufficiently small so that the above is bounded above by $\eta$.
Therefore, we have shown that for almost every $t\in[0,T)$ and every $\eta>0$, there is a sufficiently small $\varepsilon$ such that
$\normexp{\Zeps_t - \cZ_t} < 3 \eta$.
In other words, $\E[|\Zeps_t - \cZ_t|^2] \rightarrow 0$ as $\varepsilon \rightarrow 0$ for  every $t$, as required.

% 
% \vspace{0.3cm}
% It remains only to show that $Z$ has a continuous version. We recall that $(\Yeps,\Zeps)$ solve a uniformly Lipschitz continuous BSDE 
% with terminal condition in $\L_{2,\alpha}$. It follows from \cite[Theorem 2.1]{gobe:makh:10} that, for every $\varepsilon > 0$,
% $\Zeps$ has a continuous version in $\cS^2$.
% Since

\vspace{0.3cm}
To prove the result without \HFd \ and \Hge, recall the mollified BSDE $(Y_M,Z_M)$ from Corollary \ref{cor:1:moll}.
Since $\Phi_M$ satisfies \Hge \ and $f_M$ satisfies \HFd, 
there is a predictable version $\cZ_M$ of $Z_M$ satisfying $\cZ_{M,t} = \E_t[\Phi_M(X_T) H^t_T + \int_t^T f_M(r,X_r,Y_{M,r},Z_{M,r}) H^t_r dr] $
for all $t\in[0,T)$ almost surely.
 {
Thanks to the point-wise convergence of $f_M$ to $f$ and $\Phi_M$ to $\Phi$, and the convergence of $(Y_M,Z_M)$ to $(Y,Z)$ in $\cS^2\times\cH^2$
from Corollary \ref{cor:1:moll}, we can use analogous limit arguments as above to complete the proof.
}
\qed

\section{Convergence rate of the Euler scheme for BSDEs}
\label{section:reg}
%For a given time-grid $\pi:=\{0=t_{0}<\ldots<t_{N} = T\}$, recall the definition of the $\L_2$-regularity given by
%Let $\cE_1(N)$ be the error due to discrete-time approximation, as defined in \eqref{eq:disc:error:gen}, of the BSDE with the Euler scheme \eqref{eq:1:dp} on the class of time-grids $\{\pib_N \ : \ N\ge 1\}$.
Throughout this section, the assumption \HFt \ is in force.
 {
Let us recall now the Euler scheme for BSDEs:
\begin{align*}
\Y NN & := \Phi(X_T), \quad \Z Ni := \frac{1}{t_{i+1} - t_i} \E_i[\Y N{i+1}(W_{t_{i+1}}-W_{t_i})^\top], \nonumber\\
\Y Ni & := \E_i[\Y N{i+1}+ f(t_i,X_{t_i}, \Y N{i+1}, \Z Ni)(t_{i+1} - t_i)].
%\label{eq:1:dp}
\end{align*}
}
We determine error estimates on  {the error of the Euler scheme, which is given by }
 {\[\cE(N) : = \max_{0\le i <N} \E[|Y_{t_i} - \Y Ni|^2] + \sum_{i=0}^{N-1} \int_{t_i}^{t_{i+1}} \E[|Z_{t} - \Z Ni|^2] dt.\] }
 {The following proposition serves as the starting point of our analysis; it allows us to estimate the error $\cE(N)$ using estimates for the so called $\L_2$-regularity, which we will do subsequently.}
\begin{proposition} 
\label{thm:gobe:lemo:ext}
Let $\beta \le \thetaL$. For the Euler scheme for BSDEs defined on the time-grids $\{\pib_N \ : \ N\ge 1\}$, there is a constant $C$ depending only on $L_f$, $ {L_X}$, $\thetaL$,  {$\theta_X$}, $\beta$, and $T$, but not on $N$, such that, for all $N\ge 1$,
\begin{align*}
\cE(N) \le CN^{-1} + C \sum_{i=0}^{N-1} \int_{t_{i}}^{t_{i+1}} \normexp{Z_{t} - \tilde Z_{t_{i}}}^2dt
\end{align*}
where
$
 \tilde Z_{t_i} := \frac{1}{\Delta_i} \E_i \big[\int_{t_i}^{t_{i+1}} Z_t dt \big].
$
\end{proposition}
The proof is analogous to the proof of \cite[Theorem 1]{gobe:lemo:06}, one must only use the result $\Delta_k /(T-t_k)^{1-\thetaL} \le T^{\thetaL}(\beta N)^{-1}$  for $\beta \le \thetaL$ (see Lemma \ref{lem:inc growth}) in order to compensate for the local Lipschtz constant of the driver.

 {The sum $ \sum_{i=0}^{N-1} \int_{t_{i}}^{t_{i+1}} \normexp{Z_{t} - \tilde Z_{t_{i}}}^2dt$ is called the $\L_2$-regularity; it's study was initiated by \cite{zhan:04}.}
 Since $(\tilde Z_{t_i}:= \frac{1}{\Delta_i} \E_i \big[\int_{t_i}^{t_{i+1}} Z_t dt \big])_i$ is the projection of $Z$ onto the space of adapted discrete processes with nodes on $\pi$ 
under the scalar product $(u,v) = \E \int_0^T ( u_s \cdot v_s ) ds$,
it follows that
\begin{equation}
\sum_{i=0}^{N-1} \int_{t_{i}}^{t_{i+1}} \normexp{Z_{t} - \tilde Z_{t_{i}}}^2dt \le \sum_{i=0}^{N-1} \int_{t_{i}}^{t_{i+1}} \normexp{Z_{t} -  Z_{t_{i}}}^2dt .
\label{eq:l2reg:deco}
\end{equation}
To bound $\cE(N)$, it follows from Proposition \ref{thm:gobe:lemo:ext} that it is sufficient to bound the term on the right-hand side of \eqref{eq:l2reg:deco}.
  { However, as in the proof of the Representation Theorem in Section \ref{section:repr}, it is not possible to do so directly for the BSDE $(Y,Z)$, so we use an approximation procedure via
%We also define
%\begin{equation}
%\gamma := (\thetaC\wedge\frac\alpha2 + \frac\thetaL2)\wedge\thetaC .
%\label{eq:gamma}
%\end{equation}
%
%
%
 the BSDE $(\Yeps,\Zeps)$, which we recall} from Definition \ref{def:1:intermediate BSDEs} in Section \ref{section:mal:BSDE}.

Throughout the remainder of this section, we work with the version of $Z$ and  $\Zeps$ given by Theorem \ref{thm:repr},  {i.e
\begin{align*}
Z_t = \e_t[\Phi(X_T) H^{t}_T + \int_{t}^T f(s,X_s,Y_s,Z_s) H^{t}_s ds] \quad 
\text{for all } t\in[0,T) \quad \P-a.s. , \\
\Zeps_t = \E_t[\Phi(X_T)H^t_T + \int_t^T \feps(r,X_r, \Yeps_r,\Zeps_r) H^t_r dr] \quad \text{for all } t\in[0,T) \quad \P-a.s.
\end{align*}
This version empowers us with the additional a priori estimates estimates developed in Section \ref{section:regularity}; we use these estimates  frequently in the analysis of this section.
}

%Since the representation Theorem \ref{thm:repr} applies to $\Zeps$, we work with the version of $\Zeps$ given by
%\[\Zeps_t = \E_t[\Phi(X_T)H^t_T + \int_t^T \feps(s,X_s,\Yeps_s,\Zeps_s) H^t_s ds ] \quad \text{for all } t\in[0,T) \ \text{almost surely} .\]
The following lemma decomposes the $\L_2$-regularity of $Z$  {-- the left hand side of equation \eqref{eq:l2reg:deco} --} into the $\L_2$-regularity of $\Zeps$ and  {a small correction} term controlled by $\varepsilon$.
\begin{lemma}
\label{lem:Z -Zeps:L2}
Let $\beta \in (0,1]$.  
Then there is a constant $C$ depending only on $L_f$, $C_M$, $\thetaL$, $\thetaC$, $\beta$, $C_f$, $K^\alpha(\Phi)$, and $T$, such that for all $ N\ge1  $
\begin{align*}
%\cE(\pib_N) \le
 \sum_{i=0}^{N-1} \int_{t_i}^{t_{i+1}} \normexp{Z_s - \tilde Z_{t_i} }^2 ds \le 
C\sum_{i=0}^{N-1} \int_{t_i}^{t_{i+1}} \normexp{\Zeps_s - \bZeps_{t_i} }^2 ds + C \varepsilon^{2\gamma}
%CN^{-2\gamma/\beta} + C\varepsilon^{2\gamma}\big(1 + \ln(N)\big).
\end{align*}
 {where $ \tilde Z_{t_i}:= \frac{1}{\Delta_i} \E_i \big[\int_{t_i}^{t_{i+1}} Z_t dt \big]$, $ \bZeps_{t_i}:= \frac{1}{\Delta_i} \E_i \big[\int_{t_i}^{t_{i+1}} \Zeps_t dt \big]$, and 
$\gamma := (\thetaC\wedge\frac\alpha2 + \frac\thetaL2)\wedge\thetaC $.}
\end{lemma}
{\bf Proof.}
In what follows,
$C$  may change in value from
line to line.
Using the Cauchy inequality and the orthogonality of the projections, $\frac12 \sum_{i=0}^{N-1} \int_{t_i}^{t_{i+1}} \normexp{Z_s - \tilde Z_{t_i} }^2 ds \le   \int_{0}^{T} \normexp{Z_s -  \Zeps_{s} }^2 ds 
+ \sum_{i=0}^{N-1} \normexp{ \tilde Z_{t_i} - \bZeps_{t_i} }^2 \Delta_i + \sum_{i=0}^{N-1} \int_{t_i}^{t_{i+1}} \normexp{\Zeps_s -  \bZeps_{t_i} }^2 ds $.
%
% Recalling that $\Zeps_s = z_s + \zeps_s$, and combining \eqref{eq:zeps:bd:2} with 
% $\sum_{i=0}^{N-1} \int_{t_i}^{t_{i+1}} \normexp{z_s - z_{t_i} }^2 ds \le CN^{-1}$, shown in
% \cite[Theorem 1.3]{gobe:makh:10}, it follows that 
% \begin{equation}
% \sum_{i=0}^{N-1} \int_{t_i}^{t_{i+1}} \normexp{\Zeps_s - \Zeps_{t_i} }^2 ds \le CN^{-1}+ \frac{CN^{-2}}{\varepsilon^{1-\thetaL}}
% \label{eq:l2 reg:Zeps}
% \end{equation}
% 
Recall from Corollary \ref{cor:cut:error} that
$\int_0^T \normexp{Z_s - \Zeps_s}^2 ds 
 \le C\varepsilon^{2\gamma}$.
 Moreover, using Jensen's inequality,
 \begin{align*}
  \sum_{i=0}^{N-1} \normexp{ \tilde Z_{t_i} - \bZeps_{t_i} }^2 \Delta_i
  =  \sum_{i=0}^{N-1} \normexp{ \frac{1}{\Delta_i} \E_i\big[\int_{t_i}^{t_{i+1}} (Z_t - \Zeps_t) dt }^2 \Delta_i
  \le \int_0^T \normexp{Z_s - \Zeps_s}^2 ds \le C \varepsilon^{2\gamma}
 \end{align*}
and this completes the proof.
%Using the bound \eqref{eq:cut:error:2} with $t=t_i$ for every $i\in\{0,\ldots,N-1\}$ gives
%\begin{align*}
%\sum_{i=0}^{N-1} & \normexp{Z_{t_i} - \Zeps_{t_i}}^2 \Delta_i 
%\le C\sum_{i=0}^{N-1} \Big( \int_{t_i\vee(T-\varepsilon)}^T \frac{ds }{(T-s)^{1-\gamma}\sqrt{s-t_i}}\Big)^2\Delta_i \\
%& \le C \Big( \int_{t_{N-1}\vee(T-\varepsilon)}^T \frac{ds }{(T-s)^{1-\gamma}\sqrt{s-t_{N-1}}}\Big)^2\Delta_{N-1} 
%%\\& \qquad 
%+C\sum_{i=0}^{N-2} \Big( \int_{t_i\vee(T-\varepsilon)}^T \frac{ds }{(T-s)^{1-\gamma}}\Big)^2 \frac{\Delta_i}{t_{N-1}-t_i}\\
%& \le C(\Delta_{N-1}^{\gamma - 1/2})^2 \Delta_{N-1} + C\varepsilon^{2\gamma} + C\varepsilon^{2\gamma}\int_0^{t_{N-2}} \frac{ds}{t_{N-1}-s} 
%%\\& 
%\le CN^{-2\gamma/\beta} + C\varepsilon^{2\gamma}\big(1+\ln(N)\big).
%\end{align*}
%The first term in the last inequality comes from substituting $\Delta_{N-1} = TN^{-1/\beta}$.
% Combining the above result with \eqref{eq:l2 reg:Zeps} implies that
% \begin{align*}
% \sum_{i=0}^{N-1} \int_{t_i}^{t_{i+1}} \normexp{Z_s - Z_{t_i}}^2 ds & \le C \int_0^T \normexp{Z_s - \Zeps_s}^2 ds
% + C\sum_{i=0}^{N-1} \normexp{Z_{t_i} - \Zeps_{t_i}}^2 \Delta_i \\
% &\qquad + C\sum_{i=0}^{N-1} \int_{t_i}^{t_{i+1}} \normexp{\Zeps_s - \Zeps_{t_i} }^2 ds \\
% & \le C\varepsilon^{2\gamma}\ln(N) + CN^{-1} + \frac{CN^{-2}}{\varepsilon^{1-\thetaL}}
% \end{align*}
\qed

\vspace{0.3cm}
%In the following proposition, we obtain a  convergence rate for the $\L_2$-regularity $\cE(\pi)$ when $\pi$ is a time grid of the form $\pib$.
%We do not use the results of Section \ref{section:regularity}.
%We will obtain a more precise convergence rate under stronger assumptions later in Theorem \ref{thm:l2 reg:z}, but this first result will serve, for pedagogical purposes, to show that
%the results of Section \ref{section:regularity} are useful for the estimation of the $\L_2$-regularity.
We now come to our first and most general estimate on the $\cE(N)$.
Later, in Theorem \ref{thm:L2 reg},  we augment this result with stronger assumptions.
\begin{theorem}
\label{prop:1:l2 reg:makh}
Let $0<\beta < (2\gamma) \wedge \alpha$
 {and $\gamma := (\thetaC\wedge\frac\alpha2 + \frac\thetaL2)\wedge\thetaC $.}
There is a constant $C$ depending only on $L_f$, $C_M$, $\thetaL$, $\thetaC$, $\beta$, $C_f$, $K^\alpha(\Phi)$, and $T$, but not on $N$, such that for all $N\ge1$,
\[ \cE(N) \le CN^{-1} \1_{[1,2]}(\alpha + \thetaL) + CN^{-2\gamma}\1_{(0,1)}(\alpha + \thetaL)\]
%C N^{0.7\big(1-(\alpha+\thetaL)\wedge 1\big)/\gamma-1}.\]
\end{theorem}
{\bf Proof.}
In what follows,
$C$ may change in value from
line to line.
From Proposition \ref{thm:gobe:lemo:ext}, it is sufficient to bound $\sum_{i=0}^{N-1} \int_{t_{i}}^{t_{i+1}} \normexp{Z_{t} - \tilde Z_{t_{i}}}^2dt $.
To start with, assume \HFd \ and \Hge.
Recall the BSDEs $(\yeps,\zeps)$ from Definition \ref{def:1:intermediate BSDEs}
% , $(\Yeps,\Zeps)$
 and $(\Ueps,\Veps)$ from \eqref{eq:1:def:ueps veps} in Section \ref{section:mal:BSDE}. 
%It was shown in  Lemma \ref{lem:mal:ueps veps}  that $\zeps_t = \Ueps_t \sigma(t,X_t)$ and 
%$(\Veps_{j,t})^\top  =  (\nabla\zeps_{j,t})^\top \sigma^{-1}(t,X_t) - \Ueps_t \nabla_x \sigma_j(t,X_t) \ m\times\P-a.e.$  
%for all $j\in\{1,\dots,q\}$, where $\Veps_{j}$ (resp. $ \nabla\zeps_{j}$) is the $j$-th column of $\Veps$ (resp. $\nabla \zeps$),
%and $\sigma_j$ is the $j$-th column of $\sigma$.
In the proof of \cite[Theorem 3.1]{gobe:makh:10}, the authors show that for any $i$ and $s\in[t_i,t_{i+1})$,
\begin{align}
 \normexp{\zeps_s - \zeps_{t_i} }
& \le C  \int_{t_i}^{s} \normexp{\aeps r} dr
+ C \int_{t_i}^{s} \normexp{\Veps_r} dr  + C \Delta_i^{1/2} 
\label{eq:zeps:bd:1}.
\end{align}
%Therefore, to determine an upper bound for $\sum_i\int_{t_i}^{t_{i+1}}\normexp{\zeps_s - \zeps_{t_i} }^2 ds$,
% the $\L_2$-regularity of $\zeps$ in terms of $N$, 
%it suffices to determine an upper bound  for
%\[\sum_{i=0}^{N-1} \int_{t_{i}}^{t_{i+1}}\{(\int_{t_i}^{t} \normexp{\aeps{r}} dr )^2 
%+ (\int_{t_i}^{t} \normexp{\Veps_{r}} dr )^2\} dt.\] 
%
Using $(\int_0^T\normexp{\aeps r} dr)^2 +  \int_{0}^{T} \normexp{\Veps_{r}}^2 dr \le C { \varepsilon^{-1+(\thetaL +  \alpha)\wedge1} }$ 
from \eqref{eq:mal:zeps:lin} in Lemma \ref{lem:mal:ueps veps}, and \eqref{eq:l2reg:deco}, it follows from Jensen's inequality that
\begin{equation*}
\sum_{i=0}^{N-1} \int_{t_i}^{t_{i+1}} \normexp{\zeps_s - \bzeps_{t_i} }^2 ds \le \frac CN + 
\frac{C\max_{0\le i \le N-1}\Delta_i}{\varepsilon^{1 - (\thetaL + \alpha)\wedge1} } \le  CN^{-1}(1  + 
\varepsilon^{(\thetaL + \alpha)\wedge1 -1}  )
\end{equation*}
where $\max_i\Delta_i \le CN^{-1}$ follows from \eqref{eq:grid bd:1} in Lemma \ref{lem:inc growth}.
Combining this estimate with
$
% \[
 \sum_{i=0}^{N-1} \int_{t_i}^{t_{i+1}} \normexp{z_s- z_{t_i} }^2 ds \le CN^{-1}
$,% \] 
 \ shown in
\cite[Theorem 1.3]{gobe:makh:10},  $\Zeps = z + \zeps$,
and the results of Lemma \ref{lem:Z -Zeps:L2}, \eqref{eq:l2reg:deco}
it follows that
\begin{equation}
 \sum_{i=0}^{N-1} \int_{t_i}^{t_{i+1}} \normexp{Z_s - \tilde Z_{t_i} }^2 ds \le CN^{-2\gamma/\beta} + C\varepsilon^{2\gamma}%\big(1 + \ln(N)\big) 
+ CN^{-1}(1  + \varepsilon^{(\thetaL + \alpha)\wedge1 -1}  ).
%+ \frac{C}{N\varepsilon^{1 - \big((\thetaL + \alpha)\wedge1\big)}} + \frac CN.
\label{eq:cE1:eps}
\end{equation}
 {To complete the proof  under \HFd \ and \Hge,  let $\varepsilon := N^{-\delta}$ in the estimate \eqref{eq:cE1:eps},
 take $\delta := 1/(2\gamma)$ if $\alpha + \thetaL \ge 1$ and $\delta := 1$ otherwise, and notice that $2\gamma \le \alpha + \thetaL$.}

In order to prove the general result, recall the BSDE $(Y_M,Z_M)$ from Corollary \ref{cor:1:moll};
%and 
%it satisfies the $\L_2$ regularity result of the proposition statement because 
its terminal condition satisfies \Hge \ and its driver satisfies \HFd.
Moreover, \cite[Lemma 3.1]{gobe:makh:10} proves $K^\alpha(\Phi_M) \le K^\alpha(\Phi)$.
Therefore, working with the version of $Z_M$ given by the representation formula  {$Z_{M,t} =  \e_t[\Phi_M(X_T) H^{t}_T + \int_{t}^T f_M(s,X_s,Y_{M,s},Z_{M,s}) H^{t}_s ds] $} from Theorem \ref{thm:repr},
%\[ 
%Z_{M,t} = \E_t[\Phi_M(X_T) H^t_T + \int_t^T f_M(r,X_r,Y_{M,r},Z_{M,r}) H^t_r dr] 
%\quad \text{for all } t\in[0,T) \ \text{almost surely,}
%\]
%the result is extended to $(Y,Z)$ using
it follows from the triangle inequality and the results obtained above that
\[
\cE(N)  \le 2 \int_{0}^{T} \normexp{Z_s - Z_{M,s} }^2 ds
+ 2 \sum_{i=0}^{N-1} \int_{t_i}^{t_{i+1}} \normexp{Z_{M,s} - \tilde Z_{M,t_i}}^2 ds 
% + CN^{-1} + C N^{\big(1-(\alpha+\thetaL)\wedge 1\big)/\gamma -1},
\]
and letting $M\rightarrow \infty$ with Corollary \ref{cor:1:moll} yields the result.
\qed

\section{ {A priori estimates under \Hge \ and \Hgh}}

\vspace{0.3cm}
At the end  this section, we give a complementary result to Theorem \ref{prop:1:l2 reg:makh} under  {stronger the conditions on the terminal condition} 
%\Hgexp \ and
\Hge \ and \Hgh, i.e. where the function $\Phi$ is bounded (and/)or H\"older continuous, respectively. 
This is achieved using the an additional a priori estimates on $\normexp{\Veps_t}$, given in Proposition \ref{prop:veps:bd} below. 
Moreover, these a priori estimates will be critical in Section \ref{section:num}, where one requires more structure than in Section \ref{section:reg}.
 {The result is proved, roughly speaking, by using a functional representation of the intermediate process $\zeps$ and show Lipschitz continuity of the said functional representation.
This adds an additional layer of interest  under \Hgh \ for the parameters  $\thetaP + \thetaL \ge 1$, where we can demonstrate that limit of the process $\zeps_s$ in $\cH^2$, i.e. the process $Z_s - \nabla_x u(s,X_s)\sigma(s,X_s)$, has a functional representation and that function is Lipschitz continuous; see Corollary \ref{cor:proxy:lip}. 
Regularity results are important for numerical schemes as they allow one to build algorithms with lower numerical complexity -- see for example \cite[Section 3.5]{gobe:turk:13b} --
and this regularity result has such implications for the proxy scheme described in the introduction of this paper.
}

 {
First, we state the result that $x \mapsto \sigma^{-1}(t,x)$  is uniformly Lipschitz continuous, and $t \mapsto  \sigma^{-1}(t,x)$ is uniformly $1/2$-H\"older continuous.
This elementary result will also be useful in Section \ref{section:num} below.
The proof is to be found in Appendix \ref{section:mx}.
\begin{lemma}
\label{lem:1:inv sigma:lip}
The right inverse matrix $\sigma(t,\cdot)^{-1}$ is Lipschitz continuous uniformly in $t$
and $\sigma^{-1}(\cdot,x)$ is $1/2$-H\"older continuous uniformly in $x$. Its Lipschitz (resp. H\"older) constant depends $\|\sigma\|_\infty$,
$\|\nabla_x \sigma \|_\infty$
and $\elip$ only, but not on $(t,x)$.  Moreover, $\| \sigma^{-1}\|_\infty \le \|\sigma\|_\infty/ \elip$.
\end{lemma}
}

We now state the main result of this section, the a priori estimates on the process $\Veps$.

\begin{proposition}
\label{prop:veps:bd}
Suppose that \HFd \ is in force and $\Phi(x) $ is not zero everywhere in $\R^d$. 
 If \Hge \ is in force, there exists version of $\Veps$ and a finite constant $C$ depending only on $L_f$, the bounds on $b$ and $\sigma$ and their partial derivatives, $\elip$, $C_M$, $\thetaL$, $\thetaC$, $C_f$, and $T$
such that for any $\varepsilon \in(0,T]$ and every $t\in[0,T)$, 
$\normexp{\Veps_t} \le 
%C (T-t)^{(\alpha-2)/2} +
 C \phi(t,\varepsilon, \thetaL
% ,\thetaC,\alpha
) $, where
\begin{equation}
 \phi(t,\varepsilon, \thetaL
% ,\thetaC,\alpha
) := 
% 1 + \frac{1}{(T-t)^{(1-2\gamma)/2}} + 
\| \Phi\|_\infty \int_t^{T-\varepsilon}\frac{dr}{(T-r)^{(3-\thetaL
% -\alpha
)/2} \sqrt{r-t} }.
\label{eq:phi:2}
\end{equation}
If \Hgh \ is in force, there exists a version of $\Veps$,
% and a finite constant $C$ depending only on $L_f$, the bounds on $b$ and $\sigma$ and their partial derivatives, $\elip$, $C_M$, $\thetaL$, $\thetaC$, $C_f$, $  {K_\Phi}$ and $T$
such that for any $\varepsilon \in(0,T]$ and every $t\in[0,T)$, 
$\normexp{\Veps_t} \le 
%C (T-t)^{(\alpha-2)/2} +
 C  \phi(t,\varepsilon, \thetaL, \thetaP
% ,\thetaC,\alpha
) $, where
\begin{equation}
 \phi(t,\varepsilon, \thetaL, \thetaP
% ,\thetaC,\alpha
) := 
% 1 + \frac{1}{(T-t)^{(1-2\gamma)/2}} + 
 K_\Phi \int_t^{T-\varepsilon}\frac{dr}{(T-r)^{(3-\thetaL - \thetaP
% -\alpha
)/2} \sqrt{r-t} }.
\label{eq:phi:3}
\end{equation}
\end{proposition}
\begin{remark*}
The integrals in (\ref{eq:phi:2},\ref{eq:phi:3}) exist and are bounded by $C\varepsilon^{-(1-\thetaL)/2}(T-t)^{(\alpha-1)/2}$.
\end{remark*}
\noindent{\bf Proof.}
In what follows,
$C$
% depends only on $L_f$, the bounds on $b$ and $\sigma$ and their partial derivatives, $\elip$, $C_M$, $\thetaL$, $\thetaC$, $C_f$, and $T$, 
%but not on $\varepsilon$, whose value
  may change from line to line.

\noindent
{\bf Step 1. Functional and BSDE setup.}
For all $(t,x)\in [0,T)\times\R^d$, consider the FBSDE
\begin{align}
%\left.
%\begin{array}{rcl}
\yM tx_s  =  \int_s^T F(r,\XM tx_r,  \yM tx_r,  \zM tx_r)dr - \int_s^T \zM tx_r dW_r,
\qquad s \in [t,T) ,
%\end{array}
%\right\}
\label{eq:markov:dyn}
\end{align}
where $F(t,x,y,z) = \feps(t,x,u(t,x)+y,(\nabla_x u(t,x) \sigma(t,x) )^\top+ z)$ and $\XM tx$ is the solution of the SDE \eqref{eq:flow}.
Note that the BSDE $(\yeps,\zeps)$ from Section \ref{section:mal:BSDE} is equal to $(\yM 0{x_0},\zM 0{x_0})$
because, thanks to Lemma \ref{lem:lin:pde:bds}, $(y,z)$ is equal to  $\big(u(\cdot,X_\cdot),\nabla_xu(\cdot,X_\cdot)\sigma(\cdot,X_\cdot)\big)$ and $X$ is equal to $\XM 0{x_0}$.
Since $\feps(t,\cdot)$ is Lipschitz continuous for all $t\in[0,T]$, $F(t,\cdot)$ is also Lipschitz continuous, with Lipschitz constant $C \1_{[0,T-\varepsilon)}(t) \varepsilon^{(\thetaL - 3)/2} $, for all $t\in[0,T)$;  {see the first paragraph of the proof of Theorem \ref{thm:repr} for detailed computations}.
%:
%for all $(x,y,z),(x',y',z')\in \R^d \times \R \times \R^q$
%\begin{align*}
%|F(t,x,y,z) - F(t,x',y',z')|  
%& \le \frac{C \1_{[0,T-\varepsilon)}(t) }{\varepsilon^{(1-\thetaL)/2} }(|x-x'| + |y-y'| + |z-z'|) \\
%& \qquad + \frac{C \1_{[0,T-\varepsilon)}(t) }{\varepsilon^{(1-\thetaL)/2} }(|u(t,x)-u(t,x')| + |\nabla_xu(t,x) - \nabla_xu(t,x')|) \\
%& \le \frac{C \1_{[0,T-\varepsilon)}(t) }{\varepsilon^{(3  -\thetaL)/2} }(|x-x'| + |y-y'| + |z-z'|)
%\end{align*}
%where the last inequality follows from the fact that $u(t,\cdot)$ and $\nabla_xu(t,\cdot)$ are differentiable,
%and their derivatives are bounded by $C(T-t)^{-1/2}$ and $C(T-t)^{-1}$ respectively, Lemma \ref{lem:lin:pde:bds}.
Now, letting 
\[
\MW t{x,s}_r := \frac{\1_{(t,T]}(s)}{r-s} (\int_s^r \sigma^{-1}(r,\XM tx_r) D_s \XM tx_r )^\top dW_r)^\top\] 
where $D_s \XM tx$ is the 
Malliavin derivative of $\XM tx_s$, 
%as defined in  Section \ref{section:mal:SDE}. 
%One can show analogously to Theorem \ref{thm:repr},
%with $\XM tx$ replacing $X$ and $\MW t{x,s}_r$ replacing $H^s_r$, 
it follows from \cite[Theorem 4.2]{ma:zhan:02}, because the terminal condition of the BSDE satisfied by $(\yM tx, \zM tx)$ is zero,
that $\zM tx_r$ is equal to $\zeps(r,\XM tx_r)$ $m \times \P$-almost everywhere, 
where  $\zeps:[0,T)\times\R^d \rightarrow (\R^q)^\top$ is a continuous, deterministic function given by
\begin{equation}
\zeps(t,x) := \E[\int_t^T F(r,\XM tx_r,\yM tx_r, \zM tx_r) \MW t{x,t}_r dr];
 \label{eq:zeps:fn}
\end{equation}
%that there is a predictable (in fact continuous) version $\cZ$ of $\zM tx$ such that
%\begin{equation}
%\cZ_s = \E_s[\int_s^T F(r,\XM tx_r,\YM tx_r, \ZM tx_r) \MW t{x,s}_r dr] \quad \text{for all } s\in[t,T) \ \text{almost surely;}
%\label{eq:zeps:markov}
%\end{equation}
we work with this version of $\zM tx$ from hereon.
 {
Additionally, we  show in Step 3 below  that the process 
\[( (\nabXM tx{s} )^\top \nabla_x\zeps(s, \X {t,x}s) ) _{0\le s \le T}\]
(the derivative here is in the weak sense) is a version of the process $(\nabla \zM tx_s)_{0\le s \le T}$, which is a part of the solution
$(\nabla \yM tx, \nabla \zM tx)$ of the BSDE
\begin{align}
\nabla \yM tx_\tau & = \int_\tau^T \nabla_x\feps(\Theta_r)\nabXM tx{r} + \nabla_y\feps(\Theta_r) \{ \nabla_x u(r,\X {t,x}r) \nabXM tx{r}
+ \nabla \yM tx_r\} dr  \nonumber \\ 
& \quad + \int_\tau^T \nabla_z\feps(\Theta_r)U(r,\X {t,x}r)^\top \nabXM tx{r} + \sum_{j=1}^q \nabla_z\feps_{j}(\Theta_r) (\nabla \zM tx_{j,r})^\top  dr 
- \sum_{j=1}^q \int_\tau^T  (\nabla \zM tx_{j,r})^\top dW_r,
\label{eq:grad:y:2}
\end{align}
}
 {
where $\Theta_r = (r,\X{t,x}r,\YM tx_r,\ZM tx_r)$; the function $U(t,x) $ is defined
\begin{equation*}
U(t,x) := \nabla_x^2 u(t,x) \sigma(t,x) + \sum_{j=1}^d (\nabla_x u)_j (t,x) \nabla_x\sigma^\top_j(t,x) 
%\label{def:deriv du s}
\end{equation*}
for the function $u$ defined in Lemma \ref{lem:lin:pde:bds}.
%which we prove in Step 3 below.
Note that the BSDE \eqref{eq:grad:y:2} is a generalization to the BSDE \eqref{eq:grad:y} -- solved by $(\nabla \yeps, \nabla \zeps)$ -- which we recall for convenience:
\begin{align*}
\nabla \yeps_t & = \int_t^T \nabla_x \feps(\Theta_r)\nabX{r} + \nabla_y\feps(\Theta_r) (\nabla_x u(r,X_r) \nabX{r}
+ \nabla \yeps_r) dr  \nonumber \\ 
& \quad + \int_t^T \nabla_z\feps(\Theta_r)U(r,X_r)^\top \nabX{r} + \sum_{j=1}^q \nabla_z\feps_{j}(\Theta_r) (\nabla \zeps_{j,r})^\top  dr 
- \sum_{j=1}^q \int_t^T  (\nabla \zeps_{j,r})^\top dW_r;
%\label{eq:grad:y}
\end{align*}
 indeed, in \eqref{eq:grad:y}, 
% process $X$ is started at   time $0$ with value $x_0$, whereas $\X{t,x}{}$ is started at time $t$ with value $x$.
% Now, we
  set $t \equiv 0$  and $x \equiv x_0$.
}

\vspace{0.3cm}
\noindent
{\bf Step 2. Proof assuming $\zeps(t,\cdot)$ is Lipschitz continuous with $|\nabla_x \zeps (t,\cdot) | \le C \phi(t, \cdot)$  and $(\nabXM tx{s} \big)^\top \nabla_x\zeps(s, \X {t,x}s)$ is a version of $\nabla \zM tx_s$.}
%To prove the proposition, it will be sufficient to show that $\zeps(t,\cdot)$ is Lipschitz continuous for all $t \in [0,T)$ with Lipschitz constant $C\phi(t,\cdot)$ for $\phi$ corresponding  to either \Hge \ or \Hgh. 
% Assuming that, for a mollified version of $\zeps(t,\cdot)$,    $|\nabla_x \zeps_r (t,\cdot) | \le C \phi(\tau, \cdot)$,
%the Cauchy-Schwarz inequality combined with Lemma \ref{lem:mal:SDE:3} and 
The hypothesis $|\nabla_x \zeps (t,\cdot) | \le C \phi(t, \cdot)$ implies that
\[
%\normexp{\big( \sigma^{-1}(s,X_s)\nabX{s} \big)^\top D_s\zeps(t,X_t)}
 { \
\normexp{\nabla \zeps _t} = 
\normexp{ ( \nabX t)^\top \nabla_x \zeps(t, X_t)}
\le
 \| \nabX t\|_2  \| \nabla_x \zeps(t, \cdot) \|_\infty
}
\le C \phi(t,\cdot) \quad \text{ for all } s.\]
%$\normexp{D_s \zeps(t,X_t)} \le C \phi(t, \varepsilon,\thetaL)$; 
%Since the bound on $\|\nabla_x \zeps_r (t,\cdot) \|_\infty$ is independent of the mollification, 
%%it follows that $\zeps(t,X_t) \in \mal$ for all $t$, and the above estimate holds.
%the result is valid when the derivative $\nabla_x$ is in the weak sense (i.e., $\zeps(t,\cdot)$ is only Lipschitz continuous).
%Using the version  of 
%%$\zeps_t$ given by $\zeps(t,X_t)$, and 
%$\nabla \zeps_t$ given by
%$ \big( \sigma^{-1}(s,X_s)\nabX{s} \big)^\top D_s\zeps (t,X_t)$,
%%from 
%%\cite[Lemma 2.4]{ma:zhan:02}
%it follows that $\normexp{\nabla \zeps_t} \le C \phi(t, \cdot)$.
Now, using Lemma \ref{lem:mal:ueps veps},
\[\|\sup_{s \le r < T} \Ueps_r \|_2 \le C\int_s^{T-\varepsilon} \normexp{\aeps{r}} dr 
\le C\int_s^{T-\varepsilon} \frac{ dr}{(T-r)^{(3-\alpha)/2} } \le C \phi(s,\cdot),\]
and
$(\Veps_{j,t})^\top =  (\nabla\zeps_{j,t})^\top \sigma^{-1}(t,X_t) - \Ueps_t \nabla_x \sigma_j(t,X_t)$, therefore
 we conclude that 
$\normexp{\Veps_t} \le  C \phi(s,\cdot)$ as required.

\vspace{0.3cm}
\noindent
 {\bf Step 3. Proving that $\big( \nabXM tx{s} \big)^\top \nabla_x\zeps(s, \X {t,x}s)$ is a version of $\nabla \zM tx_s$.}
We make use of Malliavin calculus -- see Section \ref{section:mal:cal}. By taking the Malliavin derivative on both the BSDE solution $(\yeps,\zeps)$ and on the functional representation $\zeps(s,\XM tx_s)$, we obtain an intermediate version that is equal for both.

 {
{\bf $\blacktriangleright$ BSDE arguments.}
There is a version (see \cite[Lemma 2.2]{gobe:makh:10} for the proof) of the processes $(D_s \yM tx_\tau, D_s \zM tx_\tau)_{s\le \tau\le T} $, the Malliavin derivatives of the processes $(\yeps, \zeps)$, solving the BSDE
\begin{align}
D_s \yM tx_\tau & = \int_\tau^T  \nabla_x \feps(\Theta_r)D_s \X{t,x}{r} + \nabla_y\feps(\Theta_r) ( \nabla_x u(r,X_r) D_s \X{t,x}{r}
+ D_s \yM tx_r ) dr  \nonumber \\ 
& \quad + \int_t^T  \nabla_z\feps(\Theta_r)U(r,\X{t,x}r)^\top D_s \X{t,x}{r} + \sum_{j=1}^q \nabla_z\feps_{j}(\Theta_r) (D_s \zM tx_{j,r})^\top   dr \nonumber \\ 
& \qquad - \sum_{j=1}^q \int_t^T  (D_s \zM tx_{j,r})^\top dW_r.
\label{eq:mall:y}
\end{align}
}
 {
We  multiply \eqref{eq:mall:y} on the right by $ \sigma^{-1}(s,\X{t,x}s)\nabXM tx{s}  $
and apply Lemma \ref{lem:mal:SDE:3} to obtain
\begin{align}
D_s \yM tx_\tau  \sigma^{-1}(s,\X {t,x}s)\nabXM tx{s}  & = \int_\tau^T \nabla_x\feps(\Theta_r)\nabXM tx{r} + \nabla_y\feps(\Theta_r) (\nabla_xu(r,\X {t,x}r) \nabXM tx{r}) dr \nonumber\\
& +  \int_\tau^T  \nabla_y\feps(\Theta_r) D_s \yM tx_r  \sigma^{-1}(s,\X{t,x}s)\nabXM tx{s} dr  \nonumber \\ 
& + \int_\tau^T \nabla_z \feps(\Theta_r)U(r,\X{t,x}r)^\top  \nabXM tx{r} dr \nonumber\\
& + \sum_{j=1}^q \int_\tau^T \nabla_z\feps_{j}(\Theta_r) (\big( \sigma^{-1}(s,\X{t,x}s)\nabXM tx{s} \big)^\top  D_s \zM tx_{j,r})^\top  dr \nonumber\\
& - \sum_{j=1}^q \int_\tau^T  (\big( \sigma^{-1}(s,\X{t,x}s)\nabXM tx{s} \big)^\top D_s \zM tx_{j,r})^\top dW_r;
\label{eq:grad:y:3}
\end{align}
comparing the  BSDE \eqref{eq:grad:y:3} to to \eqref{eq:grad:y:2}
}
 {
term by term, it is clear that 
\[(D_s \yM tx_{\tau}  \sigma^{-1}(s,\X {t,x}s)\nabXM tx{s}, \big( \sigma^{-1}(s,\X{t,x}s)\nabXM tx{s} \big)^\top D_s \zM tx_{\tau})_{s \le \tau \le T},\]
a version of the solution to \eqref{eq:grad:y:3}, is a version of $(\nabla \yM tx_\tau, \nabla \zM tx_\tau)_{s \le \tau \le T}$, the solution to \eqref{eq:grad:y:2},
for all $s \in [0,T]$.
}

 {
{\bf $\blacktriangleright$ Functional arguments.}
We start by assuming that $\zeps(t,\cdot)$ is smooth (or by taking a mollification).
The chain-rule of Malliavin calculus -- Lemma \ref{lem:mal:chain rule} -- yields $D_s \zeps(\tau,\X {t,x}\tau)$ equals $ ( D_s \X {t,x}\tau )^\top \nabla_x \zeps(\tau,\X {t,x}\tau)$,
and, applying Lemma \ref{lem:mal:SDE:3},  $\big( \sigma^{-1}(s,\X{t,x}s)\nabXM tx{s} \big)^\top  D_s\zeps(\tau,\X {t,x}\tau)$ is equal to $(\nabXM tx \tau)^\top \nabla_x \zeps(\tau, \X {t,x}\tau)$.
The result follows for $\zeps(\tau, \cdot)$ only Lipschitz continuous by standard limiting arguments.
Since $(\zeps(\tau,\X{t,x}\tau))_{0\le \tau \le T}$ 
is a version of $(\zM tx_\tau)_{0\le \tau \le T}$, it follows that $(D_s \zeps(\tau,\X{t,x}\tau)))_{s \le \tau \le T}$ is a version of $(D_s \zM tx_\tau)_{s \le \tau \le T}$, and therefore that 
\[
\big( \sigma^{-1}(s,\X{t,x}s)\nabXM tx{s} \big)^\top D_s \zM tx_{\tau})_{s \le \tau \le T} = ((\nabXM tx \tau)^\top \nabla_x \zeps(\tau, \X {t,x}\tau))_{s \le \tau \le T}
\]
is a version of $(\big( \sigma^{-1}(s,\X{t,x}s)\nabXM tx{s} \big)^\top D_s \zM tx_{\tau})_{s \le \tau \le T}$ for all $s \in [0,T]$.
}

\vspace{0.3cm}
 {
We now combine the BSDE arguments and the functional arguments from above.
Thanks to the intermediate version $(\big( \sigma^{-1}(s,\X{t,x}s)\nabXM tx{s} \big)^\top D_s \zM tx_{\tau})_{0 \le \tau \le T}$, it follows that
\[((\nabXM tx \tau)^\top \nabla_x \zeps(\tau, \X {t,x}\tau))_{0 \le \tau \le T}\] is a version of $( \nabla \zM tx_\tau)_{0 \le \tau \le T}$.
%$\big( \sigma^{-1}(s,\X{t,x}s)\nabXM tx{s} \big)^\top  D_s\zeps(\tau,\X {t,x}\tau) $ is a version of $\big( \sigma^{-1}(s,\X{t,x}s)\nabXM tx{s} \big)^\top  D_s\zM tx_\tau$.
%Finally, we have that $(\nabXM tx \tau)^\top \nabla_x \zeps(\tau, \X {t,x}\tau)$ is a version of $\nabla \zM tx_\tau$.
}

%Now, using Lemma \ref{lem:mal:ueps veps},
%\[\|\sup_{s \le r < T} \Ueps_r \|_2 \le C\int_s^{T-\varepsilon} \normexp{\aeps{r}} dr 
%\le C\int_s^{T-\varepsilon} \frac{ dr}{(T-r)^{(3-\alpha)/2} } \le C \phi(s,\varepsilon, \thetaL,\thetaC,\alpha),\]
%and
%$(\Veps_{j,t})^\top =  (\nabla\zeps_{j,t})^\top \sigma^{-1}(t,X_t) - \Ueps_t \nabla_x \sigma_j(t,X_t)$ from  Lemma \ref{lem:mal:ueps veps},
% we conclude that 
%$\normexp{\Veps_t} \le  \phi(s,\varepsilon, \thetaL,\thetaC,\alpha)$ as required.
% 
% it follows that 
% 
% Additionally, 
% and $\Veps_t =  D_t \zeps_t \sigma(t,X_t) - \Ueps_t \nabla_x \sigma(t,X_t) \; dt\times\P-a.e.$,,%+ C(T-t)^{(\alpha-2)/2}$.
%\qed

%Moreover, since the terminal condition of the BSDE satisfied by $(\yM tx, \zM tx)$ is zero, there is a continuous function $\zeps:[0,T)\times\R^d \rightarrow (\R^q)^\top$ given by
%\begin{equation}
%\zeps(t,x) = \E[\int_t^T F(r,\XM tx_r,\yM tx_r, \zM tx_r) \MW t{x,t}_r dr]
% \label{eq:zeps:fn}
%\end{equation}
%for all $(t,x)\in[0,T)\times\R^d$ and $\zeps_t = \zeps(t,X_t) $ for all $t\in [0,T)$ almost surely. 
%The proof of this can be found in \cite[Theorem 4.2]{ma:zhan:02}; the conditions of \cite[Theorem 3.1]{ma:zhan:02} are satisfied because the
%terminal condition of the BSDE $(\yeps,\zeps)$ is zero.

\vspace{0.3cm}
\noindent
{\bf Step 4. Proving $\zeps(t,\cdot)$ is Lipschitz continuous.}
%Let us now prove that $\zeps(t,\cdot)$ is Lipschitz continuous. 
Fix $s \in [t,T)$. Using the representation \eqref{eq:zeps:fn} of $\zM tx$, it follows that  
\begin{align*}
 \normexp{\zM t{x_1}_s - \zM t{x_2}_s} 
% |\zeps(t,x_1) - \zeps(t,x_2)|
& \le  \|\E_s[  \int_s^T  F(r,\XM t{x_1}_r,\yM t{x_1}_r,\zM t{x_1}_r) \MW t{x_1,s}_r dr]  \\
&\quad -  \E_s[\int_s^T  F(r,\XM t{x_2}_r,\yM t{x_2}_r,\zM t{x_2}_r) \MW t{x_1,s}_r dr ]\|_2 \\
&\qquad + \| \E_s[ \int_s^T  F(r,\XM t{x_2}_r,\yM t{x_2}_r,\zM t{x_2}_r) (\MW t{x_1,s}_r - \MW t{x_2,s}_r) dr ]\|_2\\
& =: \cA_{1} + \cA_{2}.
\end{align*}
We start with an estimate for $\cA_2$. Using the Cauchy-Schwarz inequality, 
% {(Corollary \ref{cor:cs:ineq})}, 
it follows that
\begin{align}
\cA_{2} & \le    \int_s^T  \|F(r,\XM t{x_2}_r,\yM t{x_2}_r,\zM t{x_2}_r)\|_4
\|\MW t{x_1,s}_r - \MW t{x_2,s}_r\|_4 dr 
\label{eq:A2:bd:1}
\end{align}

$\blacktriangleright$
{\bf Bounding $\|\MW t{x_1,s}_r - \MW t{x_2,s}_r\|_4$.}
Using the same techniques as in the proof of Lemma \ref{lem:mal:weight}, one shows that%but replacing the It\^o isometry 
%with the Burkholder-Davis-Gundy inequality,
\begin{align}
\|\MW t{x_1,s}_r - \MW t{x_2,s}_r\|_4 &  \le C_4 \frac{\normS{8}{\sigma^{-1}(s,\XM t{x_1}_s) 
- \sigma^{-1}(s,\XM t{x_2}_s)}
\E[\sup_{s\le u \le T} |D_s\XM t{x_1}_u |^8]^{1/8}  }{\sqrt{r-s}} 
\nonumber \\
& \qquad + C_4 \frac{\|\sigma^{-1}\|_\infty \E[\sup_{s\le u \le T} |D_s\XM t{x_1}_u - D_s\XM t{x_2}_u|^8]^{1/8}  }{\sqrt{r-t}} .
\label{eq:H:diff:bd1}
\end{align}
where $C_4$ is the constant coming from the BDG inequality. 
%Moreover, using that
%\[\E[\sup_{s\le u \le T} |D_s\XM t{x_1}_u - D_s\XM t{x_2}_u|^8] \le C \e[|\sigma(s,\XM t{x_1}_s) - \sigma(s,\XM t{x_2}_s)|^8]
%\le C \normS{8}{\XM t{x_1}_s - \XM t{x_2}_s}  \]
Thanks to \cite[Theorem IX.2.4]{revu:yor:01}, we have that
\begin{equation}
 \normS{8}{\XM t{x_1}_s - \XM t{x_2}_s}% +  \E[\sup_{s\le u \le T} |D_s\XM t{x_1}_u - D_s\XM t{x_2}_u|^8]^{1/8} 
 \le C|x_1 - x_2|
\label{eq:diff:lip} 
\end{equation}
The function $\sigma^{-1}(t,\cdot)$ is Lipschitz continuous uniformly in $t$
with Lipschitz constant as given in  Lemma \ref{lem:1:inv sigma:lip}
for all $s\in[t,T)$. Moreover, Lemma \ref{lem:mal:SDE:3} gives that
\[\E[\sup_{s\le u \le T} |D_s\XM t{x_1}_u |^8]^{1/8} \le C \quad \text{and}  \quad \E[\sup_{s\le u \le T} |D_s\XM t{x_1}_u - D_s\XM t{x_2}_u|^8]^{1/8} \le C|x_1 - x_2|.\]
%due to the uniform elliptic
%condition \He.
Combining these estimates, it follows that $\|\MW t{x_1,t}_r - \MW t{x_2,t}_r\|_4 \le C|x_1 -x_2|/\sqrt{r-t}$.

\vspace{0.3cm}
$\blacktriangleright$
{\bf Bounding $\|F(r,\XM t{x_2}_r,\yM t{x_2}_r,\zM t{x_2}_r)\|_4$.}
%In order to find a bound for , 
%We show that $\yM t{x_2}_r$ and $\zM t{x_2}_r$ are in $\L_4(\cF_r)$. 
%Once this has been shown, 
We take advantage of the local Lipschitz continuity  and boundedness \eqref{eq:loclip:driver}
of $f$,
and the uniform bounds on $u$ and its partial derivatives from Lemma \ref{lem:lin:pde:bds},
in order to show that
\begin{align}
 | F(r,\XM t{x_2}_r, 0 , 0 ) |& \le | f(r,\XM t{x_2}_r, 0 , 0 ) |+ L_f \frac{|u(r,\XM t{x_2}_r) |
+ \|\sigma\|_\infty | \nabla_x u(r,\XM t{x_2}_r) | }{(T-r)^{(1-\thetaL)/2} }
\nonumber 
\\&
 \le \frac{C_f}{(T-r)^{1-\thetaC} } + \frac{C{ B_r(\Phi)} }{(T-r)^{1-\thetaL/2} } 
 \le C { \frac{ B_r(\Phi) }{(T-r)^{1-\thetaC\wedge\frac\thetaL2} } }
\label{eq:1:bd F:0}
\end{align}
where
 {
 \[B_r(\Phi):= 
\left\{
\begin{array}{ll}
\sqrt{T-t} + C &\text{if $\Phi$ is constant,} \\
\E_r[|\Phi(\XM t{x_2}_T) -\E_r[\Phi(\XM t{x_2}_T)]|^2]^{1/2} + C & \text {else}.
\end{array}
\right.
\]
}
 {Without loss of generality, we will consider the setting where $\Phi$ is not constant, because, for constant $\Phi$, the arguments will be analogous to the arguments under \Hgh \ with $\thetaP \equiv 1$.}
It follows from the triangle inequality, the local Lipschitz continuity  \eqref{eq:loclip:driver} and the inequality \eqref{eq:1:bd F:0} that
\begin{align*}
|F(r,\XM t{x_2}_r,\yM t{x_2}_r,\zM t{x_2}_r)| & \le |F(r,\XM t{x_2}_r, 0 , 0 )|
 + L_f \frac{|\yM t{x_2}_r | + | \zM t{x_2}_r | }{(T-r)^{(1-\thetaL)/2} }
%  \nonumber \\
%& \le  \frac{C \| \Phi\|_\infty }{(T-r)^{(3-\thetaL)/2} } 
%+ L_f \frac{\|\yM t{x_2}_r \|_4 + \| \zM t{x_2}_r\|_4 }{(T-r)^{(1-\thetaL)/2} }.
%\label{eq:1:l4 f}
\end{align*}
%To show that show that $\yM t{x_2}_r$ and $\zM t{x_2}_r$ are in $\L_4(\cF_r)$, apply Lemma \ref{lem:1:apriori:cond}
But $\yM t{x_2}_r$ and $\zM t{x_2}_r$ are bounded in $\L_4$: 
applying  Proposition \ref{prop:1:apriori}
with $(Y_1,Z_1) = (0,0)$ and $(Y_2,Z_2) = (\yM t{x_2}, \zM t{x_2})$, combined with inequality \eqref{eq:1:bd F:0} and Lemma \ref{lem:integration:1} to obtain that
\begin{align}
|\yM t{x_2}_r| & \le C \int_r^{T-\varepsilon} \E_r[|F(u,\XM t{x_2}_u, 0, 0)|^2]^{1/2} du 
\le C B_r(\Phi) (T-r)^{\thetaC \wedge \thetaL/2} , \nonumber \\ 
|\zM t{x_2}_r| & \le C \int_r^{T-\varepsilon} \E_r[|F(u,\XM t{x_2}_u, 0, 0)|^2]^{1/2}  (u-r)^{-1/2} du 
\le  C  B_r(\Phi)  (T-r)^{ \thetaC\wedge \frac\thetaL2 -\frac12} 
\label{eq:z0}
\end{align}
for all $r \in [t,T)$.
Therefore, 
%\[ 
$\|F(r,\XM t{x_2}_r,\yM t{x_2}_r,\zM t{x_2}_r) \|_4 $
%\le
is bounded above by
%\frac{
$C (T-r)^{\thetaC\wedge\frac\thetaL2 -1 }$.
Now, both \Hge \ and \Hgh imply that $\|B_r(\Phi)\|_4$ is bounded above uniformly in $r$ by C.
%} \]  
% \frac{C \| \Phi\|_\infty }{(T-r)^{(3 - \thetaL)/2} } + C\| \Phi\|_{\infty}\frac{\int_r^{T-\varepsilon} (T-u)^{(\thetaL  - 3)/2} (u -r)^{-1/2} du }{ (T-r)^{(1-\thetaL)/2 } } .\]
% One can show, analogously to \eqref{eq:bd:f} in Corollary \ref{cor:mom bd},  that
% \[\normexp{F(r,\XM t{x_2}_r,\yM t{x_2}_r,\zM t{x_2}_r)} \le C(T-r)^{\thetaC\wedge(\thetaL/2)-1}.\]
%Combining the above estimates and Lemma \ref{lem:integration:1}, it follows that
Substituting this and the bound on $\|\MW t{x_1,t}_r - \MW t{x_2,t}_r\|_4$  into \eqref{eq:A2:bd:1}
\begin{align*}
\cA_2 & 
%\le  
%\int_t^T  \normexp{F(r,\XM t{x_2}_r,\yM t{x_2}_r,\zM t{x_2}_r)} \normexp{\MW t{x_1,t}_r - \MW t{x_2,t}_r} dr \\
%& 
\le C B_4(\Phi) |x_1 - x_2| \int_s^T \frac{dr}{(T-r)^{1- \thetaC\wedge\frac\thetaL2} \sqrt{r-s}}
 \le \frac{CB_4(\Phi) |x_1 - x_2|}{(T-s)^{\frac12- \thetaC \wedge \frac\thetaL2}}
\end{align*}

Now, we estimate $\cA_1$.
%Using the conditional Fubini's theorem, Lemma \ref{lem:cond:fubini}, and the Cauchy-Schwarz inequality,
 {Using  Corollary \ref{cor:cs:ineq} (with $\MW t{x_1,s}_r $ in the place of $H^s_r$), it follows that}
\begin{align*}
\cA_1 & \le  \| \int_s^{T-\varepsilon} (\E_s[| F(r,\XM t{x_1}_r,\yM t{x_1}_r,\zM t{x_1}_r) - F(r,\XM t{x_2}_r,\yM t{x_2}_r,\zM t{x_2}_r) |^2])^{1/2}
 (\E_s [|\MW t{x_1,s}_r |^2]) ^{1/2} dr \|_2 
\end{align*}
Analogously to Lemma \ref{lem:mal:weight},  
$(\E_s [|\MW t{x_1,s}_r |^2])^{1/2}\le C_M(r-s)^{-1/2}$, therefore 
%Minkowski's inequality implies
\begin{align*}
\cA_1 & 
 { \le  C \normexp{ \int_s^{T-\varepsilon} 
\frac{(\E_s[ |F(r,\XM t{x_1}_r,\yM t{x_1}_r,\zM t{x_1}_r) - F(r,\XM t{x_2}_r,\yM t{x_2}_r,\zM t{x_2}_r)|^2 ] )^{1/2} }{ \sqrt{r-s} }
dr } }
\\
&  { \le  C { \int_s^{T-\varepsilon} 
\frac{\normexp{ (\E_s[ |F(r,\XM t{x_1}_r,\yM t{x_1}_r,\zM t{x_1}_r) - F(r,\XM t{x_2}_r,\yM t{x_2}_r,\zM t{x_2}_r)|^2 ] )^{1/2} } }{ \sqrt{r-s} }
dr } }
\\
& =  C \int_s^{T-\varepsilon} 
\frac{\normexp{ F(r,\XM t{x_1}_r,\yM t{x_1}_r,\zM t{x_1}_r) - F(r,\XM t{x_2}_r,\yM t{x_2}_r,\zM t{x_2}_r) } }{ \sqrt{r-s} }
dr
\end{align*}
 {where we have  used Minkowski's inequality to take the norm $\normexp{\cdot }$ into the Lebesgue integral.}
By applying the Lipschitz continuity of $\feps(r,\cdot)$, $\cA_1$ is bounded by
\begin{align*}
&C \int_s^{T-\varepsilon} \frac{ \| \XM t{x_1}_r-\XM t{x_2}_r \|_2   }{ {(T-r)^{1-\theta_X/2}}\sqrt{r-t}}dr \\
& \quad + C \int_s^{T-\varepsilon} \frac{\| \sigma \|_\infty \| u(r,\XM t{x_1}_r) - u(r,\XM t{x_2}_r) \|_2 + \|\nabla_x\sigma\|_\infty \| \nabla_xu(r,\XM t{x_1}_r) - \nabla_xu(r,\XM t{x_2}_r) \|_2}
{(T-r)^{(1-\thetaL)/2}\sqrt{r-t}}dr \\
 & \qquad +  C \int_s^{T-\varepsilon} \frac{ \| \yM t{x_1}_r - \yM t{x_2}_r \|_2 +\| \zM t{x_1}_r -\zM t{x_2}_r\|_2 }
{(T-r)^{(1-\thetaL)/2}\sqrt{r-t}}dr
\end{align*}
Using the differentiability of $u(s,\cdot)$, it follows that
\begin{align*}
& \| u(r,\XM t{x_1}_r) - u(r,\XM t{x_2}_r) \|_2 
\le 
\|\cR( u,r,\XM t{x_1}_r,\XM t{x_2}_r)\|_2,
\\
& \| \nabla_xu(r,\XM t{x_1}_r) - \nabla_xu(r,\XM t{x_2}_r) \|_2 \le 
\|\cR(\nabla_x u,r,\XM t{x_1}_s,\XM t{x_2}_s)\|_2
\end{align*}
for all $r\in[t,T)$, where, for a differentiable function $g$, $\cR(g,r,x,x')$ is the remainder from the first order Taylor expansion of $g(r,x) -  g(r,x') $:
in the case of $g$ taking values in $\R$, this is equal to
\begin{equation}
\cR(g,r,x,x') = \{ \int_0^1 \nabla_x  g(\delta x + (1-\delta) x' ) d \delta \} (x -x');
\label{eq:taylor}
\end{equation}
in the multidimensional case, the expansion \eqref{eq:taylor} is defined component-wise.
Denote by $\cR(r)$ the sum of the normed residuals 
%\[ \cR(t) :=
$ \|\cR( u,r,\XM t{x_1}_r,\XM t{x_2}_r)\|_2+ \|\cR(\nabla_x u,r,\XM t{x_1}_s,\XM t{x_2}_s)\|_2 $.
Therefore, using the notation 
 $\Theta_r := \normexp{\yM t{x_1}_r - \yM t{x_2}_r} +\normexp{\zM t{x_1}_r -\zM t{x_2}_r}$,
the final bound on $\cA_1$ is
\[
 \cA_1 \le 
  {C \int_s^{T-\varepsilon} \frac{ \| \XM t{x_1}_r-\XM t{x_2}_r \|_2  dr }{ {(T-r)^{1-\theta_X/2}}\sqrt{r-t}}}
 +C \int_s^{T-\varepsilon} \frac{\cR(r) dr}{(T-r)^{(1-\thetaL)/2}\sqrt{r-t}}
+ C\int_s^{T-\varepsilon} \frac{\Theta_r dr }
{(T-r)^{(1-\thetaL)/2}\sqrt{r-t}}.
\]
It follows from  Lemma \ref{lem:lin:pde:bds} that
\[
\cR(r) \le
\left\{
\begin{array}{rl}
C \|\Phi\|_\infty \normexp{\XM t{x_1}_r - \XM t{x_2}_r } (T-r)^{-1} \le C \|\Phi\|_\infty |x_1 - x_2| (T-r)^{-1}    & \text{under \Hge,} \\
& \\
C K_\Phi  \normexp{\XM t{x_1}_r - \XM t{x_2}_r } (T-r)^{\frac\thetaP2 -1 } \le C K_\Phi |x_1 - x_2|(T-r)^{\frac\thetaP2-1} & \text{under \Hgh},
\end{array}
\right.
\] 
The bound $\|\XM t{x_1}_r - \XM t{x_2}_r\|_2 \le C|x_1 - x_2|$ is obtained from \cite[Theorem IX.2.4]{revu:yor:01},
 {which also implies that
\[
\int_s^{T-\varepsilon} \frac{ \| \XM t{x_1}_r-\XM t{x_2}_r \|_2   }{ {(T-r)^{1-\theta_X/2}}\sqrt{r-t}}dr \le \int_s^{T-\varepsilon} \frac{C |x_1-x_2|   }{ {(T-r)^{1-\theta_X/2}}\sqrt{r-t}}dr.
\]
}
It is clear from the bounds above that the integral in $\cR(r)$ in the bound of $\cA_1$ dominates the upper bound on $\cA_2$, and  {also the integral
\[
\int_s^{T-\varepsilon} \frac{ \| \XM t{x_1}_r-\XM t{x_2}_r \|_2   }{ {(T-r)^{1-\theta_X/2}}\sqrt{r-t}}dr.
\]}
Therefore,
\begin{align*}
\normexp{\zM t{x_1}_s - \zM t{x_2}_s}
% |\zeps(t,x_1) - \zeps(t,x_2)|
 &\le %\frac{C|x_1 - x_2|}{(T-t)^{(1-(2\thetaC)\wedge\thetaL)/2} } + C|x_1 - x_2|(T-t)^{\thetaL/2} \\ 
 C \int_s^{T-\varepsilon} \frac{\cR(r) dr}{(T-r)^{(1-\thetaL)/2}\sqrt{r-t}}
+ C\int_s^{T-\varepsilon} \frac{\Theta_r }
{(T-r)^{(1-\thetaL)/2}\sqrt{r-t}}dr.
\end{align*}
Since 
%$\yM tx_s = \e_s[\int_s^T \feps(r,\XM tx_r,\YM tx_r,\ZM tx_r)dr ]$, 
$\yM tx_s = \e_s[\int_s^T F(r,\XM tx_r,\yM tx_r,\zM tx_r)dr ]$
similar estimates yield
\[
 \Theta_s \le C   \int_s^{T-\varepsilon} \frac{\cR(r) dr}{(T-r)^{(1-\thetaL)/2}\sqrt{r-s}}
+ C\int_s^{T-\varepsilon} \frac{\Theta_r}{(T-r)^{(1-\thetaL)/2}\sqrt{r-s} }dr
\]
for all $s \in [t,T-\varepsilon)$.
Let \Hge \ be in force.
Applying Lemma \ref{lem:iteration:2} with 
\[w_r := C  \|\Phi\|_\infty  |x_1 - x_2|\int_r^{T-\varepsilon} (T-u)^{(\thetaL-3)/2}(u-r)^{-1/2} du\] and $u_r := \Theta_r$, it follows that
\begin{align*}
 \Theta_s &\le 
C  \|\Phi\|_\infty |x_1 - x_2| \int_s^{T-\varepsilon} \frac{dr}{(T-r)^{(3-\thetaL)/2}\sqrt{r-s}} 
\\
& \quad
+ C  \|\Phi\|_\infty |x_1 - x_2| \int_s^{T-\varepsilon} \frac{ 
\int_s^{u} (u-r)^{\thetaL/2 -1}(r-s)^{-1/2} dr
}{(T-u)^{(3-\thetaL)/2} } du 
 + C\int_s^{T-\varepsilon} \frac{\Theta_r}{(T-r)^{(1-\thetaL)/2} } dr \\
& \le C  \|\Phi\|_\infty  |x_1 - x_2| \int_s^{T-\varepsilon} \frac{dr}{(T-r)^{(3-\thetaL)/2}\sqrt{r-s}} 
+ C \int_s^{T-\varepsilon} \frac{\Theta_r}{(T-r)^{(1-\thetaL)/2} } dr,
\end{align*}
where we have used Lemma \ref{lem:integration:1} to bound the integral $\int_s^{u} (u-r)^{\thetaL/2 -1}(r-s)^{-1/2} dr$.
%by $C (u-s)^{(\thetaL -1)s/2}$.
Then, applying Lemma \ref{lem:iteration:3} to bound the integral $\int_s^{T-\varepsilon} {\Theta_r}{(T-r)^{(\thetaL-1)/2} } dr$,
final bound on $\normexp{\zM t{x_1}_s - \zM t{x_2}_s}$ for all $t \in [0,T)$, $(x_1,x_2) \in (\R^d)^2$, and $s \in [t,T)$ is
\begin{align}
 \normexp{\zM t{x_1}_s - \zM t{x_2}_ s} & \le 
 C  \|\Phi\|_\infty |x_1 - x_2| \int_s^{T-\varepsilon} \frac{dr}{(T-r)^{(3-\thetaL)/2}\sqrt{r-s}}  \nonumber \\
& \qquad + C  \|\Phi\|_\infty |x_1 - x_2|
\int_s^{T-\varepsilon} \frac{
\int_r^{T-\varepsilon} (T-u)^{(\thetaL-3)/2}(u-r)^{-1/2} du
}{(T-r)^{(1-\thetaL)/2} } dr \nonumber 
%\\
%& \le C \phi(t,\varepsilon,\thetaL) |x_1 - x_2|. % \int_s^{T-\varepsilon} (T-r)^{(\thetaL-3)/2}(r-s)^{-1/2} dr.
%\label{eq:1:zeps:lip const}
\end{align}
 and application of Lemma \ref{lem:integration:1} yields the final upper bound $ C \phi(t,\varepsilon,\thetaL) |x_1 - x_2|$.
 Therefore, setting $t =s$ yields that the function $\zeps(t,\cdot)$ is Lipschitz continuous with Lipschitz constant $C \phi(t,\varepsilon,\thetaL)$, as required. The proof under \Hgh \ is analogous.
\qed

\vspace{0.3cm}
 {
We now come to the regularity result advertised at the beginning of this section; this result is not used in the remainder of this paper, but may hold some interest for other works.
\begin{corollary}
\label{cor:proxy:lip}
Let \Hgh \ and \HFd \ be in force, and  let $\thetaL + \thetaP \ge 1$.
Then there exists a function $z : [0,T) \times \R^d \to (\R^q)^\top$ such that $\nabla_x u(s,X_s)\sigma(s,X_s) + z(s,X_s)$ is a version of $Z_s$.
Moreover, recalling the function $ \phi(t,\varepsilon,\thetaL,\thetaP)$ from \eqref{eq:phi:3} for all $t \in [0,T)$, $x \mapsto z(t,x)$ is Lipschitz continuous with Lipschitz constant equal to 
\[\lim_{\varepsilon \to 0}\phi(t,\varepsilon,\thetaL,\thetaP) \le {C \over (T-t)^{1 - (\thetaL + \thetaP)/2} } \]
for some finite constant $C$ depending only on $K_\Phi$, $L_f$, the bounds on $b$ and $\sigma$ and their partial derivatives, $\elip$, $C_M$, $\thetaL$, $\thetaC$, $C_f$, and $T$.
\end{corollary}
{\bf Proof.}
Let $(Y^{(t,x)},Z^{(t,x)})$ be the solution of 
\[
Y^{(t,x)}_s = \Phi(\X{t,x}T) + \int_s^T f(\tau,\X{t,x}\tau,Y^{(t,x)}_\tau,Z^{(t,x)}_\tau) d\tau - \int_s^T Z^{(t,x)}_\tau dW_\tau,
\]
and set
\[
z(t,x) := \E_t[\int_s^T f(\tau,\X{t,x}\tau,Y^{(t,x)}_\tau,Z^{(t,x)}_\tau) \MW t{x,s}_\tau d\tau]
\]
for
\[
\MW t{x,s}_r := \frac{\1_{(t,T]}(s)}{r-s} (\int_s^r \sigma^{-1}(r,\XM tx_r) D_s \XM tx_r )^\top dW_r)^\top\] 
Recall the function $\zeps : [0,T) \times \R^d \to (\R^q)^\top$ from \eqref{eq:zeps:fn}.
One shows  $z(t,x) = \lim_{\varepsilon\to 0} \zeps(t,x) $ by mimicking the proof of Theorem \ref{thm:repr}.
Since $Z$ is the limit of $\Zeps$ as $\varepsilon \to 0$ in $\cH^2$, 
and $\zeps(s,X _s)$ is a version of $\Zeps_s -\nabla_x u(s,X_s)\sigma(s,X_s)$,
it follows that $z(s,X_s)$ is a version of $Z_s -\nabla_x u(s,X_s)\sigma(s,X_s)$, as required.
Finally, to prove the Lipschitz continuity of $z(t,\cdot)$, we observe that, for $\thetaL + \thetaP \ge 1$,
\begin{equation*}
\lim_{\varepsilon \to 0} \phi(t,\varepsilon, \thetaL, \thetaP) = 
 K_\Phi \int_t^{T}\frac{dr}{(T-r)^{(3-\thetaL - \thetaP
)/2} \sqrt{r-t} } \le {C \over (T-t)^{1 - (\thetaL + \thetaP)/2} }
\end{equation*}
 thanks to Lemma \ref{lem:integration:1},
and proceed as in {\bf Step 4} of the proof of Proposition \ref{prop:veps:bd} (with $z(t,\cdot)$ in the place of $\zeps(t,\cdot)$);
the upper bound on the limit $\lim_{\varepsilon \to 0} \phi(t,\varepsilon, \thetaL, \thetaP) $  comes from Lemma \ref{lem:integration:1}. \qed
}

\vspace{0.3cm}

In order to make use of Proposition \ref{prop:veps:bd}, 
it is is necessary 
%decompose the $\L_2$-regularity of 
to approximate $Z$
% into the $\L_2$-regularity of an 
 by an intermediate process $Z_M$ which satisfies the hypotheses of Proposition \ref{prop:veps:bd}.
\begin{lemma}
\label{lem:1:mol:proxy}
Assume that \Hgexp \ is in force.
Recall the BSDE $(Y_M,Z_M)$ defined in Corollary \ref{cor:1:moll}.
Take the version of $Z_M$ given by Theorem \ref{thm:repr}.
%satisfying $Z_{M,t} = \E_t[ \Phi_M(X_T) H^t_T + \int_t^T f_M(s,X_s,Y_{M,s}, Z_{M,s} ) H^t_s ds ]$ for all $t\in[0,T)$ almost surely.
For $M =   (3 \ln(N) )^{1/4}$ and $R(M)$ equal to $3 L_f e^{M^{2}/2}  $, there is a finite constant $C$ depending only on $L_f$, $C_M$, $\thetaL$,  $C_\xi$ and $T$, but not on $N$,
such that for all $ N\ge1$
\begin{align*}
\sum_{i=0}^{N-1} \int_{t_i}^{t_{i+1}} \normexp{Z_s - \tilde Z_{t_i} }^2 ds \le 
C \sum_{i=0}^{N-1} \int_{t_i}^{t_{i+1}} \normexp{Z_{M,s} - \tilde Z_{M,t_i} }^2 ds + C N^{-1}
\end{align*}
 {where
$ \tilde Z_{t_i} := \frac{1}{\Delta_i} \E_i \big[\int_{t_i}^{t_{i+1}} Z_t dt \big]$ and $ \tilde Z_{M,t_i} := \frac{1}{\Delta_i} \E_i \big[\int_{t_i}^{t_{i+1}} Z_{M,t} dt \big]$.}
\end{lemma}

 {\bf Proof.}
 In what follows, $C$ may change from line to line.
Using Cauchy's inequality and the orthogonality of the projections,
 \[
 %\begin{align*}
 \frac 12 \sum_{i=0}^{N-1} \int_{t_i}^{t_{i+1}} \normexp{Z_s - \tilde Z_{t_i} }^2 ds 
% & 
 \le  \int_{0}^{T} \normexp{Z_s -   Z_{M,s} }^2 ds 
  +  \sum_{i=0}^{N-1}  \normexp{\tilde Z_{t_i} -  \tilde Z_{M,t_i} }^2 \Delta_i 
%  \\
%  &\qquad  + \sum_{i=0}^{N-1} \int_{t_i}^{t_{i+1}} \normexp{Z_{M,s} - \bar Z_{M,s} }^2 ds 
%  +  \sum_{i=0}^{N-1}  \normexp{Z_{M,t_i} - \bar Z_{M,t_i} }^2 \Delta_i \\
%  & \qquad \qquad 
  + \sum_{i=0}^{N-1} \int_{t_i}^{t_{i+1}} \normexp{Z_{M,s} - \tilde Z_{M,t_i} }^2 ds.
  \]
From Jensen's inequality, it follows that $ \sum_{i=0}^{N-1}  \normexp{\tilde Z_{t_i} -  \tilde Z_{M,t_i} }^2 \Delta_i \le \int_{0}^{T} \normexp{Z_s - \bar  Z_{M,s} }^2 ds 
$.
Proposition \ref{prop:1:apriori} with $(Y_1,Z_1) := (Y_M,Z_M)$ and $(Y_2,Z_2) := (Y,Z)$ yields that, for any $s\in[0,T)$,
\begin{align}
\int_0^T \normexp{Z_t - Z_{M,t}}^2 dt \le C \normexp{\Phi(X_T) - \Phi_M(X_T)}^2 + C\int_0^T \normexp{f(t,X_t, Y_{M,t},  Z_{M,t})  - f_M(t,X_t,  Y_{M,t},  Z_{M,t} )}^2 dt ,
\label{eq:moll:diff}
\end{align}
It follows from Markov's exponential inequality and \Hgexp \ that 
\begin{equation}
\normexp{\Phi(X_T) - \Phi_M(X_T) }^2 = \int_{M^2}^\infty \P( |\Phi(X_T)|^2 \ge x) dx \le C_\xi \int_{M^2}^\infty e^{-\sqrt{x} } dx
= 2 C_\xi (1+M^4) e^{-M^4}
\le CN^{-2}.
\label{eq:moll:term}
\end{equation}
The last inequality is obtained by substituting the value of $M$.
On the other hand, the basic properties of the mollifier in Definition \ref{def:1:moll} yields
\begin{align}
& | f(t,X_t, Y_{M,t},  Z_{M,t})  - f_M(t,X_t,  Y_{M,t},  Z_{M,t} )  | 
\nonumber\\
&  \le \int_{\R^d \times \R \times (\R^q)^\top } | f(t,X_t, Y_{M,t},  Z_{M,t} )  - f(t,X_t - x, Y_{M,t} - y ,  Z_{M,t} - z )  | \phi_{R(M)}(x, y,  z)  d(x,y,z) 
\nonumber \\
& \le \frac{L_f }{(T-t)^{(1-\thetaL)/2} } \int_{\{ |x|^2 + |y|^2 + |z|^2 \le R(M)^{-2} \} }( |x| + |y| + |z| ) \phi_{R(M)}(x, y,  z)  d(x,y,z) 
%\\& 
\le \frac{3 L_f }{(T-t)^{(1-\thetaL)/2} R(M) }
\label{eq:moll:driver}
\end{align}
Substituting the value of $R(M)$ then gives $\normexp{ f(t,X_t,\bar Y_{M,t} , \bar Z_{M,t})  - f_M(t,X_t,\bar Y_{M,t}, \bar Z_{M,t}  } \le (T-t)^{(\thetaL-1)/2}N^{-1/2}$
for all $t \in [0,T)$.
%Applying Proposition \ref{prop:1:apriori} with $(Y_1,Z_1) = (\bar Y_M,\bar Z_M)$ and $(Y_2,Z_2) = (Y_M,Z_M)$,
%the bound on $\normexp{ f(t,X_t,\bar Y_{M,t} , \bar Z_{M,t})  - f_M(t,X_t,\bar Y_{M,t}, \bar Z_{M,t}  }$, 
Substituting \eqref{eq:moll:term} and \eqref{eq:moll:driver} into \eqref{eq:moll:diff} Lemma \ref{lem:integration:1} then yields
\[
\int_{0}^{T} \normexp{Z_s -  Z_{M,s} }^2 ds 
 \le CN^{-1} + CN^{-2} \sum_{i=0}^{N-1} \frac{\Delta_i}{T-t_i},
\]
The sum on the right hand side above is bounded by $1 + \int_0^{t_{N-1}} (T-t)^{-1} dt = 1 + C\ln(N)$,
whence the proof is complete.
%bound the terms in $Z_M$ and $\bar Z_M$ in the decomposition by $CN^{-1}$.
%it follows that
%% there is a (possibly different) finite $C \ge 0$ depending only on $L_f$, $\thetaL$ $C_M$, and $T$ such that
%\begin{align*}
%\int_0^T \normexp{Z_{M,t} - \bar Z_{M,t} }^2 dt & \le C \int_0^T \normexp{ f(t,X_t,\bar Y_{M,t} , \bar Z_{M,t})  - f_M(t,X_t,\bar Y_{M,t}, \bar Z_{M,t}  } ^2dt, 
%\\
%\normexp{Z_{M,t_i} - \bar Z_{M,t_i} } & \le  C \int_{t_i}^T \frac{  \normexp{ f(t,X_t, \bar Y_{M,t}, \bar Z_{M,t} )  - f_M(t,X_t,\bar Y_{M,t} ,\bar Z_{M,t})  } }{ \sqrt{t-t_i} } dt \quad \forall i.
%\end{align*}
%%For $M \ge M_0 := \inf \{ i \in \N \ : \ i> \|\Phi(X_T) \}$, $\normexp{ \Phi(X_T) - \Phi_M(X_T)  }$ is equal to $0$.
%%Using the definition of the function $\phi_{R,M}$, in Definition \ref{def:1:moll}, and the definition of $f_M$, it follows that
%Setting   the function  $R(t,M)$ to be equal to $3 L_f e^{M^{4}/2} (T-t)^{(\thetaL -1)/2} $ ensures that
%$\int_0^T \normexp{Z_{M,t} - \bar Z_{M,t}}^2 dt \le CN^{-1}$ and $\normexp{Z_{M,t_i} - \bar Z_{M,t_i} } \le CN^{-1/2}(T-t_i)^{1/2}$, whence the result follows.
\qed

\vspace{0.3cm}
%We now come to the main result of this section, where a convergence rate for the $\L_2$ regularity is computed 
We now provide an extension to Theorem \ref{prop:1:l2 reg:makh} under \Hgh \ with the aid of Proposition \ref{prop:veps:bd}.
\begin{theorem}
\label{thm:L2 reg}
Let  \Hgh \ be in force and $0<\beta < (2\gamma)\wedge(\alpha\wedge\thetaL)$.
There is a constant $C$ depending only on $L_f$, $C_M$, $\thetaL$, $\thetaC$, $\beta$, $C_f$, $K_\Phi$ and $T$, but not on $N$,
such that for all $ N \ge 1$, 
\begin{equation}
%\cE(\pib_N) \le  \sum_{i=0}^{N-1} \int_{t_i}^{t_{i+1}} \normexp{Z_s - Z_{t_i}}^2 ds 
\cE(N)
\le CN^{-1} \1_{[1,3]}(\thetaP + \beta + 2\gamma) + CN^{-2\gamma} \1_{(0,1)}(\thetaP + \beta + 2\gamma)
\label{thm:l2 reg:z}
\end{equation}
 {for $\gamma := (\thetaC\wedge\frac\alpha2 + \frac\thetaL2)\wedge\thetaC $.}
%where $\delta_N \ge \1_{[3,\infty)}(N) \ln\ln(N) / \ln(N)$.
\end{theorem}
{\bf Proof.}
%From Proposition \ref{}, it is sufficient to bound $ $.
In what follows,
 $C$ may change from line to line.
To start with, we assume that \HFd \ and \Hge \ are in force.
Recall \eqref{eq:zeps:bd:1}. From the bounds $\normexp{\aeps{r}} \le C(T-r)^{(\alpha+\thetaL -3)/2}$ 
in the proof on Lemma \ref{lem:mal:ueps veps}
the first sum $\sum_{i=0}^{N-1} \int_{t_{i}}^{t_{i+1}} (\int_{t_i}^{t} \normexp{\aeps{r}} dr )^2 dt$
is bounded above by
\begin{align}
%&C\sum_{i=0}^{N-1} \int_{t_{i}}^{t_{i+1}} 
%\{ 
%(\int_{t_i}^{t} \normexp{\aeps{r}} dr )^2 dt \\
& C \sum_{i=0}^{N-2} \int_{t_{i}}^{t_{i+1}} \frac{ \big( \int_{t_i}^{t}\frac{dr}{(t-r)^{(1-\thetaL)/2} } \big)^2 }{(T-t)^{2-\alpha} } dt
+C \int_{t_{N-1}}^{T} \frac{ \big( \int_{t_{N-1}}^{t}\frac{dr}{(t-r)^{1-\thetaL/2} } \big)^2 }{(T-t)^{1-\alpha} } dt 
\nonumber \\
& \le C\sum_{i=0}^{N-2} \int_{t_{i}}^{t_{i+1}} \frac{(T-t)^{\beta-1} (t-t_{i})^{1+\thetaL} }{(T-t)^{1 + \beta-\alpha} } dt
+C\int_{t_{N-1}}^{T} \frac{(t-t_{N-1})^{\thetaL} }{(T-t)^{1-\alpha} } dt \nonumber\\
&\le  C\sum_{i=0}^{N-2} \frac{\Delta_{i}^{1+\thetaL} }{(T-t_{i+1})^{1-\beta}} \int_{t_{i}}^{t_{i+1}} \frac{ dt }{(T-t)^{1+\beta-\alpha} } 
+C\Delta_{N-1}^{\thetaL+\alpha} .
\label{eq:aeps:bd:1}
\end{align}
Using \eqref{eq:grid bd:2} from Lemma \ref{lem:inc growth}, ${\Delta_{i}} \le C{\Delta_{i+1}}$ for $i<N-1$, which,
combined with \eqref{eq:grid bd:1}, yields
\[\max_{0\le i\le N-2}\frac{\Delta_{i}^{1+\thetaL} }{(T-t_{i+1})^{1-\beta}} 
\le C \max_{0\le i\le N-1}\frac{\Delta_{i}^{1+\thetaL} }{(T-t_{i})^{1-\beta}} <CN^{-1-\thetaL}. \] 
Additionally, $\beta < \alpha$ implies that  $\Delta_{N-1}^{\alpha + \thetaL} = C N^{-(\alpha+\thetaL)/\beta} \le CN^{-1}$. 
Substituting these results into \eqref{eq:aeps:bd:1} gives 
\begin{equation}
\sum_{i=0}^{N-1} \int_{t_{i}}^{t_{i+1}} (\int_{t_i}^{t} \normexp{\aeps{r}} dr )^2 dt 
\le CN^{-1}.
\label{eq:aeps:bd:2}
\end{equation}
The refined estimates
 { -- $\normexp{\Veps_r} \le 
C\phi(r,\varepsilon, \thetaL, \thetaP
) $  for  all $r \in [0,T)$ --}
from Proposition \ref{prop:veps:bd}
are used to bound $\sum_{i=0}^{N-1} \int_{t_{i}}^{t_{i+1}} (\int_{t_i}^{t} \normexp{\Veps_{r}} dr )^2 dt$ from above  by 
$\sum_{i=0}^{N-1} \int_{t_{i}}^{t_{i+1}} \Big( \int_{t_{i}}^{t}\phi(r,\varepsilon,\thetaL,\thetaP
) dr \Big)^2 dt 
$, which itself is bounded above by
\begin{align}
 \sum_{i=0}^{N-1} \int_{t_{i}}^{t_{i+1}} & \Big( \int_{t_{i}}^{t}\phi(r,\varepsilon,\thetaL,\thetaP
% ,\thetaC,\alpha
) dr \Big)^2 dt 
%\nonumber \\&
 { \ =  \sum_{i=0}^{N-1} \int_{t_{i}}^{t_{i+1}} \Big( \int_{t_{i}}^{t} \Big\{ K_\Phi \int_r^{T-\varepsilon} \frac{du}
{(T-u)^{\frac{3  -\thetaP -\thetaL}{2}}\sqrt{u-r}} \Big\} dr \Big)^2 dt } \nonumber \\
&  {= K_\Phi^2 \sum_{i=0}^{N-1} \int_{t_{i}}^{t_{i+1}} \Big( \int_{t_{i}}^{t}\Big\{\int_r^{T-\varepsilon} \frac{(T-u)^{(\beta-1)/2}du}
{(T-u)^{\frac{ 2 + \beta -\thetaP-\thetaL}{2}}\sqrt{u-r}} \Big\} dr \Big)^2 dt } \nonumber \\
&\le 
 \frac{K_\Phi^2}{\varepsilon^{1-\thetaP-\beta} }\sum_{i=0}^{N-1} \int_{t_{i}}^{t_{i+1}} \Big( \int_{t_{i}}^{t}\Big\{\int_r^{T-\varepsilon} \frac{du}
{(T-u)^{1-(\thetaL-\beta)/2}\sqrt{u-r}} \Big\} dr \Big)^2 dt .\nonumber 
\end{align}
Now, using Lemma \ref{lem:integration:1}  to obtain  {an upper bound ${(T-r)^{-\big(1- (\thetaL-\beta) \big)/2}}$ on the inner integral $\int_r^{T-\varepsilon} 
{(T-u)^{-\big(1-(\thetaL-\beta)/2 \big)}({u-r})^{-1/2}} du$},
\begin{align}
 \sum_{i=0}^{N-1} \int_{t_{i}}^{t_{i+1}} & \Big( \int_{t_{i}}^{t}\phi(r,\varepsilon,\thetaL,\thetaP
% ,\thetaC,\alpha
) dr \Big)^2 dt 
\le  \frac{C}{\varepsilon^{1-\thetaP-\beta} }\sum_{i=0}^{N-1} \int_{t_{i}}^{t_{i+1}} 
\Big( \int_{t_{i}}^{t} \frac{ dr}{(T-r)^{\big(1- (\thetaL-\beta) \big)/2}} \Big)^2 
dt \nonumber \\
&\le \frac{C  (\max_{0\le i \le N-1}\Delta_i)^2}{\varepsilon^{1- \thetaP - \beta} } \int_{0}^{T}  \frac{ dr}{(T-t)^{1-(\thetaL-\beta)}} dt 
\le \frac{C  N^{-2}}{\varepsilon^{1-\thetaP-\beta}}
\label{eq:veps:bd:2} 
\end{align}
 {where we have used  Jensen's inequality to get
\[
\Big( \int_{t_{i}}^{t} \frac{ dr}{(T-r)^{\big(1- (\thetaL-\beta) \big)/2}} \Big)^2 \le C \int_{0}^{T}  \frac{ dr}{(T-t)^{1-(\thetaL-\beta)}} dt 
,\]
}
and then \eqref{eq:grid bd:1} in Lemma \ref{lem:inc growth} for the bound $\max_{0\le i \le N-1}\Delta_i\le CN^{-1}$.
Substituting \eqref{eq:aeps:bd:2} and \eqref{eq:veps:bd:2} into \eqref{eq:zeps:bd:1} finally yields
\begin{equation}
\sum_{i=0}^{N-1} \int_{t_i}^{t_{i+1}} \normexp{\zeps_s - \zeps_{t_i} }^2 ds \le CN^{-1} + \frac{C N^{-2}}{\varepsilon^{1-\thetaP-\beta}}.
\label{eq:zeps:bd:2}
\end{equation}
Then, using $\Zeps = \zeps + z$, Lemma \ref{lem:Z -Zeps:L2}, and  
$\sum_{i=0}^{N-1} \int_{t_i}^{t_{i+1}} \normexp{z_s- z_{t_i} }^2 ds \le CN^{-1}$ as shown in \cite[Theorem 1.3]{gobe:makh:10},
it follows that
\begin{align}
\sum_{i=0}^{N-1} \int_{t_i}^{t_{i+1}} \normexp{Z_s - Z_{t_i}}^2 ds & \le  CN^{-1} + 
{C} {N^{-2}\varepsilon^{\thetaP + \beta-1}} + C  {N^{-2\gamma/\beta}} + C\varepsilon^{2\gamma}.%\big(1 + \ln(N) \vee 1 \big).
\label{eq:cE1:eps:2}
\end{align}
%Let $\delta \in [0,1]$ be sufficiently large so that $N^{-\delta}\ln(N) \le 1$; 
For $N \in\{1,2\}$, 
let $\delta := 0$, and for $N > 2$, $\delta := \ln\ln(N)/\ln(N)$.
%notice that the value of $\delta$ to be $\delta_N$ given in the theorem statement.
Set $\varepsilon = N^{-(1+\delta)/(2\gamma)}$.
Recalling further that $2\gamma < \beta$, this implies that, under \Hge \  and \HFd,
\[
\sum_{i=0}^{N-1} \int_{t_i}^{t_{i+1}} \normexp{Z_s - Z_{t_i}}^2 ds \le  CN^{-1} + 
{C} N^{-2 - (\thetaP +\beta-1)(1+\delta_N)/(2\gamma) } .
\]
To obtain the general result, recall the BSDE $(Y_M,Z_M)$ from Corollary \ref{cor:1:moll}.
The driver of $(Y_M,Z_M)$ satisfies assumptions \HFd.
%, and it follows that
%\[
%\sum_{i=0}^{N-1} \int_{t_i}^{t_{i+1}} \normexp{Z_{M,s} - Z_{M,t_i}}^2 ds  \le  CN^{-1} + C M N^{-2 + (\beta-1)(1+\delta_N)/(2\gamma) }.
%\]
The proof is complete by taking $M$ equal to $(3\ln(N))^{1/4}$,  $R(M)$ equal to $3 L_f e^{M^{4}/2} $,  and applying Lemma \ref{lem:1:mol:proxy}.
\qed

\section{Convergence rate of the Malliavin weights scheme}
\label{section:num}
 {In this section, we treat the Malliavin weights scheme
\begin{align*}
 \bY NN & := \Phi(X_T), \quad \bY N{i} :=  \e_i[\Phi(X_{T}) + \sum_{j=i}^{N-1}f(t_j, X_{t_j},\bY N{j+1}, \bZ N{j})(t_{j+1} - t_j)], \nonumber \\
 \bZ N{i} & := \e_i[\Phi(X_{T})H^{i}_{N} + \sum_{j=i+1}^{N-1} f(t_j , X_{t_j},\bY N{j+1}, \bZ N{j})H^{i}_{j}(t_{j+1} - t_j)]
%\label{eq:1:disc:malscheme}
\end{align*}
}
Recall the
% flow process $\nabX{}$ and its inverse $\inX{}$ defined by $\nabXM 0{x_0}{}$ and $\inXM t{x_0}{}$, respectively, 
Malliavin derivative of the the marginals of the process $X$ in Section \ref{section:mal:SDE}.
In the definition of the Malliavin weights scheme \eqref{eq:1:disc:malscheme}, we use the following discrete-time approximation 
%-  built with the Markov chain approximations of $X$, $\nabX{ }$ and $\inX{ }$ -
 of the Malliavin weights \eqref{eq:malweights}:
\begin{equation}
H^i_j := \frac{1}{t_j-t_i} \big( \sum_{k=i}^{j-1} %(\sigma^{-1}(t_k,X_{t_k})
% \nabla X_{t_k} \nabla X_{t_i}^{-1} 
D_{t_i} X_{t_k}
 \sigma(t_i,X_{t_i}) )^\top \Delta W_k \big)^\top
\label{eq:disc:malweights}
\end{equation}
Notice that $H^i_j$ satisfies $\E_i[H^i_j] = 0 $ and $\E_i[|H^i_j|^2] \le C_M (t_j-t_i)^{-1}$;
the latter property is proved exactly like Lemma \ref{lem:mal:weight}.
If the marginals of $X$ and $D_{t_i} X$ are not known explicitly, one can use an SDE scheme to provide approximations, but this is beyond the scope of this work;
 {some work has been done on this in the zero driver case ($f\equiv0$), in particular we refer the reader to Section 3 (and the sequel) of \cite{gobe:muno:05}.}
 { In what follows, we use the version of $Z$ given by Theorem \ref{thm:repr}, in other words 
 \[
 Z_t  = \e_t[\Phi(X_T) H^{t}_T + \int_{t}^T f(s,X_s,Y_s,Z_s) H^{t}_s ds] \quad 
\text{for all } t\in[0,T) \quad \P-a.s.
 \]
 }
We start with some preliminary results.

\begin{lemma}
\label{lem:mw:diffs}
There is a constant $C$ depending only on the bound on $b$ and it's derivatives, the bound on $\sigma$ and it's derivatives, $\elip$, $L_f$, $\thetaL$, $C_f$, $\thetaC$, $K^\alpha(\Phi)$ and $T$ such that, for any $0 \le i < j \le N$, 
\begin{align*}
\normexp{\E_i[\Phi(X_T) (H^{t_i}_{t_j} - H^i_j)]} & \le \frac{CN^{-1/2}}{(T-t_i)^{(1-\alpha)/2}} , \\
 \normexp{\E_i[f_j(X_{t_j},Y_{t_{j+1}},Z_{t_j}) (H^{t_i}_{t_j} - H^i_j)]} & \le \frac{CN^{-1/2}}{(T-t_j)^{1-\gamma}\sqrt{t_j-t_i}}
\end{align*}
 {where $\gamma := (\thetaC\wedge\frac\alpha2 + \frac\thetaL2)\wedge\thetaC $.}
\end{lemma}
{\bf Proof.}
For any $j >i$ and $t\ge t_{j}$, define
$N^{t_{i}}_{t} := \sigma^{-1}(t,X_{t}) D_{t_i} X_t \sigma(t_i,X_{t_{i}})$.
%\text{ and } \quad N^{i}_{j} := \sigma^{-1}(t_{j},X_{j}) \nabX{j} \inX{i} \sigma(t_i,X_{i}).\] 
%Observe the  decompositions into telescopic sums
Using the decomposition
\begin{align*}
N^{t_{i}}_{t} - N^{t_{i}}_{t_{j}} & =  \sigma^{-1}(t,X_{t})  (D_{t_i}X_{t} - D_{t_i}X_{t_j})   \sigma(t_i,X_{t_{i}}) \\
& \quad   + (\sigma^{-1}(t,X_{t}) - \sigma^{-1}(t_j ,X_{t_j}) ) D_{t_i}X_{t_j} \sigma(t_i,X_{t_{i}}),
\end{align*}
it follows from the boundedness and Lipschitz continuity  of $\sigma$ and $\sigma^{-1}$ (Lemma \ref{lem:1:inv sigma:lip})
that for any $j>i$ and $t \in[t_{j},t_{j+1}]$,
\begin{align*}
\E_i[|H^{t_{i}}_{t_{j}} - H^i_j|^2] 
& = \frac{\sum_{k=i}^{j-1}\int_{t_{k}}^{t_{k+1}}\E_i[|N^{t_{i}}_{t} - N^{t_i}_{t_{k}}|^2] dt }{(t_{j} - t_i)^2}
\le \frac{C \sum_{k=i}^{j-1}\int_{t_{k}}^{t_{k+1}} \E_i[|D_{t_i}X_{t} - D_{t_i}X_{t_k} |^2 + |X_t - X_{t_k}|^2]  dt}{(t_j-t_i)^2} .
\end{align*}
It now follows from Lemma \ref{lem:mal:SDE:3} the usual bound $ \E_i[|X_t - X_{t_j}]|^2] \le C(t-t_j)$
%, which can be easily  shown from \Hs, 
and Lemma \ref{lem:inc growth} that 
\begin{equation}
\label{eq:MW:est}
\E_i[|H^{t_{i}}_{t_{j}} - H^i_j|^2] \le C\max_k \Delta_k (t_j - t_i)^{-1} \le  CN^{-1} (t_j - t_i)^{-1}.
\end{equation}

 {Since $\E_i[H^{t_i}_{t_j} - H^i_j] = 0$, it follows that 
\[\normexp{\E_i[\Phi(X_T)  (H^{t_i}_{t_j} - H^i_j)]} = \normexp{\E_i[\{\Phi(X_T)  { - \E_i[\Phi(X_T)]}\} (H^{t_i}_{t_j} - H^i_j)]}.\]}
The upper bound 
\[\normexp{\E_i[\{\Phi(X_T)  { - \E_i[\Phi(X_T)]}\} (H^{t_i}_{t_j} - H^i_j)]} \le \normexp{ (\E_i[|\{\Phi(X_T) { - \E_i[\Phi(X_T)]}|^2])^{1/2} (\E_i[|H^{t_i}_{t_j} - H^i_j)|^2])^{1/2}  } \] 
follows from the conditional Cauchy-Schwarz inequality  {(Corollary \ref{cor:cs:ineq})}. 
Therefore, \eqref{eq:MW:est} and $ {\normexp{\Phi(X_T) { - \E_i[\Phi(X_T)]}}} \le K^\alpha(\Phi)(T-t)^{\alpha/2}$ (from \Hg) together imply that
\[
\normexp{\E_i[\Phi(X_T)  (H^{t_i}_{t_j} - H^i_j)]} \le \frac{CN^{-1/2}}{(T-t_i)^{(1-\alpha)/2}}
\]
as required.
The upper bound on $\normexp{\E_i[f_j(X_{t_j},Y_{t_{j+1}},Z_{t_j}) (H^{t_i}_{t_j} - H^i_j)]} $ follows from the Cauchy-Schwarz inequality  {(Corollary \ref{cor:cs:ineq})}, i.e.
 {
\[
\normexp{\E_i[f_j(X_{t_j},Y_{t_{j+1}},Z_{t_j}) (H^{t_i}_{t_j} - H^i_j)]} \le
\normexp{ (\E_i[|f_j(X_{t_j},Y_{t_{j+1}},Z_{t_j})|^2])^{1/2} (\E_i[|H^{t_i}_{t_j} - H^i_j)|^2])^{1/2}  };
\]
from here, one applies the estimate \eqref{eq:MW:est} and the fact that,  similarly to \eqref{eq:bd:f},}
$\normexp{f_j(X_{t_j},Y_{t_{j+1}},Z_{t_j})} \le C(T-t_j)^{\gamma -1}$.
\qed

\begin{lemma}
\label{lem:mal:int}
For all $t_{i},t_{j}\in\pi$ such that $t_{i}\le t_{j}$ and $r\in[t_{j},T]$,
\begin{align}
\e_i[f(r,X_{r},Y_{r},Z_{r}) H^{t_{i}}_{r}] &= \e_{i}[f(r,X_{r},Y_{r},Z_{r})H^{t_{i}}_{t_{j}}].
\label{eq:mal:integrand}
\end{align}
Moreover,
\begin{align}
 \e_i[\int_{t_i}^{T} f(r,X_{r},Y_{r},Z_{r}) H^{t_i}_{r} dr] 
%\nonumber\\ &
& = \e_i[ \sum_{j=i+1}^{N-1} f_j ( X_{t_{j}},Y_{t_{j+1}},Z_{t_j}) H^{t_i}_{t_j} \Delta_j]
+\e_i[\int_{t_i}^{t_{i+1}} f(r,X_{r},Y_r,Z_r)H^{t_i}_r dr] \nonumber \\
& + \e_i[\sum_{j=i+1}^{N-1} \int_{t_j}^{t_{j+1}} (f(r,X_{r},Y_r,Z_r) - f_j (X_{t_{j}},Y_{t_{j+1}},Z_{t_j}) ) H^{t_i}_{t_j} dr] .
\label{eq:mal:weights:sum}
\end{align}
\end{lemma}
{\bf Proof.}
%The arguments are the same as for the proof of Theorem \ref{thm:repr},
%except that in \eqref{eq:f:int} the integral with respect to $dt$ is over $v \in [t_{i},t_{j}]$:
First let \HFd \ be in force and recall, as argued in the proof of Theorem \ref{thm:repr}, that the BSDE solved by $(\yeps , \zeps)$ in Definition \ref{def:1:intermediate BSDEs}
satisfies the conditions of \cite[Theorem 4.2]{ma:zhan:02}.
A key element of the proof of that Theorem is to show that, for almost all $v \in[0,r)$,
\begin{align*}
%\frac{1}{t_{j}-t_{i}} 
%
%\int_{t_{i}}^{t_{j}} 
D_v \feps(\Theta_r) \sigma^{-1}(v,X_v)  \nabX v 
%dv  \nonumber \\
%&
= \nabla_x \feps( & \Theta_r)  \nabX{r} + \nabla_y\feps(\Theta_r) (u(r,X_r) \nabX{r} + \nabla \yeps_r)  \nonumber \\ 
&\qquad + \nabla_z \feps(\Theta_r)(U(r,X_r)\nabX{r} + \nabla \zeps_r ) \qquad \qquad m\times\P-a.e;
\end{align*}
where $U(r,x)$ is defined in \eqref{def:deriv du s}; see the equality just above equation (4.19) in \cite{ma:zhan:02}.
Integrating with respect to $v$ over $v \in [t_i,t_j)$, on the one hand, and between  $v \in [t_i,r)$, on the other,
which yields
\[
\frac1{t_j - t_i} \int_{t_i}^{t_j} D_v \feps(\Theta_r) \sigma^{-1}(v,X_v)  \nabX v dv  = \frac1{r-t_i} \int_{t_i}^r D_v \feps(\Theta_r) \sigma^{-1}(v,X_v)  \nabX v dv.
\]
One then follows the proof of \cite[Theorem 4.2]{ma:zhan:02}, which essentially uses integration-by-parts for Malliavin calculus -- Lemma \ref{lem:mal:ibp} -- to show that
$\e_i[\feps(r,X_{r},\Yeps_{r},\Zeps_{r}) H^{t_{i}}_{r}] = \e_{i}[\feps(r,X_{r},\Yeps_{r}, \Zeps_{r})H^{t_{i}}_{t_{j}}]$.
One extends to the general case \eqref{eq:mal:integrand} by convergence arguments as in the proof of Theorem \ref{thm:repr}
The relation \eqref{eq:mal:weights:sum} is now straightforward to obtain from \eqref{eq:mal:integrand}.
\qed

\begin{lemma}
\label{lem:disc:extra:bounds}
There is a finite constant  $C$ depending only on the bound on $b$ and its derivatives, the bound on $\sigma$ and its derivatives, $L_f$, $\thetaL$, $C_f$, $\thetaC$, and $T$
such that, for all   $i\in\{0,\ldots,N-1\}$,
\begin{align}
&\normexp{\E_{t_i}[\int_{t_i}^{t_{i+1}} f(r,X_{r},Y_r,Z_r)H^{t_i}_r dr] }  \le 
C \int_{t_i}^{t_{i+1}} \frac{dr}{(T-r)^{1- \gamma} \sqrt{r-t_{i}}},
\label{eq:mal:lin} \\
&\normexp{\sum_{j=i+1}^{N-1} \E_{t_i}[(f_j(Y_{t_{j+1}},Z_{t_j}) - f_j(\bar Y_{j+1}, \bar Z_j)) H^{t_i}_{t_j} ]\Delta_{j} }
\nonumber \\
&\qquad \qquad 
\le C \sum_{j=i+1}^{N-1}\frac{ \normexp{Y_{t_{j+1}} - \bar Y_{j+1}} + \normexp{Z_{t_j} - \bar Z_j} }{(T - t_j)^{(1-\thetaL)/2}
\sqrt{t_j-t_i}} \Delta_j,
\label{eq:disc:Z:2} \\
&\normexp{\sum_{j=i+1}^{N-1} \E_{t_i}[ \int_{t_j}^{t_{j+1}} (f(r,X_{r},Y_r,Z_r) - f(t_j,X_{t_j},Y_{t_{j+1}},Z_{t_j}) ) H^{t_i}_{t_j} dr ] }\nonumber \\
& \qquad \qquad \le  { CN^{-1/2} \over (T-t_i)^{(1-\theta_X)/2}} + C\sum_{j=i+1}^{N-1} \frac{\int_{t_j}^{t_{j+1}} \{  \normexp{ Y_r - Y_{t_{j+1}} }
 + \normexp{ Z_r-Z_{t_j} }  \} dr}{(T-t_{j})^{(1-\thetaL)/2}\sqrt{t_j-t_i}} ,
\label{eq:disc:Z:3}
\end{align}
 {where $\gamma := (\thetaC\wedge\frac\alpha2 + \frac\thetaL2)\wedge\thetaC $.}
\end{lemma}
{\bf Proof.}
In what follows, $C$ may change from line to line.
 {
Using the conditional Cauchy-Schwarz inequality (Corollary \ref{cor:cs:ineq}), 
\[
\E_{i}[\int_{t_i}^{t_{i+1}} f(r,X_{r},Y_r,Z_r)H^{t_i}_r ] dr \le \sqrt{C_M}\int_{t_i}^{t_{i+1}} {( \E_{i}[|f(r,X_{r},Y_r,Z_r)|^2])^{1/2} \over \sqrt{r-t_i} } dr ;
\]
}
then, Minkowski's inequality and the moment bound \eqref{eq:bd:f} of Corollary \ref{cor:mom bd} imply that
\begin{align*}
& \normexp{\E_{t_i}[\int_{t_i}^{t_{i+1}} f(r,X_{r},Y_r,Z_r)H^{t_i}_r  dr]} \le \sqrt{C_M}\int_{t_i}^{t_{i+1}} {\normexp{f(r,X_{r},Y_r,Z_r)} \over \sqrt{r-t_i} } dr
 \le C \int_{t_i}^{t_{i+1}} \frac{ dr}{(T-r)^{1- \gamma} \sqrt{r-t_{i}}}.
\end{align*}
Using the Lipschitz continuity of $f$, Minkowski's inequality, and Lemma \ref{lem:mal:weight}, 
\begin{align*}
\normexp{ & \sum_{j=i+1}^{N-1} \E_{t_i}[ \big(f_j (Y_{t_{j+1}},Z_{t_j}) - f_j(\bar Y_{j+1},\bar Z_j) \big) H^{t_i}_{t_j} ] \Delta_{j-1} } 
\le C \sum_{j=i+1}^{N-1}\frac{\normexp{Y_{t_{j+1}} - \bar Y_{j+1} } + \normexp{ Z_{t_j} - \bar Z_j } }{(T - t_j)^{(1-\thetaL)/2}\sqrt{t_j-t_i}} \Delta_j
\end{align*}
For \eqref{eq:disc:Z:3}, the $t$-H\"older continuity of $f$ in \HFt, the Cauchy-Schwarz inequality  {(Corollary \ref{cor:cs:ineq})}, Minkowski's inequality, and H\"older's inequality are needed:
\begin{align*}
&\normexp{\sum_{j=i+1}^{N-1} \E_{t_i}[ \int_{t_j}^{t_{j+1}} \big(f(r,X_{r},Y_r,Z_r) - f_j (X_{t_{j}}, Y_{t_{j+1}},Z_{t_j}) \big) H^{t_i}_{t_j} dr ]} \nonumber \\
& \le C\sum_{j=i+1}^{N-1} \frac{\int_{t_j}^{t_{j+1}} \normexp{ f(r,X_{r},Y_r,Z_r) - f_{j}(X_{r},Y_r,Z_r) } dr }{\sqrt{t_j-t_i}} 
%\nonumber \\
%&\qquad
 + C \sum_{j=i+1}^{N-1} \frac{\int_{t_j}^{t_{j+1}} \normexp{ f_{j}(X_{r},Y_r,Z_r) - f_j(X_{t_{j}},Y_{t_{j+1}},Z_{t_j}) } dr}{\sqrt{t_j-t_i}} \nonumber \\
&\le C \sum_{j=i+1}^{N-1} \frac{\int_{t_j}^{t_{j+1}} \sqrt{r-t_j}  dr}{
% { (T-t_j)^{(1-\thetaL)/2} }
\sqrt{t_j-t_i}} 
%\nonumber \\&\quad
+ C\sum_{j=i+1}^{N-1} \frac{\int_{t_j}^{t_{j+1}} \normexp{X_{r} - X_{t_{j}}} dr } {(T-t_{j})^{1-\theta_X/2}\sqrt{t_j-t_i}}
+ C\sum_{j=i+1}^{N-1} \frac{\int_{t_j}^{t_{j+1}} \{\normexp{ Y_r - Y_{t_{j+1}} }
 + \normexp{ Z_r-Z_{t_j} }  \} dr}{(T-t_{j})^{(1-\thetaL)/2}\sqrt{t_j-t_i}}  
\nonumber 
\end{align*}
The usual upper bound 
$\normexp{X_{r} - X_{t_{j}}} \le C \sqrt{r-t_j}$ implies that 
\[
\int_{t_j}^{t_{j+1}} \normexp{X_{r} - X_{t_{j}}} dr \le C \int_{t_j}^{t_{j+1}}  \sqrt{r-t_j} dr. % \le C (\Delta_j)^{3/2}. 
\]
Now, we obtain the upper bound $ \int_{t_j}^{t_{j+1}} \sqrt{r-t_j}  dr = \frac23 \Delta_j^{3/2} \le C N^{-1/2} \Delta_j$  from Lemma \ref{lem:inc growth}, and substitute it to
the already acquired estimates to obtain
\begin{align*}
&\normexp{\sum_{j=i+1}^{N-1} \E_{t_i}[ \int_{t_j}^{t_{j+1}} \big(f(r,X_{r},Y_r,Z_r) - f_j (X_{t_{j}}, Y_{t_{j+1}},Z_{t_j}) \big) H^{t_i}_{t_j} dr ]} \nonumber \\
&\le C N^{-1} \sum_{j=i+1}^{N-1} \frac{\Delta_j }{
\sqrt{t_j-t_i}} 
%\nonumber \\&\quad
+ C N^{-1}\sum_{j=i+1}^{N-1} \frac{ \Delta_j } {(T-t_{j})^{1-\theta_X/2}\sqrt{t_j-t_i}}
+ C\sum_{j=i+1}^{N-1} \frac{\int_{t_j}^{t_{j+1}} \{\normexp{ Y_r - Y_{t_{j+1}} }
 + \normexp{ Z_r-Z_{t_j} }  \} dr}{(T-t_{j})^{(1-\thetaL)/2}\sqrt{t_j-t_i}}  
\nonumber 
\end{align*}

Applying Lemma \ref{lem:integration:1} to bound the  sums without the integrals is then sufficient to complete the proof. \qed

\vspace{0.3cm}
In the following proposition, we obtain a bound for the error terms on the right hand side of \eqref{eq:disc:Z:3};
these error terms are intrinsically related to the discritization error of the Malliavin weights scheme.
% that is intrinsically related to approximation error caused by these terms.
Proposition \ref{prop:veps:bd} will be essential in the proof of this result.

\begin{proposition}
\label{prop:Z sums:num er}
 {Recall the definition $\gamma := (\thetaC\wedge\frac\alpha2 + \frac\thetaL2)\wedge\thetaC $.}
Let either \Hgexp \ or \ \Hgh \ be in force and
suppose that $0<\beta < (2\gamma)\wedge\alpha\wedge\thetaL$. 
For $\delta, K >0$, define
%\[
$\cC(\delta,K) := K(N^{-1/2}\1_{[1,3]}(\delta)+  N^{-\gamma} \1_{(0,1)}(\delta)),$ 
%\]
and, for $j\in\{0,\ldots,N-1\}$, 
\[\Psi_j := \int_{t_j}^{t_{j+1}} \{  \normexp{ Y_r - Y_{t_{j+1}} }  + \normexp{ Z_r-Z_{t_j} } \} dr .\]
There is a constant $C$ depending only on $L_f$, $\thetaL$, $C_f$, $\thetaC$, $\beta$, $\elip$, the bound on $b$ and its derivatives, the bound on $\sigma$ and it's derivatives, and $T$, but not on $N$, such that, for all $N \ge 1$,
%\begin{align*}
%\sum_{j=0}^{N-1} \frac{\int_{t_j}^{t_{j+1}} \Psi_j }{(T-t_{j})^{(1-\thetaL)/2} } 
%& + \sum_{j=i+1}^{N-1} \frac{\int_{t_j}^{t_{j+1}} \{  \normexp{ Y_r - Y_{t_{j+1}} }
% + \normexp{ Z_r-Z_{t_j} }  \} dr}{(T-t_{j})^{(1-\thetaL)/2} \sqrt{t_j-t_i}} 
%% \nonumber \\
%%\nonumber \\
%% \le 
%% & .
%%\label{eq:1:Z sums:num er}
%\end{align*}
%is bounded above by 
\begin{align*}
\sum_{j=0}^{N-1}  \frac{ \Psi_j }{(T-t_{j})^{(1-\thetaL)/2} }&  \le C (T-t_i)^{(1 + \thetaL - \beta)/2}\cC(\beta +2\gamma,\ln(N)^{1/4} \vee 1) ,\\
\sum_{j=i+1}^{N-1} \frac{  \Psi_j   }{(T-t_{j})^{(1-\thetaL)/2} \sqrt{t_j-t_i}} & \le C(T-t_i)^{-(1+\beta - \thetaL)/2} \cC(\beta +2\gamma,\ln(N)^{1/4} \vee 1)
%C(T-t_i)^{-(1+\beta - \thetaL)/2}(N^{-1/2}\1_{[1,3]}()() +  N^{-\gamma} \1_{(0,1)}(\beta+\thetaP +2\gamma))\] 
\end{align*}
in the case of \Hgexp, and  
\begin{align*}
\sum_{j=0}^{N-1}  \frac{\Psi_j }{(T-t_{j})^{(1-\thetaL)/2} } & \le C (T-t_i)^{(1 + \thetaL - \beta)/2}\cC(\beta +\thetaP+2\gamma,K_\Phi) ,\\
\sum_{j=i+1}^{N-1} \frac{  \Psi_j   }{(T-t_{j})^{(1-\thetaL)/2} \sqrt{t_j-t_i}} & \le C(T-t_i)^{-(1+\beta - \thetaL)/2} \cC(\beta +\thetaP+2\gamma,K_\Phi)
%C(T-t_i)^{-(1+\beta - \thetaL)/2}(N^{-1/2}\1_{[1,3]}()() +  N^{-\gamma} \1_{(0,1)}(\beta+\thetaP +2\gamma))\] 
\end{align*}
%\[C(T-t_i)^{-(1+\beta - \thetaL)/2} (K_\Phi N^{-1/2}\1_{[1,3]}(\beta + \thetaP+2\gamma) +  N^{-\gamma} \1_{(0,1)}(\beta+\thetaP +2\gamma))\]
 in the case of \Hgh.
\end{proposition}

{\bf Proof.}
We will prove the bounds for 
\[
\sum_{j=i+1}^{N-1} \frac{\Psi_j}{(T-t_{j})^{(1-\thetaL)/2} \sqrt{t_j-t_i}} .
\]
The bounds for the $\sum_{j=0}^{N-1}  \frac{ \Psi_j }{(T-t_{j})^{(1-\thetaL)/2} } $
are obtained analogously.
Moreover, we will only prove the result for the terms in $Z$. The bound for the terms in $Y$ are also obtained analogously.
In what follows, $C$  may change from line to line.
%\vspace{0.3cm}
We first prove the result under \HFd \ and \Hge, and then obtain the general result by means of mollification.
Fix $\varepsilon \le \Delta_{N-1}$ and recall the BSDE $(\Yeps,\Zeps)$ from Definition \ref{def:1:intermediate BSDEs} in Section \ref{section:mal:BSDE}.
We use the version of $\Zeps$ provided by Theorem \ref{thm:repr}.
%The triangle inequality yields 
First, apply the triangle inequality to the integrand in order to obtain
$
\normexp{Z_{t} - Z_{t_i}} \le \normexp{Z_{t} - \Zeps_{t}} + \normexp{Z_{t_i} - \Zeps_{t_i}} + \normexp{\Zeps_{t} - \Zeps_{t_i}}  .
$
%This implies that
%\[
%\sum_{j=i+1}^{N-1} \frac{\int_{t_j}^{t_{j+1}}  \normexp{ Z_r-Z_{t_j} }  dr}{(T-t_{j})^{(1-\thetaL)/2}\sqrt{t_j-t_i}} \le
%\sum_{j=i+1}^{N-1} \frac{\int_{t_j}^{t_{j+1}}  \{\normexp{Z_{t} - \Zeps_{t}} + \normexp{Z_{t_i} - \Zeps_{t_i}}
% + \normexp{\Zeps_{t} - \Zeps_{t_i}} \}  dr}{(T-t_{j})^{(1-\thetaL)/2}\sqrt{t_j-t_i}} .% \cA_{i,1} + \cA_{i,2} + \cA_{i,3}.
%\]
% We bound $\sum_{i=0}^{N-2} \cA_i^2 \Delta_i$; $\cA_{N-1} = 0$. 

%  { 
% Starting with the first term, H\"older's inequality yields
% \begin{align*}
% &\sum_{i=0}^{N-2}\bigg(\sum_{j=i+1}^{N-1} \frac{\int_{t_j}^{t_{j+1}}  \normexp{Z_{t} - \Zeps_{t}}   dr}
% {(T-t_{j})^{(1-\thetaL)/2}\sqrt{t_j-t_i}} \bigg)^2 \Delta_i \\
% & \le \sum_{i=0}^{N-1}\bigg(\sum_{j=i+1}^{N-1} \frac{(\int_{t_j}^{t_{j+1}}  \normexp{Z_{t} - \Zeps_{t}}^2 dr)^{1/2} \Delta_j^{1/2}}
% {(T-t_{j})^{(1-\thetaL)/2}} \bigg)^2 \\
% & \le \max_{0\le j \le N-1}\frac{\Delta_j}{ (T-t_{j})^{1-\thetaL}}N^2  \int_{0}^{T} \normexp{Z_{t} - \Zeps_{t}}^2 dr \le CN\varepsilon^{2\gamma}
% \end{align*}
% Taking $\varepsilon = CN^{-1/\gamma}$ is sufficient to ensure that the above sum is bounded by $CN^{-1}$.
% Without loss of generality, we can assume that $\varepsilon < \Delta_{N-1} = CN^{-1/\beta}$.
% }
To bound the terms in $Z -\Zeps$, recall the bound \eqref{eq:cut:error:2} from  Corollary \ref{cor:cut:error}. 
%gives  
%\[\normexp{Z_{t} - \Zeps_{t}} \le C\int_{t\vee(T - \varepsilon)}^T \frac{dr}{(T-r)^{1-\gamma}\sqrt{r-t}} \qquad \forall \ t\in[0,T). \]
%  Assume that $\varepsilon \le \Delta_{N-1} = TN^{-1\beta}$.
For $j\le N-2$, the bound on $\normexp{Z_{t} - \Zeps_{t}}$ implies that 
\begin{align*}
\int_{t_j}^{t_{j+1}} \normexp{Z_{t} - \Zeps_{t}} dt & 
\le C \int_{t_j}^{t_{j+1}} \frac{\int_{T - \varepsilon}^T (T-r)^{\gamma-1}dr }{\sqrt{t_{N-1}-t} } dt 
%\le C\varepsilon^\gamma \Delta_j 
\le C\varepsilon^\gamma \int_{t_j}^{t_{j+1}} \frac{dt }{\sqrt{t_{N-1}-t} } .
\end{align*}
Lemma \ref{lem:simple:int} yields $\int_{t_j}^{t_{j+1}} (t_{N-1}-t)^{-1/2} dt  \le 2 \Delta_j(t_{N-1} - t_j)^{-1/2}$.
%\begin{align*}
%\int_{t_j}^{t_{j+1}} \frac{dt }{\sqrt{t_{N-1}-t} } & = 2 \big\{ (t_{N-1}-t_j)^{1/2} - (t_{N-1} - t_{j+1})^{1/2}\big\} \\
%& \le 2\Big\{ \frac{t_{N-1}-t_j}{(t_{N-1}-t_j)^{1/2}} - \frac{t_{N-1} - t_{j+1}}{(t_{N-1}-t_j)^{1/2}} \Big\}
% = \frac{2\Delta_j}{(t_{N-1}-t_j)^{1/2}}
%\end{align*}
Therefore, applying Lemma \ref{lem:integration:1} implies that
\begin{align}
\sum_{j=i+1}^{N-2} \frac{  \int_{t_j}^{t_{j+1}} \normexp{Z_{t} - \Zeps_{t}} dt }{(T-t_{j})^{(1-\thetaL)/2}\sqrt{t_j-t_i}}
& \le C \varepsilon^\gamma \sum_{j=i+1}^{N-2} \frac{  \Delta_j }{(t_{N-1}-t_{j})^{1-\thetaL/2}\sqrt{t_j-t_i}} 
\le \frac{C\varepsilon^\gamma}{(t_{N-1} - t_i)^{(1-\thetaL)/2} }.
\label{eq:num er:1}
\end{align}
{Then, use  $(t_{N-1} - t_i)^{-1/2} \le 2 (T - t_i)^{-1/2}$ on the denominator on the right hand side. }
For the outstanding term, $j=N-1$, we implement  Lemma \ref{lem:integration:1} to show that
\begin{align*}
\int_{t_{N-1}}^{T} \normexp{Z_{t} - \Zeps_{t}} dt & 
\le C \int_{t_{N-1}}^{T} \Big\{ \int_{t}^T (T-r)^{\gamma-1}(r-t)^{-1/2}dr \Big\} dt 
\le C \Delta_{N-1}^{1/2+\gamma},
\end{align*}
whence it follows that
\begin{equation}
\frac{  \int_{t_{N-1}}^{T} \normexp{Z_{t} - \Zeps_{t}} dt }{\Delta_{N-1}^{(1-\thetaL)/2}\sqrt{t_{N-1}-t_i}}
\le \frac{C \Delta_{N-1}^{\gamma + \thetaL/2} }{\sqrt{t_{N-1}-t_i}} 
{ \le \frac{C \Delta_{N-1}^{\gamma + \thetaL/2} }{\sqrt{T-t_i}} }
\label{eq:num er:2} 
\end{equation}
Combining \eqref{eq:num er:1} and \eqref{eq:num er:2}, it follows that
\begin{align}
\sum_{j=i+1}^{N-1} \frac{\int_{t_j}^{t_{j+1}}  \normexp{Z_{t} - \Zeps_{t}}   dr}
{(T-t_{j})^{(1-\thetaL)/2}\sqrt{t_j-t_i}} 
\le {\frac{C\varepsilon^\gamma}{(T - t_i)^{(1-\thetaL)/2} } + \frac{C N^{-1} }{\sqrt{T-t_i}} }
%\nonumber \\
%& \le 
%C\varepsilon^{2\gamma}\sum_{i=0}^{N-2} \frac{\Delta_i }{(t_{N-1} - t_i)^{1-\thetaL} } 
%+ C\Delta_{N-1}^{2\gamma + \thetaL}\sum_{i=0}^{N-2} \frac{\Delta_i }{t_{N-1}-t_i} 
% \nonumber \\ & 
%\le C\varepsilon^{2\gamma} + CN^{-2}\big(1 + \ln(N) \big)
\label{eq:num er:3}
\end{align}
where we have used that $ \Delta_{N-1}^{2\gamma +\thetaL} = TN^{(2\gamma +\thetaL)/\beta}$ 
and $\beta < (2\gamma)\wedge\thetaL$. Analogously, we can also show that
\begin{align}
&\sum_{j=i+1}^{N-1} \frac{  \normexp{Z_{t_j} - \Zeps_{t_j}}   \Delta_j }
{(T-t_{j})^{(1-\thetaL)/2}\sqrt{t_j-t_i}} 
\le  \frac{C\varepsilon^\gamma}{(T - t_i)^{(1-\thetaL)/2} } + \frac{C N^{-1} }{\sqrt{T-t_i}} 
%\le C\varepsilon^{2\gamma} + CN^{-2}\big(1 + \ln(N) \big).
\label{eq:num er:4}
% & \le C\sum_{i=0}^{N-2}\bigg(C\int_{T-\varepsilon}^T \frac{dr}{(T-r)^{1-\gamma}}\sum_{j=i+1}^{N-2} \frac{\Delta_j}
% {(t_{N-1}-t_{j})^{1-\thetaL/2}\sqrt{t_j-t_i}} \bigg)^2 \Delta_i \\
% & \quad + C\sum_{i=0}^{N-2}\bigg( \frac{C\int_{t_{N-1}}^T \frac{dr}{(T-r)^{1-\gamma}\sqrt{r-t_{N-1} } }\Delta_{N-1}}
% {\Delta_{N-1}^{(1-\thetaL)/2}\sqrt{t_{N-1}-t_i}} \bigg)^2 \Delta_i \\
% & \le C\varepsilon^{2\gamma} \sum_{i=0}^{N-2} \frac{\Delta_i}
% {(t_{N-1}-t_{i})^{1-\thetaL}}  + C\sum_{i=0}^{N-2}\bigg( \frac{C\Delta_{N-1}^{\gamma + \thetaL/2} }
% {\sqrt{t_{N-1}-t_i}} \bigg)^2 \Delta_i 
\end{align}
% Taking $\gamma \le \beta \le \gamma + \thetaL/2$ implies that the above is bounded by $CN^{-1}$.
Recalling the BSDEs $(y,z)$ and $(\yeps,\zeps)$ from Definition \ref{def:1:intermediate BSDEs} and that $\Zeps = z + \zeps$, the triangle inequality yields 
$
\normexp{\Zeps_{t} - \Zeps_{t_i}} \le \normexp{z_t - z_{t_i}} + \normexp{\zeps_t - \zeps_{t_i}}.
$
In the proof of \cite[Theorem 1.1]{gobe:makh:10}, in bounding the terms $E_1$ and $E_2$, it is shown that, for all $t\in[0,T]$, 
\begin{align*}
\normexp{z_t - z_{t_i}}^2 \le \frac{C(t-t_i)}{(T-t)^{1-\alpha}} + C\int_{t_i}^t \normexp{\nabla_x^2 u(r,X_r)}^2 dr .
% \le \frac{C(t-t_i) + C}{(T-t)^{1-\alpha}}.
\end{align*}
Lemma \ref{lem:lin:pde:bds} implies $\int_{t_i}^t \normexp{\nabla_x^2 u(r,X_r)}^2 dr \le C \int_{t_i}^t (T-r)^{\alpha-2}dr $.
Now, applying Jensen's inequality, Lemma \ref{lem:simple:int}, Lemma \ref{lem:integration:1}, and the above bound, one obtains
\begin{align}
& 
\sum_{j=i+1}^{N-1} \frac{\int_{t_j}^{t_{j+1}}  \normexp{z_{r} - z_{t_j}}   dr}
{(T-t_{j})^{(1-\thetaL)/2}\sqrt{t_j-t_i}} 
%\nonumber \\
%& 
\le \sum_{j=i+1}^{N-1}  \frac{C\Delta_j\int_{t_j}^{t_{j+1}} (T - r)^{(\alpha-1)/2} dr 
+ C\int_{t_j}^{t_{j+1}}\Big(\int_{t_j}^r (T-t)^{\alpha-2}dt\Big)^{1/2} dr }
{(T-t_{j})^{(1-\thetaL)/2}\sqrt{t_j-t_i}} 
\nonumber \\
& \qquad \qquad \qquad \le \sum_{j=i+1}^{N-1}  \frac{C\Delta_j^{2}(T-t_j)^{(\alpha-1)/2}
+ C\Delta_j^{1/2}\Big(\int_{t_j}^{t_{j+1}}(t_{j+1}-t)(T-t)^{\alpha-2} dt\Big)^{1/2} }
{(T-t_{j})^{(1-\thetaL)/2}\sqrt{t_j-t_i}} .
\label{eq:num er:5}
\end{align}
For $j\le N-2$, one can apply Lemma \ref{lem:inc growth} and Lemma \ref{lem:simple:int} to show that
\[
\int_{t_j}^{t_{j+1}}(t_{j+1}-t)(T-t)^{\alpha-2} dt \le \frac{\Delta_j}{(T-t_{j+1})^{1-\beta}}\int_{t_j}^{t_{j+1}}(T-t)^{\alpha-\beta -1} dt
\le \frac{CN^{-1}\Delta_j}{ (T-t_j)^{1- \alpha + \beta}}.
\]
%where we have used Lemma \ref{lem:inc growth} to show that $(\Delta_j /\Delta_{j+1}) \Delta_{j+1} (T-t_{j+1})^{\beta -1 } \le CN^{-1}$ and the direct computation
%\begin{align*}
%\int_{t_j}^{t_{j+1}} \frac{dt }{(T - t)^{1 - \alpha + \beta} } & = \frac 1{\alpha-\beta} \big\{ (T-t_j)^{\alpha-\beta} 
%- (T - t_{j+1})^{\alpha-\beta}\big\} \\
%& \le\frac 1{\alpha-\beta} \Big\{ \frac{T-t_j}{(T-t_j)^{1-\alpha+\beta}} - \frac{T - t_{j+1}}{(T - t_j)^{1-\alpha+\beta}} \Big\}
% = \frac{\Delta_j}{(\alpha-\beta)(T-t_j)^{1-\alpha+\beta}}
%\end{align*}

On the other hand, for $j=N-1$, since $\beta < \alpha$,
\[
\int_{t_{N-1}}^{T}(T-t)(T-t)^{\alpha-2} dt = \frac1\alpha\Delta_{N-1}^\alpha = \frac T\alpha N^{-\alpha/\beta} \le  \frac T\alpha N^{-1}.
\]
Substituting these bounds into \eqref{eq:num er:5} and implementing Lemma \ref{lem:inc growth} and Lemma \ref{lem:integration:1}, we obtain
\begin{align}
%&
\sum_{j=i+1}^{N-1} \frac{\int_{t_j}^{t_{j+1}}  \normexp{z_{r} - z_{t_j}}   dr}
{(T-t_{j})^{(1-\thetaL)/2}\sqrt{t_j-t_i}} 
%\nonumber\\
& \le
% C N^{-1}
% \\ & \quad 
 C N^{-1/2}
% N^{-1} 
\sum_{j=i+1}^{N-1}  \frac{
\Delta_j }
{(T-t_{j})^{1+(\beta-\alpha-\thetaL)/2}\sqrt{t_j-t_i}} 
% \\ &\qquad 
+CN^{-1/2} 
\frac{\Delta_{N-1}^{\thetaL/2}  }
{\sqrt{t_{N-1}-t_i}} 
\nonumber\\
& \le \frac{ C N^{-1/2} }{(T-t_{i})^{(1+\beta-\alpha-\thetaL)/2}} + \frac{CN^{-1}}{\sqrt{T-t_i}}.
%& \le CN^{-1} + CN^{-1} \sum_{i=0}^{N-2} \frac{\Delta_i}{(t_{N-1}-t_i)^{1-\thetaL}}  \le CN^{-1}.
\label{eq:num er:6}  
\end{align}
In the bounds \eqref{eq:zeps:bd:1}, we used the inequality
\begin{equation}
\normexp{\zeps_{r} - \zeps_{t_i}} \le C \int_{t_i}^{r} \normexp{\aeps t}dt + C \int_{t_i}^{r} \normexp{\Veps_t}dt + C \Delta_i^{1/2}.
\label{eq:num er:7} 
\end{equation}
% The sum $ \sum_{j=i+1}^{N-1} \frac{\int_{t_j}^{t_{j+1}}  \normexp{\zeps_{r} - \zeps_{t_j}}   dr}
% {(T-t_{j})^{(1-\thetaL)/2}\sqrt{t_j-t_i}}$ is therefore, using Lemma \ref{prop:veps:bd}, bounded by }
Using $\normexp{\aeps t} \le C (T-t)^{(\alpha + \thetaL -3)/2}$ as shown Lemma \ref{lem:mal:ueps veps}, 
Lemma \ref{lem:inc growth}, Lemma \ref{lem:simple:int}, and Lemma \ref{lem:integration:1},
it follows that
\begin{align}
& \sum_{j=i+1}^{N-1} \frac{\int_{t_j}^{t_{j+1}} \{ \int_{t_j}^r \normexp{\aeps{t}} dt \}  dr}{(T-t_{j})^{(1-\thetaL)/2}\sqrt{t_j-t_i}} 
% \\
\le C \sum_{j=i+1}^{N-1} \frac{\int_{t_j}^{t_{j+1}} \{\int_{t_j}^r(T-t)^{(\thetaL+\alpha-3)/2}dt \}dr}
{(T-t_{j})^{(1-\thetaL)/2}\sqrt{t_j-t_i}} 
\nonumber \\
& \le C \sum_{j=i+1}^{N-1}  \frac{ \int_{t_j}^{t_{j+1}} (T-r)^{(\alpha-2)/2} dr  \int_{t_i}^{t_{j+1}} (T-t)^{(\thetaL-1)/2}dt }
{(T-t_{j})^{(1-\thetaL)/2}\sqrt{t_j-t_i}}
%+ C \frac{\int_{t_{N-1} }^{T}   (T-r)^{(\thetaL + \alpha - 1)/2}dr   } 
%{(T-t_{N-1})^{(1-\thetaL)/2}\sqrt{t_{N-1}-t_i}} 
%\\
%&
\le C \sum_{j=i+1}^{N-1}  \frac{ \Delta_j^{(3 + \thetaL)/2} }
{(T-t_{j})^{(3-\alpha-\thetaL)/2}\sqrt{t_j-t_i}}
\nonumber \\
%\end{align*}
%where we have used
%\[ \int_{t_{N-1}}^T \{ \int_{t_{N-1}}^r (T-t)^{(\alpha + \thetaL -3)/2} dt \} dr =  
%\frac{\int_{t_{N-1}}^T (T-t)^{(\alpha + \thetaL -1)/2} dt }{ (1- \alpha - \thetaL)/2 }.
%\]
%As before, we use direct computation to show
%\begin{align*}
%\int_{t_j}^{t_{j+1}} \frac{dt }{(T - t)^{1 - \alpha/2} } & = \frac 2\alpha \big\{ (T-t_j)^{\alpha/2} - (T - t_{j+1})^{\alpha/2}\big\} \\
%& \le\frac 2\alpha \Big\{ \frac{T-t_j}{(T-t_j)^{1-\alpha/2}} - \frac{T - t_{j+1}}{(T - t_j)^{1-\alpha/2}} \Big\}
% = \frac{2\Delta_j}{\alpha(T-t_j)^{1-\alpha/2}}
%\end{align*}
%whence it follows that
%\begin{align}
%\sum_{j=i+1}^{N-1} \frac{\int_{t_j}^{t_{j+1}} \{ \int_{t_j}^r \normexp{\aeps{t}} dt \}  dr}{(T-t_{j})^{(1-\thetaL)/2}\sqrt{t_j-t_i}}
% &\le C \sum_{j=i+1}^{N-2}  \frac{\Delta_j  \Delta_j^{(1+\thetaL)/2}  }
% {(T-t_{j+1})^{(3-\thetaL-\alpha)/2}\sqrt{t_j-t_i}}
% + C \frac{ \Delta_{N-1}^{(1+\alpha + \thetaL)/2} } 
% {(T-t_{N-1})^{(1-\thetaL)/2}\sqrt{t_{N-1}-t_i}} \\
& \le C\Big( \max_{0\le i \le N-1 }\frac{ \Delta_i^{1+\thetaL}}{(T-t_i)^{1-\beta} }\Big)^{1/2} \sum_{j=i+1}^{N-1}  \frac{\Delta_j   }
{(T-t_j)^{(2+\beta-\thetaL-\alpha)/2}\sqrt{t_j-t_i}} 
%\nonumber \\
%& \qquad + C \frac{ \Delta_{N-1}^{(1+\alpha + \thetaL)/2} } 
%{(T-t_{N-1})^{(1-\thetaL)/2}\sqrt{t_{N-1}-t_i}} 
%\nonumber\\
% &\qquad + C \frac{ \Delta_{N-1}^{(\alpha + 2\thetaL)/2} } 
% {\sqrt{t_{N-1}-t_i}} \\
%& 
\le \frac{CN^{-1/2}}{(T-t_i)^{(1+\beta-\thetaL - \alpha)/2} } 
%+  \frac{ C N^{-(\alpha + 2\thetaL)/(2\beta)} } 
%{\sqrt{t_{N-1}-t_i}} .
\label{eq:num er:8}
\end{align}
On the other hand, we obtain bounds for $\normexp{\Veps_t}$
from Proposition \ref{prop:veps:bd} under \Hgexp \ or \Hgh.
Let us work under \Hgexp.
It follows from Lemma \ref{lem:simple:int}  and Lemma \ref{lem:integration:1} that, for all $j$  and $r\in[t_j,t_{j+1}]$,
\begin{align*}
\int_{t_j}^r \normexp{\Veps_{t}} dt & \le C \| \Phi\|_\infty  \int_{t_j}^r \big\{ \int_t^{T-\varepsilon} (T-s)^{(\beta-1)/2} (T-s)^{(
% \alpha + 
\thetaL - \beta - 2)/2} (s-t)^{-1/2} ds\big\} dt\\
& \le C \| \Phi\|_\infty  \varepsilon^{(\beta-1)/2} \int_{t_j}^r \big\{ \int_t^{T} (T-s)^{\thetaL - \beta  - 1} (s-t)^{-1/2} ds\big\} dt \\
& \le C \| \Phi\|_\infty  \varepsilon^{(\beta-1)/2} \int_{t_j}^r (T-t)^{(\thetaL - \beta -1)/2} dt 
\le C\| \Phi\|_\infty   \frac{ \varepsilon^{(\beta-1)/2} \Delta_j}{(T-t_j)^{(1+\beta-\thetaL)/2} }
\end{align*}
%where we have used the direct computations\begin{align*}
%\int_{t_j}^{r} \frac{dt }{(T - t)^{(1 - \thetaL/4)/2} } & = \frac 2{1+ \thetaL/4} \big\{ (T-t_j)^{(1+ \thetaL/4)/2}
% - (T - r )^{(1+ \thetaL/4)/2}\big\} \\
%& \le \frac 2{1+ \thetaL/4} \Big\{ \frac{T-t_j}{(T-t_j)^{(1+ \thetaL/4)/2}} - \frac{T - r }{(T - t_j)^{(1+ \thetaL/4)/2}} \Big\}
% \le \frac{C\Delta_j}{(T-t_j)^{(1+ \thetaL/4)/2}}
%\end{align*}
Therefore, using Lemma \ref{lem:inc growth}, Lemma \ref{lem:integration:1} and the above bound,
\begin{align}
\sum_{j=i+1}^{N-1} \frac{\int_{t_j}^{t_{j+1}} \{ \int_{t_j}^r \normexp{\Veps_{t}} dt \}  dr}{(T-t_{j})^{(1-\thetaL)/2}\sqrt{t_j-t_i}}
& \le C \| \Phi\|_\infty  \varepsilon^{(\beta-1)/2} \max_i \Delta_i \sum_{j=i+1}^{N-1} \frac{\Delta_j}{(T-t_{j})^{1- ( \thetaL - \beta)/2}\sqrt{t_j-t_i}}
\nonumber \\
 & \le  \frac{ C \| \Phi\|_\infty  \varepsilon^{(\beta-1)/2} N^{-1} }{(T-t_{i})^{(1 +  \beta - \thetaL)/2} }.
\label{eq:num er:9}
\end{align}
Now, substituting \eqref{eq:num er:8} and \eqref{eq:num er:9} into \eqref{eq:num er:7}, it follows that
\begin{align}
\sum_{j=i+1}^{N-1} \frac{\int_{t_j}^{t_{j+1}}  \normexp{\zeps_{r} - \zeps_{t_j}}   dr}
{(T-t_{j})^{(1-\thetaL)/2}\sqrt{t_j-t_i}}
\le
\frac{CN^{-1/2}}{(T-t_i)^{(1+\beta-\thetaL - \alpha)/2} } 
+ \frac{ C \| \Phi\|_\infty  \varepsilon^{(\beta-1)/2} N^{-1} }{(T-t_{i})^{(1 +  \beta - \thetaL)/2} }
%& \le 
%\sum_{i=0}^{N-2} \frac{CN^{-1}\Delta_i}{(T-t_i)^{1+\beta-\thetaL - \alpha} }
%%+\sum_{i=0}^{N-2}\frac{ C N^{-3}\Delta_i } {{t_{N-1}-t_i}} \nonumber \\
%%& \qquad
% +\sum_{i=0}^{N-2}\frac{C \| \Phi\|_\infty^2  \varepsilon^{\beta-1} N^{-2}\Delta_i}{(T-t_{i})^{1 + \beta - \thetaL }} 
%\nonumber \\
%& \le CN^{-1} 
%%+ CN^{-3}\big(1+\ln(N)\big) 
%+ C \| \Phi\|_\infty ^2 \varepsilon^{\beta-1} N^{-2} .
\label{eq:num er:10}
\end{align}
Combining \eqref{eq:num er:3}, \eqref{eq:num er:4}, \eqref{eq:num er:6} and \eqref{eq:num er:10} yields
\begin{align*}
&\sum_{j=i+1}^{N-1} \frac{\int_{t_j}^{t_{j+1}}  \normexp{Z_{r} - Z_{t_j}}   dr} 
{(T-t_{j})^{(1-\thetaL)/2}\sqrt{t_j-t_i}} \\
&\le \frac{CN^{-1/2}}{(T-t_i)^{(1+\beta-\thetaL - \alpha)/2} } 
+\frac{CN^{-1}  }{\sqrt{T-t_i}}
+ \frac{C\varepsilon^\gamma}{(T - t_i)^{(1-\thetaL)/2} }
+ \frac{ C \| \Phi\|_\infty  \varepsilon^{(\beta-1)/2} N^{-1} }{(T-t_{i})^{(1 +  \beta - \thetaL)/2} }
%\le CN^{-1} + C\varepsilon^{2\gamma} + C \| \Phi\|_\infty ^2 \varepsilon^{\beta-1} N^{-2}
\end{align*}
and we take $\varepsilon = N^{-1/(2\gamma)}$ if $1-\beta -2\gamma < 0$ and $\varepsilon = N^{-1}$ otherwise to complete the proof under \HFd \ and \Hge.
The proof under \HFd \ and \Hgh is analogous.

To prove the result without \HFd \ or \Hge, recall the mollified BSDE $(Y_M,Z_M)$ from Corollary \ref{cor:1:moll}.
Set $M =   (3 \ln(N) )^{1/4}$ and $R(M)$ equal to $3 L_f e^{M^{2}/2}  $. Substituting equations \eqref{eq:moll:term} and \eqref{eq:moll:driver} into \eqref{eq:moll:diff},
$\normexp{Z_s - Z_{M,s}} \le CN^{-1} (T-s)^{-1/2} $ for all $s\in[0,T)$,
whence the triangle inequality and Lemma \ref{lem:simple:int} imply
\begin{align*}
\int_{t_j}^{t_{j+1}}  \normexp{Z_{r} - Z_{t_j}}   dr & \le \int_{t_j}^{t_{j+1}}  \normexp{Z_{r} - Z_{M,r} }   dr +   \normexp{Z_{t_j} - Z_{M,t_j}}   \Delta_j + \int_{t_j}^{t_{j+1}}  \normexp{Z_{M,r} - Z_{M,t_j}}   dr \\
& \le C N^{-1} \Delta_j (T-t_j)^{-1/2} + \int_{t_j}^{t_{j+1}}  \normexp{Z_{M,r} - Z_{M,t_j}}   dr .
\end{align*}
The proof is then completed with  use of Lemma \ref{lem:integration:1}.
\qed

\vspace{0.3cm}
We come to the main result of this section, namely the error estimation for the Malliavin weights scheme.
\begin{theorem}
\label{thm:1:disc er}
 {Recall the definition $\gamma := (\thetaC\wedge\frac\alpha2 + \frac\thetaL2)\wedge\thetaC $.}
Let \Hgexp \ or \Hgh \ be and force and suppose that $0<\beta < \gamma \wedge \alpha \wedge \thetaL$.
For $\delta, K >0$, define
$\cC(\delta,K) := KN^{-1/2}\1_{[1,3]}(\delta)+  N^{-\gamma} \1_{(0,1)}(\delta)$.
There is a constant $C$ depending only on $L_f$, $\thetaL$, $C_f$, $\thetaC$, $\beta$, $\elip$, the bound on $b$ and its derivatives, the bound on $\sigma$ and it's derivatives, $K^\alpha(\Phi)$ and $T$, but not on $N$, such that, for all  $ N\ge 1$, 
\begin{align*}
%\cE^D(N) 
%:= \max_{0\le i \le N-1} \normexp{Y_{t_i} - Y_i}^2 + \sum_{i=0}^{N-1} \normexp{ Z_{t_i} - Z_i }^2\Delta_i 
%\le  
\left.
\begin{array}{l}
\normexp{Y_{t_i} - \bY Ni} \le  C (T-t_i)^{(1 + \thetaL - \beta)/2}\cC(\beta +2\gamma,\ln(N)^{1/4} \vee 1)  ,\\
\normexp{Z_{t_i} - \bZ Ni } \le C (T-t_i)^{-(1 + \beta - \thetaL )/2}\cC(\beta +2\gamma,\ln(N)^{1/4} \vee 1) \\
\hspace{3cm}  {+ CN^{-1/2} (T-t_i)^{-(1  -\alpha\wedge(2\gamma) \wedge\theta_X )/2}}.
\end{array}
\right\}
\text{ in the case of \Hgexp,} 
\\ \\
\left.
\begin{array}{l}
\normexp{Y_{t_i} -  \bY Ni} \le  C (T-t_i)^{(1 + \thetaL - \beta)/2}\cC(\beta + \thetaP+2\gamma,K_\Phi)  ,\\
\normexp{Z_{t_i} - \bZ Ni } \le C (T-t_i)^{-(1 + \beta - \thetaL )/2}\cC(\beta +\thetaP + 2\gamma,K_\Phi) \\
\hspace{3cm}  {+ CN^{-1/2} (T-t_i)^{-(1  -\alpha\wedge(2\gamma) \wedge\theta_X )/2}}
\end{array}
\right\} 
\text{ in the case of \Hgh.}
\end{align*}
\end{theorem}
{\bf Proof.}
In what follows, $C$ may change from line to line.
For simplicity, we omit the process $X$ from the driver, so that $f(t,y,z) := f(t,X_t,y,z)$ and $f_j(y,z) := f_j(X_{t_j},y,z)$.
Fix $i\in\{0,\ldots,N-1\}$. 
Using the estimates from Lemma \ref{lem:mw:diffs} and Lemma \ref{lem:disc:extra:bounds},
and \eqref{eq:mal:weights:sum} from Lemma \ref{lem:mal:int}, it follows that
\begin{align}
& \normexp{ Z_{t_i} - \bZ Ni }  =\normexp{\e_i[\Phi(X_T) H^{t_i}_T - \Phi(X_T)H^i_N + \int_{t_i}^T f(t,Y_t,Z_t)H^{t_i}_t dt
-\sum_{j=i+1}^{N-1} f_j(\bY N{j+1}, \bZ Nj)H^i_j \Delta_j  ] } \nonumber \\
&\le \normexp{\e_i[\int_{t_i}^{t_{i+1}} f(r,Y_r,Z_r)H^{t_i}_r dr] } + \normexp{\e_i[\Phi(X_T)(H^{t_i}_T - H^i_N)] } 
+\normexp{\sum_{j=i+1}^{N-1} \E_i[ f_j(Y_{t_{j+1}},Z_{t_j}) (H^{t_i}_{t_j}-H^i_j) ] \Delta_j } \nonumber \\
& + \normexp{\sum_{j=i+1}^{N-1} \e_i[(f_j(Y_{t_{j+1}},Z_{t_j}) - f_j(\bY N{j+1},\bZ Nj)) H^{i}_{j} ] \Delta_j } 
%\nonumber \\
%& \qquad
 + \normexp{\sum_{j=i+1}^{N-1} \E_i[  \int_{t_j}^{t_{j+1}} (f(r,Y_r,Z_r) - f_j(Y_{t_{j+1}},Z_{t_j}) ) H^{t_i}_{t_j} dr ] }
 \nonumber \\
& \le \frac{C  N^{-1/2}}{(T-t_i)^{(1-\alpha\wedge(2\gamma) {\wedge\theta_X})/2}} 
+ C \int_{t_i}^{t_{i+1}} \frac{dr}{(T-r)^{1-\gamma} \sqrt{r-t_i} } + C\cH(i) 
+ C \sum_{j=i+1}^{N-1}\frac{\Theta_j \Delta_j }{(T - t_j)^{(1-\thetaL)/2}\sqrt{t_j-t_i}} 
\label{eq:disc er:2}
\end{align}
where $\Theta_j := \normexp{Y_{t_{j+1}} - \bY N{j+1}} + \normexp{Z_{t_j} - \bZ Nj} $ and 
\[
\cH(i) := \sum_{j=i+1}^{N-1} \frac{\Psi(j) }{(T-t_{j})^{(1-\thetaL)/2}\sqrt{t_j-t_i}}, 
\quad \Psi(j) :=  \int_{t_j}^{t_{j+1}} \{  \normexp{ Y_r - Y_{t_{j+1}} } + \normexp{ Z_r-Z_{t_j} }  \} dr.
\]
In \eqref{eq:disc er:2}, we have estimated $\normexp{\sum_{j=i+1}^{N-1} \E_i[ f_j(Y_{t_{j+1}},Z_{t_j}) (H^{t_i}_{t_j}-H^i_j) ] \Delta_j } $
by $C N^{-1/2}\sum_{j=i+1}^{N-1} (T-t_j)^{\gamma-1}(T_j - t_i)^{-1/2} \Delta_j$ using Lemma \ref{lem:mw:diffs}, and the latter sum by $C (T-t_i)^{\gamma - 1/2}$ using Lemma \ref{lem:integration:1}.
Using a similar technique,  $\normexp{Y_{t_i} - \bY Ni}$ is bounded above by
\begin{align}
& C \int_{t_i}^{t_{i+1}} \frac{dr}{(T-r)^{1-\gamma}} + 
C \sum_{j=i}^{N-1}\frac{ \Psi(j)}{(T - t_j)^{(1-\thetaL)/2}}  
+ C \sum_{j=i}^{N-1}\frac{ \Theta_j \Delta_j  }{(T - t_j)^{(1-\thetaL)/2}} \nonumber\\
& \qquad \le \quad C N^{-1/2}  + C \sum_{j=i}^{N-1}\frac{ \Psi(j)}{(T - t_j)^{(1-\thetaL)/2}}  
+ C \sum_{j=i}^{N-1}\frac{ \Theta_j \Delta_j  }{(T - t_j)^{(1-\thetaL)/2}} 
\label{eq:disc er:3}
\end{align}
where we have used Lemma \ref{lem:inc growth} and Lemma \ref{lem:simple:int}  on the first integral to obtian
$ \int_{t_i}^{t_{i+1}} (T-r)^{\gamma-1} dr \le C \Delta_j(T-t_j)^{\gamma - 1} \le CN^{-1/2}$. 
%where we have used the following direct computation to bound the first integral term:
%\begin{align*}
%\int_{t_i}^{t_{i+1}} \frac{dt }{(T - t)^{1 - \gamma} } & = \frac 1 \gamma \big\{ (T-t_i)^{\gamma} - (T - t_{i+1})^{\gamma}\big\} \\
%& \le\frac 1 \gamma \Big\{ \frac{T-t_i}{(T-t_i)^{1-\gamma}} - \frac{T - t_{i+1}}{(T - t_i)^{1-\gamma}} \Big\}
% = \frac{\Delta_j}{\gamma(T-t_i)^{1-\gamma}} \le CN^{-1/2}.
%\end{align*}
%
%Define $\Theta_i := \normexp{Y_{t_i} - Y_i} + \normexp{Z_{t_i} - Z_i}$.
It follows from \eqref{eq:disc er:2} and \eqref{eq:disc er:3} that
\begin{align*}
\Theta_i & \le \frac{C N^{-1/2}}{(T-t_i)^{(1-\alpha\wedge(2\gamma)  {\wedge\theta_X})/2} }
+ C \int_{t_i}^{t_{i+1}} \frac{dr}{(T-r)^{1-\gamma} \sqrt{r-t_i} } + C \cH(i) 
+ C \sum_{j=i+1}^{N-1}\frac{ \Theta_j \Delta_j }{(T - t_j)^{(1-\thetaL)/2}\sqrt{t_j-t_i}} 
\end{align*}
Letting $U_i := \Theta_i$, $\Gamma(i)  :=  N^{-1/2}(T-t_i)^{(\alpha\wedge(2\gamma) {\wedge\theta_X}-1)/2}$,  $\Xi(i) := \int_{t_i}^{t_{i+1}} \frac{dr}{(T-r)^{1-\gamma} \sqrt{r-t_i} } $, and 
\begin{equation}
W_i :=   \Gamma(i) + \Xi(i) + \cH(i),
\label{eq:1:W:i}
\end{equation}
it follows from Lemma \ref{lem:iteration:2} that
\begin{align}
\Theta_i & \le C W_i +
C \sum_{j=i+1}^{N-1}\frac{ W_j \Delta_j }{(T - t_j)^{(1-\thetaL)/2}\sqrt{t_j-t_i}} 
+ C \sum_{j=i+1}^{N-1}\frac{ \Theta_j \Delta_j }{(T - t_j)^{(1-\thetaL)/2}} 
\label{eq:1:theta:1}
\end{align}
Therefore, using Lemma \ref{lem:iteration:3} in \eqref{eq:disc er:2} and \eqref{eq:disc er:3},
\begin{align}
\normexp{Z_{t_i} - \bZ Ni} & \le C W_i + C \sum_{j=i+1}^{N-1}\frac{ W_j \Delta_j }{(T - t_j)^{(1-\thetaL)/2}\sqrt{t_j-t_i}}, 
\label{eq:disc er:4}
\\
\normexp{Y_{t_i} - \bY Ni} & \le CN^{-1/2}  + C \sum_{j=i}^{N-1}\frac{ \Psi(j)}{(T - t_j)^{(1-\thetaL)/2}} 
+ C \sum_{j=i}^{N-1}\frac{ W_j \Delta_j }{(T - t_j)^{(1-\thetaL)/2}}.
\label{eq:disc er:5}
\end{align}

Let us consider the sum in the $W$ terms. Firstly, remark that we only need consider the sums for $i < N-1$.
 Recall the terminology of \eqref{eq:1:W:i}.  Using Lemma \ref{lem:integration:1},
\begin{align}
\sum_{j=i+1}^{N-1}  \frac{\Gamma(j) \Delta_j}{(T - t_j)^{(1-\thetaL)/2}\sqrt{t_j-t_i}} & =
C N^{-1/2}\sum_{j=i+1}^{N-1}  \frac{\Delta_j}{(T - t_j)^{1-\thetaL/2}\sqrt{t_j-t_i}} 
%\nonumber \\
%& 
\le CN^{-1/2}(T-t_i)^{(\thetaL-1)/2}.
\label{eq:1:gamma:1}
%\\
%& \nonumber \\
\end{align}
Using the fact that $\Delta_j \le \Delta_{j-1}$ to show that $\sqrt{t_{j+1}-t_i} / \sqrt{t_{j}-t_i} \le 2$,
Lemma \ref{lem:inc growth} to show that $\Delta_j/\Delta_{j+1} \le C$ and $\max_j \Delta_j (T-t_j)^{2\gamma-1} \le N^{-1}$,
one can apply Lemma \ref{lem:integration:1} to bound the sum in $\Xi(j)$ as follows:
\begin{align}
\sum_{j=i+1}^{N-1} & \frac{\Xi(j) \Delta_j}{(T - t_j)^{(1-\thetaL)/2}\sqrt{t_j-t_i}} \nonumber \\
& \le
C\frac{\int_{t_{N-1}}^{T} (T-r)^{\gamma-1}(r-t_{N-1})^{-1/2} dr \Delta_{N-1} }{\Delta_{N-1}^{(1-\thetaL)/2}\sqrt{t_{N-1} - t_i} } 
%\nonumber \\
%& \quad 
+ C\sum_{j=i+1}^{N-2}  \frac{ \Delta_{j+1}^{3/2} (\Delta_j/\Delta_{j+1})^{3/2} \sqrt{t_{j+1}-t_i} / \sqrt{t_{j}-t_i}}
{(T - t_{j+1})^{(3-\thetaL-2\gamma)/2}\sqrt{t_{j+1}-t_i}} 
\nonumber \\
& \le \frac{C\Delta_{N-1}^{(1+\gamma + \thetaL)/2 } \sqrt{\Delta_{N-2} } }{\sqrt{t_{N-1} - t_i} }
 + C\max_j \sqrt{\frac{\Delta_j}{(T-t_j)^{1-2\gamma} }  }\sum_{j=i+1}^{N-2}  \frac{ \Delta_j }
{(T - t_{j+1})^{1-\thetaL/2}\sqrt{t_j-t_i}}
\nonumber \\
& \le CN^{-3/2} + C N^{-1/2} (T-t_i)^{(\thetaL - 1)/2}
\label{eq:1:xi:1} 
%\\
%& \nonumber \\
\end{align}
 In order to deal with the sum in $\cH(j)$, we change the order of summation and apply Lemma \ref{lem:integration:1} to obtain
\begin{align}
\sum_{j=i+1}^{N-1}  \frac{\cH(j) \Delta_j}{(T - t_j)^{(1-\thetaL)/2}\sqrt{t_j-t_i}} & =
\sum_{j=i+1}^{N-1}  \frac{  \sum_{k=j+1}^{N-1}  \frac{\Psi(k) }{(T - t_k)^{(1-\thetaL)/2}\sqrt{t_k-t_j}}
 \Delta_j}{(T - t_j)^{(1-\thetaL)/2}\sqrt{t_j-t_i}} 
\nonumber \\
& = \sum_{k=i+2}^{N-1}  \frac{  \sum_{j=i+1}^{k-1}  \frac{ \Delta_j}{(T - t_j)^{1-\thetaL/2}\sqrt{t_j-t_i}}
  \Psi(k) }{(T - t_k)^{(1-\thetaL)/2}} 
\nonumber\\
& \le C \sum_{j=i+1}^{N-1}  \frac{\Psi(j) }{(T - t_j)^{(1-\thetaL)/2}(t_j-t_i)^{(1-\thetaL)/2} } = C\cH(i).
\label{eq:1:e:2}
\end{align}
%where we have used 
%\[
%\sum_{j=i+1}^{N-2}  \frac{\Delta_j}{(T - t_{j+1})^{1-\thetaL/2}\sqrt{t_j-t_i}} \le C (T-t_i)^{(\thetaL-1)/2},
%\]
%a result which can be proved analogously to Lemma \ref{lem:integration:1}.
Combining \eqref{eq:1:gamma:1} - \eqref{eq:1:e:2}, the bound on the sum in $W_j$ is
\begin{equation}
\sum_{j=i+1}^{N-1}  \frac{W_j \Delta_j}{(T - t_j)^{(1-\thetaL)/2}\sqrt{t_j-t_i}} \le C N^{-1/2} (T-t_i)^{(\thetaL-1)/2}  + C \cH(i).
\label{eq:1:W:2}
\end{equation}
By analogous calculations, one  shows that 
\begin{equation}
\sum_{j=i+1}^{N-1}  \frac{W_j \Delta_j}{(T - t_j)^{(1-\thetaL)/2}} \le C N^{-1/2} (T-t_i)^{(\thetaL+1)/2}  + C \sum_{j=i}^{N-1}\frac{ \Psi(j)}{(T - t_j)^{(1-\thetaL)/2}} .
\label{eq:1:W:3}
\end{equation}
The proof is completed by substituting \eqref{eq:1:W:2} into \eqref{eq:disc er:4}, \eqref{eq:1:W:3} into \eqref{eq:disc er:5}, and using Proposition \ref{prop:Z sums:num er} to bound the remaining terms.
\qed

\appendix 
%\section{Appendix}
\section{Stochastic analysis}
The following {\em conditional} Fubini's theorem is a consequence of the Monotone Class Theorem.
\begin{lemma}
\label{lem:cond:fubini}
Let $f_s \in \L_2([0,T]\times \Omega)$. Then, for all $t\in [0,T]$, there exists a $\cB([0,T])\otimes\cF_t$-measurable processes $F_t$
% measurable function 
% \[
% F_t : ([0,T] \times \Omega , \cB([0,T])\otimes\cF_t) \rightarrow (\bar\R,\cB(\bar\R))
% \]
belonging to $L_2([0,T]\times\Omega)$ such that $(\omega,s)\mapsto F_t(s)$ is a version of $(\omega,s)\mapsto \E_t[f_s]$  and
\[
\E_t[\int_0^T f_s ds ] = \int_0^T F_t(\cdot,s) ds \quad \text{almost surely.}
\]
\end{lemma}

\vspace{0.3cm}
We need the following generalization of the a priori estimates \cite[Proposition 3.2]{bria:dely:hu:pard:stoi:03}:
\begin{proposition}
\label{prop:briand:et:al}
Let $k$ be an integer, and $p$ be an integer greater than or equal to $2$.
Let $f : \Omega\times[0,T)\times (\R^k)^\top \times \R^{ q\times k} \rightarrow (\R^k)^\top$  be 
$\cP \times \cB\big((\R^k)^\top \big) \otimes \cB(\R^{q \times k})$-measurable,
and $\xi$ be an $(\R^k)^\top$-valued random variable in $\L_p(\cF_T)$.
Let  $(f_t)_{t\in [0,T]}$ be  non-negative,  predictable process,
$\mu \in \L_1 ([0,T];m)$ and $\lambda \in \L_2([0,T];m)$ be $\R$-valued non-negative.
% , and assume that 
% $f_t$, $\mu_t$ and $\lambda_t$ take finite values for all $t$ in $[0,T)$.
Additionally, assume that $\e[(\int_0^T f_t dt)^p] <\infty$.
For any $(y_1,y_2) \in (\R^k)^2$, define the scalar product $(y_1,y_2) := \sum_{j=1}^k y_{1,j}y_{2,j}$ and
assume that, for all $(t,y,z)\in [0,T) \times(\R^k)^\top\times\R^{k \times q}$, $(\omega,t,y,z) \mapsto f(\omega, t,y,z)$ 
satisfies
\begin{equation}
\big( |y|^{-1}y \1_{|y|>0} , f(\omega, t,y,z) \big)\le f_t(\omega) + \mu_t |y| + \lambda_t |z| \quad  \text{almost surely.}
\label{eq:bria et al} 
\end{equation}
Let $(Y,Z)$ be a solution to the $\big( (\R^k)^\top, \R^{q\times k}\big)$-valued BSDE 
\[
Y_t = \xi + \int_t^T f(r,Y_r,Z_r) dr - \sum_{j=1}^q\int_t^T (Z_{j,r})^\top dW_{j,r}.
\]
in the space $\cS^p \times \cH^p$, where $\cH^p$ is the space of 
predictable processes $X$ such that $\E[ ( \int_0^T |X_s|^2 ds) ^{p/2}]$ is finite; $Z_j$ denotes the $j$-th column of $Z$.

Then, there exists a constant $C_p$, depending only on $p$, such that, for any $\eta_t \ge \mu_t + \lambda_t^2/(p-1) $ in $L_1(\R;dt)$,
\[ \e[\sup_t e^{p \int_0^t \eta_r dr}|Y_t|^p + (\int_0^T e^{2\int_0^t\eta_rdr}|Z_t|^2 dt)^{p/2}]
\le C_p \e[  e^{ p \int_0^T\eta_rdr} |\xi|^p + (\int_0^T e^{\int_0^t\eta_rdr}f_tdt)^p]. \]
\end{proposition}
{\bf Proof.}
Consider the processes $\tilde Y_t =  e^{\int_0^t\eta_rdr} Y_t$ and $\tilde Z_t =  e^{\int_0^t\eta_rdr} Z_t$.
Then $(\tilde Y, \tilde Z)$ satisfies a BSDE with terminal condition $\tilde \xi = e^{\int_0^T\eta_rdr} \xi$
and driver $\tilde f(t,y,z) = e^{\int_0^t\eta_rdr} f(t, e^{-\int_0^t\eta_rdr}y, e^{-\int_0^t\eta_rdr}z) - \eta_t y$.
Moreover, for all $(t,y,z) \in[0,T) \times \R^k \times \R^{k\times q}$, $\tilde f(\omega,y,z)$ satisfies 
\[
 \big( |y|^{-1}y \1_{|y|>0} , \tilde f(\omega, t,y,z) \big)\le \tilde f_t(\omega) + \tilde \mu_t |y| + \tilde \lambda_t |z|
\quad  \text{almost surely.}
\]
with $\tilde f_t = e^{-\int_0^t\eta_rdr} f_t$, $\tilde \mu_t = \mu_t - \eta_t$, and $\tilde \lambda_t = \lambda_t$.
The rest of the proof follows exactly as the proof of \cite[Proposition 3.2]{bria:dely:hu:pard:stoi:03}. \qed

\section{Time-grids}
\begin{lemma}
\label{lem:inc growth}
The time grid  $\pib = \{0 = t_0 <\ldots < t_N =T \; : \; t_i=T-T(1-i/N)^{1/{\beta}}\}$ 
 with $\beta \in (0,1]$ satisfies 
 \begin{align}
&\max_{0\le i<N}\frac{\Delta_k }{(T-t_k)^{1-\theta} }
 \leq \frac{ T^{\theta} }{\beta }\frac{ 1 }{N^{1\land \frac{\theta  }{\beta }} }, 
\label{eq:grid bd:1}\\
 & \max_{0 \le i \le N-2} \frac{\Delta_k}{\Delta_{k+1}} \le \frac{1}{\beta}\bigg( 1 \vee \big(\frac{ 1 }{2\beta }\big)^{\frac{1  }{\beta}-1} \bigg),
\label{eq:grid bd:2}
 \end{align}
for all  $\theta\in(0,1]$.
\end{lemma}

\section{Integral estimates}
The following is a trivial result that will come in useful.
\begin{lemma}
\label{lem:simple:int}
For finite $\delta > 0$, and $s<r \le R<\infty$, $\int_s^r (R -t)^{\delta -1} dt \le \frac1\delta (r-s)(T-s)^{\delta-1}$.
\end{lemma}
{\bf Proof.}
Direct computation of the integral term yields
\begin{align*}
\int_{s}^{r} \frac{dt }{(R-t)^{1-\delta} } & = \frac1\delta \big\{ (R-s)^{\delta} - (R - r)^{\delta}\big\} 
%\\
%&
 \le \frac1\delta\Big\{ \frac{R - s}{(R - s)^{1-\delta}} - \frac{R-r}{(R-s)^{1-\delta}} \Big\}
 = \frac{r-s}{\delta (R-s)^{1 - \delta}} .
\end{align*}

\qed

The following three lemmas and their proofs can be found in Section 2.1 of \cite{gobe:turk:13b};
the results on the integrals are proved exactly as the results on the sums.
\begin{lemma}
\label{lem:integration:1}
Let $\delta,\rho \in(0,1]$. Then for $B_{\delta,\rho} := \int_0^1 (1-r)^{\delta-1}r^{\rho-1} dr$, 
for any $0\le s < t\le T$,
 {\[\int_{t}^{s} {dr \over (s - r)^{1-\delta }(r-s)^{1-\rho } } \le B_{\delta,\rho}(s - t)^{\delta + \rho -1}.\]}
Moreover, on the time-grid $\pib = \{0 = t_0 <\ldots < t_N =T \; : \; t_i=T-T(1-i/N)^{1/{\beta}}\}$,
for any $0 \le i < k \le N$,
\[
\sum_{j=i+1}^{k-1}(t_k - t_j)^{\delta -1}(t_j-t_i)^{\rho -1}\Delta_j \le 2 B_{\delta,\rho}(t_k - t_i)^{\delta + \rho -1}.
\]
\end{lemma} 
%The proof of the integral part of Lemma \ref{lem:integration:1} is done by changing variables. The proof of the sum part of
%Lemma \ref{lem:integration:1} is exactly the same as the proof of Lemma \ref{lem:integrals} in Section \ref{proof:lem:integrals} later.
%We note that, in the context of the notation of that Paper, $R_\pi$ is equal to 1, because the time-grid $\pib$ has the property
%$\Delta_{i+1} \le \Delta_i$ for all $i$.

\begin{lemma}
\label{lem:iteration:2}
Let $\delta \in (0,1/2]$, $ \rho > 0$ and $t\in[0,T)$.
Suppose that, for a positive constant $\cuM$, the finite positive real functions $u : [t,T] \mapsto [0,\infty)$ 
and $w : [t,T] \mapsto [0,\infty)$ satisfy
\begin{align}
u_t \le  w_t +\cuM \int_{t}^{T} \frac{u_r dr}{(T-r)^{{\frac 1 2 -  \delta}}(r-t)^{{\frac 1 2-\rho}}}.
\label{eq:iteration:1}
\end{align}
Then, for  constants $\cwM$ and $\cuhM$  depending only on $\cuM, T, \delta$ and $\rho$, 
\begin{align}
u_t \le \cwM w_t +  \cwM \int_{t}^{T} \frac{w_r dr}{(T-r)^{{\frac 1 2 -  \delta}}(r - t)^{{\frac 1 2-\rho}}}
+\cuhM \int_t^{T} \frac{u_r dr}{(T-r)^{{\frac 1 2 -  \delta}}}.
\label{eq:iteration:2}
\end{align}
Moreover, on the time-grid $\pib = \{0 = t_0 <\ldots < t_N =T \; : \; t_i=T-T(1-i/N)^{1/{\beta}}\}$,
suppose that the real functions $U : \pib \mapsto [0,\infty)$ and $W: \pib \mapsto [0,\infty)$ satify
\begin{align}
U_i \le  W_i +\cuM \sum_{j=i+1}^{N-1} \frac{U_j\Delta_j}{(T-t_j)^{{\frac 1 2 -  \delta}}(t_j-t_i)^{{\frac 1 2-\rho}}}
\label{eq:iteration:3}
\end{align}
for all $i\in\{0,\ldots,N-1\}$. It follows that 
\begin{align*}
U_i \le 2\cwM W_i +  2\cwM \sum_{j=i+1}^{N-1} \frac{W_j\Delta_j}{(T-t_j)^{{\frac 1 2 -  \delta}}(t_j - t_i)^{{\frac 1 2-\rho}}}
+2 \cuhM \sum_{j=i+1}^{N-1} \frac{U_j \Delta_j}{(T-t_j)^{{\frac 1 2 -  \delta}}}
\end{align*}
for all $i\in\{0,\ldots,N-1\}$.
\end{lemma}
%The proof of Lemma \ref{lem:iteration:2} is analogous to the proof of Lemma \ref{lem:iteration:gen:1} in Paper 4. 
%The integral part is proved in the same way as the sum part, with, of course, integrals replacing the sums.
\begin{lemma}
\label{lem:iteration:3}
Let $\delta \in(0,1/2]$, $ \rho > 0$ and $t\in[0,T)$.
Suppose that the finite positive real functions $u : [t,T] \mapsto [0,\infty)$ 
and $w : [t,T] \mapsto [0,\infty)$ satisfy
\eqref{eq:iteration:2} for some positive constants $\cwM$ and $\cuhM$.
Then, for $\nu> 0$, there is a positive constant $\cg{\nu}$ (depending only on $\cwM, \cuhM, T, \delta, \rho, \nu$) such that
\begin{align}
\int_t^T \frac{u_r dr}{(T-r)^{\frac12 - \delta}(r -t)^{1-\nu} }
& 
\le \cg\nu \int_t^T \frac{w_r dr}{(T-r)^{\frac12 - \delta}(r -t)^{1-\nu} }
\label{eq:integration:3}
\end{align}
Moreover, on the time-grid $\pib = \{0 = t_0 <\ldots < t_N =T \; : \; t_i=T-T(1-i/N)^{1/{\beta}}\}$,
suppose that the real functions $U : \pib \mapsto [0,\infty)$ and $W: \pib \mapsto [0,\infty)$ satify \eqref{eq:iteration:3} for all
$i\in\{0,\ldots,N-1\}$. It follows that
\begin{align*}
\sum_{j=i+}^{N-1} \frac{U_j\Delta_j}{(T-t_j)^{\frac12 - \delta}(t_j -t_i)^{1-\nu} }
& 
\le 2\cg\nu \sum_{j=i+1}^{N-1} \frac{W_j \Delta_j}{(T-t_j)^{\frac12 - \delta}(t_j -t_i)^{1-\nu} }
\end{align*}
for all $i\in\{0,\ldots,N-1\}$.
\end{lemma}
%The proof of Lemma \ref{lem:iteration:3} is analogous to the proof of Lemma \ref{lem:iteration:gen:2} in Paper 4. 
%The integral part is proved in the same way as the sum part, with, of course, integrals replacing the sums.

\section{Regularity results for inverse matrices}
\label{section:mx}
\begin{lemma}
\label{lem:1:mx inverse}
Let $\xi >0$ be finite and  $A : \R^n \rightarrow \R^{l\times l}$ be symmetric and such that $\eta ^\top A(x) \eta \ge \xi |\eta|^2$
for all $x \in  \R^n$ and $\eta \in \R^l$.
Then, for every $x\in\R^d$, the matrix $A(x)$ is invertible and $|A^{-1}(x) | \le 1/\xi$.
Moreover, if $x\mapsto A(x)$ is $\gamma$-H\"older continuous,% with respect to the matrix 2-norm, 
then it's inverse $x \mapsto A^{-1}(x)$ 
is also $\gamma$-H\"older continuous.% with respect to the matrix 2-norm.
\end{lemma}
{\bf Proof.}
Due to the condition $\eta ^\top A(x) \eta \ge \xi |\eta|^2$, it follows that $A(x)$ is positive definite for every $x\in\R^n$. 
This implies that the singular values of $A(x)$ are all greater than $\xi$ \cite[Theorem 8.1.2]{golu:vanl:96},  and so $A(x)$ is invertible.
Using the singular value decomposition of $A(x)$ to construct the inverse as in \cite[Section 5.5.4]{golu:vanl:96}, 
the maximal sigular value of $A^{-1}(x)$ is less than $1/\xi$ and so, using \cite[Section 2.5.2]{golu:vanl:96}  combined with the 
singular value decomposition of $A^{-1}(x)$, the matrix 2-norm of $A^{-1}(x)$ is equal to its maximal singular value, i.e.
$|A^{-1}(x)| \le 1/\xi$ for all $x \in \R^d$.
Now, let $x$ and $y$ be elements in $\R^d$. 
Since $A^{-1}(y) - A^{-1}(x)$ is equal to  \[ - A(x)^{-1} \big( A(y) - A(x) \big) A(y)^{-1},\]
it follows that 
\[
|A^{-1}(y) - A^{-1}(x)| \le |A^{-1}(x)| | A(y) - A(x) | |A^{-1}(y)| \le \frac {L_A}{\xi^2} |x-y|^\gamma,
\]
where $L_A$ is the H\"older constant of $A$. \qed

\vspace{0.3cm}
 {
{\bf Proof of Lemma \ref{lem:1:inv sigma:lip}.}
Let $t\in[0,T)$ be fixed,
and define $A : [0,T) \times \R^d \rightarrow \R^{d\times d}$ by $A(t,x) = \sigma(t,x) \sigma(t,x)^\top$.
It can be computed directly that $\sigma^{-1}(\cdot) = \sigma(\cdot)^\top A^{-1}(\cdot)$, whether  or not $d$ equals $q$.
It follows from uniform ellipticity \He \ and Lemma \ref{lem:1:mx inverse} that  $|A^{-1}(t,x)| \le 1/\elip$ for all $(t,x) \in [0,T)\times\R^d$.
Due to the differentiability condition \Hs \ on $\sigma(t,\cdot)$, $\sigma(t,\cdot)$ is Lipschitz continuous uniformly in $t$ 
with Lipschitz constant $\|\nabla_x \sigma \|_\infty$, and, using additionally the equality 
$A(t,x) - A(t,y) = \sigma(y)(\sigma(t,x)^\top - \sigma(t,y)^\top) + (\sigma(t,x) - \sigma(t,y)) \sigma(t,x)^\top$,
$A(t,\cdot)$ is Lipschitz continuous uniformly in $t$ with Lipschitz constant
$2\|\sigma\|_\infty\|\nabla_x\sigma\|_\infty$.
Using Lemma \ref{lem:1:mx inverse}, it follows that $A^{-1}(t,\cdot)$ is Lipschitz continuous uniformly in $t$ with Lipschitz constant
$2\|\sigma\|_\infty\|\nabla_x\sigma\|_\infty/\elip^2$.
For any $(x,y)\in(\R^d)^2$, $\sigma(t,x)^{-1} - \sigma(t,y)^{-1}$ is equal to 
$\big(\sigma(t,x)^\top - \sigma(t,y)^\top \big) A^{-1}(t,x) + \sigma(t,y)^\top \big( A^{-1}(t,x) - A^{-1}(t,y) \big)$.
 and therefore
\[
|\sigma(t,x)^{-1} - \sigma(t,y)^{-1}| \le \frac{\|\nabla_x \sigma \|_\infty}{\elip} |x-y| 
+ \frac{2\|\sigma\|_\infty \|\nabla_x\sigma \|_\infty}{\elip^2} |x-y|.
\]
The proof that $\sigma^{-1}(\cdot,x)$ is $1/2$-H\"older continuous is essentially the same and we do not include it.
\\
 \qed
 }

%\section{References}
\bibliographystyle{alpha}
%\bibliography{../../Firstdraft/discretization}

\begin{thebibliography}{BDH{\etalchar{+}}03}

\bibitem[BDH{\etalchar{+}}03]{bria:dely:hu:pard:stoi:03}
P.~Briand, B.~Delyon, Y.~Hu, E.~Pardoux, and L.~Stoica.
\newblock {$\L_p$} solutions of backward stochastic differential equations.
\newblock {\em Stochastic Processes and their Applications}, 108(1):109--129,
  2003.

\bibitem[BL13]{bria:laba:13}
P.~Briand and C.~Labart.
\newblock Simulation of {BSDEs} by {W}iener {C}haos {E}xpansion.
\newblock {\em To appear in Annals of Applied Probability}, 2013.

\bibitem[BT04]{bouc:touz:04}
B.~Bouchard and N.~Touzi.
\newblock {Discrete time approximation and {M}onte {C}arlo simulation of
  backward stochastic differential equations}.
\newblock {\em Stochastic Processes and their Applications}, 111:175--206,
  2004.

\bibitem[CD12]{cris:dela:12}
D.~Crisan and F.~Delarue.
\newblock Sharp derivative bounds for solutions of degenerate semi-linear
  partial differential equations.
\newblock {\em Journal of Functional Analysis}, 263(10):3024--3101, 2012.

\bibitem[CR14]{chas:rich:14}
J.~F. Chassagneux and A.~Richou.
\newblock Numerical simulation of quadratic {BSDEs}.
\newblock {\em Available on http://arxiv.org/abs/1307.5741}, 2014.

\bibitem[DG06]{dela:guat:06}
F.~Delarue and G.~Guatteri.
\newblock Weak existence and uniqueness for forward-backward {SDE}s.
\newblock {\em Stochastic Processes and their Applications},
  116(12):1712--1742, 2006.

\bibitem[EKPQ97]{elka:peng:quen:97}
N.~El~Karoui, S.~Peng, and M.~C. Quenez.
\newblock Backward stochastic differential equations in finance.
\newblock {\em Mathematical Finance. An International Journal of Mathematics,
  Statistics and Financial Economics}, 7(1):1--71, 1997.

\bibitem[FJ12]{fan:jian:12}
S.~J. Fan and L.~Jiang.
\newblock {$\L_p$} solutions of finite and infinite time interval {BSDE}s with
  non-{L}ipschitz coefficients.
\newblock {\em Stochastics}, 84(4):487--506, 2012.

\bibitem[Fri64]{frie:64}
A.~Friedman.
\newblock {\em Partial differential equations of parabolic type}.
\newblock Prentice-Hall Inc., Englewood Cliffs, N.J., 1964.

\bibitem[GGG12]{geis:geis:gobe:12}
C.~Geiss, S.~Geiss, and E.~Gobet.
\newblock Generalized fractional smoothness and {$\L_p$}-variation of {BSDE}s
  with non-{L}ipschitz terminal condition.
\newblock {\em Stochastic Processes and their Applications}, 122(5):2078--2116,
  2012.

\bibitem[GL06]{gobe:lemo:06}
E.~Gobet and J.~P. Lemor.
\newblock Numerical simulation of {BSDE}s using empirical regression methods:
  theory and practice.
\newblock In {\em Proceedings of the Fifth Colloquium on BSDEs (29th May - 1st
  June 2005, Shangai) - Available on
  http://hal.archives-ouvertes.fr/hal-00291199/fr/}, 2006.

\bibitem[GL07]{gobe:laba:07}
E.~Gobet and C.~Labart.
\newblock Error expansion for the discretization of backward stochastic
  differential equations.
\newblock {\em Stochastic Processes and their Applications}, 117(7):803--829,
  2007.

\bibitem[GL10]{gobe:laba:10}
E~Gobet and C.~Labart.
\newblock Solving {BSDE} with adaptive control variate.
\newblock {\em SIAM Journal on Numerical Analysis}, 48(1):257--277, 2010.

\bibitem[GM{\etalchar{+}}05]{gobe:muno:05}
E.~Gobet, R.~Munos, et~al.
\newblock Sensitivity analysis using {It\^o}-malliavin calculus and
  martingales, and application to stochastic optimal control.
\newblock {\em SIAM Journal on control and optimization}, 43(5):1676--1713,
  2005.

\bibitem[GM10]{gobe:makh:10}
E.~Gobet and A.~Makhlouf.
\newblock {${\bf L}_2$}-time regularity of {BSDE}s with irregular terminal
  functions.
\newblock {\em Stochastic Processes and their Applications}, 120(7):1105--1132,
  2010.

\bibitem[GT13a]{gobe:turk:13b}
E.~Gobet and P.~Turkedjiev.
\newblock Approximation of backward stochastic differential equations using
  malliavin weights and least-squares regression.
\newblock {\em To appear in Bernoulli, available on
  http://hal.archives-ouvertes.fr/hal-00855760}, 2013.

\bibitem[GT13b]{gobe:turk:13a}
E.~Gobet and P.~Turkedjiev.
\newblock Linear regression {MDP} scheme for discrete backward stochastic
  differential equations under general conditions.
\newblock {\em Available on http://hal.archives-ouvertes.fr/hal-00642685},
  2013.

\bibitem[GVL96]{golu:vanl:96}
G.~H. Golub and C.~F. Van~Loan.
\newblock {\em Matrix computations}.
\newblock Johns Hopkins Studies in the Mathematical Sciences. Johns Hopkins
  University Press, Baltimore, MD, third edition, 1996.

\bibitem[HIM05]{hu:imke:mull:05}
Y.~Hu, P.~Imkeller, and M.~M{\"u}ller.
\newblock Utility maximization in incomplete markets.
\newblock {\em The Annals of Applied Probability}, 15(3):1691--1712, 2005.

\bibitem[HNS11]{hu:nual:song:11}
Y.~Hu, D.~Nualart, and X.~Song.
\newblock {Malliavin calculus for backward stochastic differential equations
  and application to numerical solutions}.
\newblock {\em Ann. Appl. Probab.}, 21(6):2379--2423, 2011.

\bibitem[IDR10]{imke:dosr:10}
P.~Imkeller and G.~Dos~Reis.
\newblock Path regularity and explicit convergence rate for {BSDE} with
  truncated quadratic growth.
\newblock {\em Stochastic Processes and their Applications}, 120(3):348--379,
  2010.

\bibitem[JS03]{jaco:shir:03}
J.~Jacod and A.~N. Shiryaev.
\newblock {\em Limit theorems for stochastic processes}, volume 288 of {\em
  Grundlehren der Mathematischen Wissenschaften [Fundamental Principles of
  Mathematical Sciences]}.
\newblock Springer-Verlag, Berlin, second edition, 2003.

\bibitem[Kus03]{kusu:03}
S.~Kusuoka.
\newblock Malliavin calculus revisited.
\newblock {\em The University of Tokyo. Journal of Mathematical Sciences},
  10(2):261--277, 2003.

\bibitem[MZ02]{ma:zhan:02}
J.~Ma and J.~Zhang.
\newblock Representation theorems for backward stochastic differential
  equations.
\newblock {\em The Annals of Applied Probability}, 12(4):1390--1418, 2002.

\bibitem[Nee11]{nee:11}
C.~Nee.
\newblock {\em Sharp Gradient Bounds for the Diffusion Semigroup}.
\newblock PhD thesis, Imperial College London, 2011.

\bibitem[Nua06]{nual:95}
D.~Nualart.
\newblock {\em The {M}alliavin calculus and related topics}.
\newblock Probability and its Applications (New York). Springer-Verlag, Berlin,
  second edition, 2006.

\bibitem[REK00]{elka:roug:01}
R.~Rouge and N.~El~Karoui.
\newblock Pricing via utility maximization and entropy.
\newblock {\em Mathematical Finance. An International Journal of Mathematics,
  Statistics and Financial Economics}, 10(2):259--276, 2000.
\newblock INFORMS Applied Probability Conference (Ulm, 1999).

\bibitem[Ric11]{rich:11}
A.~Richou.
\newblock Numerical simulation of {BSDEs} with drivers of quadratic growth.
\newblock {\em The Annals of Applied Probability}, 21(5):1933--1964, 2011.

\bibitem[Ric12]{rich:12}
A.~Richou.
\newblock Markovian quadratic and superquadratic {BSDE}s with an unbounded
  terminal condition.
\newblock {\em Stochastic Processes and their Applications}, 122(9):3173--3208,
  2012.

\bibitem[RY99]{revu:yor:01}
D.~Revuz and M.~Yor.
\newblock {\em Continuous martingales and {B}rownian motion}, volume 293 of
  {\em Grundlehren der Mathematischen Wissenschaften [Fundamental Principles of
  Mathematical Sciences]}.
\newblock Springer-Verlag, Berlin, third edition, 1999.

\bibitem[Zha04]{zhan:04}
J.~Zhang.
\newblock A numerical scheme for {BSDE}s.
\newblock {\em The Annals of Applied Probability}, 14(1):459--488, 2004.

\bibitem[Zha05]{zhan:05}
J.~Zhang.
\newblock Representation of solutions to {BSDE}s associated with a degenerate
  {FSDE}.
\newblock {\em The Annals of Applied Probability}, 15(3):1798--1831, 2005.

\end{thebibliography}

\newcommand{\etalchar}[1]{$^{#1}$}

\end{document}